\documentclass[12pt,twoside]{article}
\usepackage[margin=0.9in]{geometry}
\usepackage{fancyhdr}
\usepackage{graphicx}
\usepackage{amssymb}
\usepackage{amsmath}
\usepackage{amsthm}
\usepackage{amsfonts}
\usepackage{mathrsfs}
\usepackage{mathtools}
\usepackage{bm}
\usepackage{color}
\usepackage{setspace}
\usepackage{exscale}
\usepackage{relsize}
\usepackage{float}
\usepackage{cite}
\usepackage{hyperref}
\usepackage{cleveref}
\usepackage{nicefrac}
\usepackage{ulem}
\usepackage{caption}
\usepackage{subcaption}
\usepackage{comment}
\usepackage{booktabs,multirow} % for much better looking tables
\usepackage{array} % for better arrays (eg matrices) in maths
\usepackage{paralist} % very flexible & customisable lists (eg. enumerate/itemize, etc.)
\usepackage{fancyhdr} % This should be set AFTER setting up the page geometry

\theoremstyle{plain}                    % use "default" font
\newtheorem{thm}{Theorem}[section]

\newtheorem{rmk}[thm]{Remark}
\newenvironment{acknowledgment}{{\flushleft \bf Acknowledgment:}}{}

\usepackage{tikz}
\usetikzlibrary{positioning}
\usepackage{xcolor}
\allowdisplaybreaks[1]
\numberwithin{equation}{section}
\numberwithin{figure}{section}
\numberwithin{table}{section}
\newcommand\eref[1]{(\ref{#1})}
\newcommand\fref[1]{Figure~\ref{#1}}

\newcommand*\xbar[1]{%
  \hbox{%
    \vbox{%
      \hrule height 0.5pt % The actual bar
      \kern0.4ex%         % Distance between bar and symbol
      \hbox{%
        \kern-0.05em%      % Shortening on the left side
        \ensuremath{#1}%
        \kern-0.00em%      % Shortening on the right side
      }%
    }%
  }%
}

\newcommand{\hf}{{\frac{1}{2}}}

\newcommand{\lph}{{\ell+\frac{1}{2}}}

\graphicspath{{Figures/}}

\hypersetup{colorlinks=true,
            breaklinks=true,
            urlcolor=blue,
            linkcolor=black,
            bookmarksopen=false,
            filecolor=black,
            citecolor=black,
            linkbordercolor=blue
}

\title{Uncertainty Quantification in Forward Problems: Balancing Accuracy and Robustness Using CWENO Interpolations}
\author{Alina Chertock\thanks{Department of Mathematics, NC State University, Raleigh, NC, USA;
{\href{mailto:chertock@math.ncsu.edu}{chertock@math.ncsu.edu}}}, Arsen S. Iskhakov\thanks{Department of Mechanical and Nuclear Engineering, 
Kansas State University, Manhattan, KS, USA; {\href{mailto:aiskhak@ksu.edu}{aiskhak@ksu.edu}}}, and
Alexander Kurganov\thanks{Department of Mathematics and Shenzhen International Center for Mathematics, SUSTech, Shenzhen 518055, China;
{\href{mailto:alexander@sustech.edu.cn}{alexander@sustech.edu.cn}}}}
\date{}

\begin{document}
\maketitle
\begin{abstract}
In this paper, we study uncertainty quantification (UQ) in forward problems. Our objective is to construct accurate and robust surrogate
models by incorporating the seventh-order central weighted essentially non-oscillatory (CWENO7) scheme into the stochastic collocation
framework. A key focus is on mitigating the oscillatory behavior often encountered in traditional spectral methods while retaining
high-order accuracy in smooth regions.

We present a systematic comparison between CWENO7-based and generalized polynomial chaos (gPC)-based approaches. Although gPC methods
achieve spectral convergence, they are prone to Gibbs-type oscillations in nonsmooth settings. By contrast, CWENO7 utilizes local stencils
to achieve a balance: non-oscillatory behavior near discontinuities and high-order convergence in smooth regions.

To validate the approach, we conduct numerical experiments on a range of one- and two-dimensional smooth and nonsmooth problems, including 
shallow water equations with random inputs. The results demonstrate that CWENO7 interpolation provides accurate estimates of probability
density functions, mean values, and standard deviations, particularly in regimes where gPC expansions exhibit strong oscillations.
Furthermore, computational tests confirm that CWENO7 interpolation is efficient and scalable, establishing it as a reliable alternative to
conventional stochastic collocation techniques for UQ in the presence of discontinuities.
\end{abstract}



%%%%%%%%%%%%%
% SECTION 1 %
%%%%%%%%%%%%%
\section{Introduction}
Many scientific and engineering problems involve inherent uncertainties arising from various sources. Quantifying these uncertainties is
crucial for assessing predictive capability and improving the accuracy of the numerical models; see, e.g., \cite{oberkampf2010}. 
Uncertainty quantification (UQ) in forward problems refers to the process of propagating input uncertainties through a mathematical model 
to quantify the resulting uncertainties in the output predictions. In general, forward models can be expressed as
\begin{equation}
{\cal M}(\bm U;\bm\xi)=\bm0,
\label{1.1}
\end{equation}
where the input uncertainties are represented by (real-valued) random variables $\bm\xi\in\Xi\subset\mathbb R^s$, defined on a probability
space $(\Xi,{\cal F},p)$ with a $\sigma$-algebra ${\cal F}$ and probability density function (PDF) $p(\bm\xi)$, and the output function
$\bm U(\bm x,t;\bm\xi)\in\mathbb R^K$ denotes the model prediction, which may depend on spatial ($\bm x\in\mathbb R^d$) and temporal
($t\in\mathbb R$) variables. The objective of UQ is to estimate statistical properties of $\bm U$, such as means, variances, or even full
PDFs, while accounting for variability in $\bm\xi$.

Various numerical techniques have been developed for UQ in forward problems, including Monte Carlo (MC) methods and stochastic collocation
with generalized polynomial chaos (gPC) expansions being a common approach within the latter framework. While MC simulations are
straightforward and non-intrusive, they are often computationally expensive due to the large number of realizations (samples) required for
convergence; see, e.g., \cite{AM17,MSS13,MSS12a}.

Stochastic collocation methods, on the other hand, evaluate the deterministic model ${\cal M}(\bm U;\bm\xi)$ at a set of collocation
points, which are typically determined by specific quadrature rules for gPC, and construct a surrogate model that can efficiently
approximate $\bm U$ at additional points. This enables the computation of statistical moments and PDF with reduced computational effort;
see, e.g., \cite{Xiu09,Xiu07,Xiu10}. The gPC approach offers spectral convergence for smooth output functions (see, e.g.,
\cite{SMS13,PIN-book}), but can exhibit Gibbs-type oscillations when the output function is discontinuous or nonsmooth (see, e.g.,
\cite{LKNG04,WanK-SISC06}). Spline-based approaches can mitigate oscillations, but may oversmear sharp features in $\bm U$ (see, e.g.,
\cite{DFS19,CIJK,chertock2025}), prompting the need for alternative approximation strategies.

To address these limitations, several techniques have been proposed within the stochastic collocation framework. They are based on either
detection of the discontinuities in the stochastic space \cite{MZ09,SNR14,wan2005} or tracking them with the help of either the level set
\cite{pettersson2019} or machine learning \cite{JCN19,TB18} methods. After the ``rough'' parts of $\bm U$ are identified in the stochastic
space, one may apply either the localized gPC expansions (leading to the multi-element gPC approach \cite{wan2005}) or piecewise surrogate
models (see, e.g., \cite{JCN19,MZ09,pettersson2019,SNR14,TB18}) there. While capable of reducing Gibbs-type oscillations in regions with
sharp gradients or discontinuities, the aforementioned methods increase computational complexity due to the need for multiple local
approximations, careful handling/tracking the interfaces between the ``rough'' and smooth parts of the approximant, and sensitivity to
additionally introduced parameters.

In this work, we use uniformly high-order accurate central weighted essentially non-oscillatory (CWENO) interpolations for constructing
surrogate models. CWENO approximations were developed in the context of finite-volume methods for hyperbolic conservation laws (see, e.g.,
\cite{CPSV,LPR99,STP,Zah09} and references therein), where they were used to obtain uniformly high-order accurate reconstructions out of
available cell averages of the computed solution. For CWENO interpolations, which are based on the set of given point values of $\bm U$, we
refer the reader to \cite{CPSV,DZP}. These interpolations achieve high-order accuracy in smooth regions while effectively suppressing
oscillations near discontinuities. Moreover, unlike the gPC expansions, which have to be constructed using a specific set of Gaussian
quadrature nodes, CWENO interpolations can be used on any set of nodes in the stochastic variables without producing any oscillations near
the end-point of the computational domain (Runge-type phenomenon), which will appear if the collocation points are arbitrarily selected for
the gPC expansion.

The main goal of this study is threefold: (i) to develop a surrogate modeling approach based on the seventh-order CWENO (CWENO7)
interpolation; (ii) to compare the performance of the CWENO7- and gPC-based stochastic collocation methods; to demonstrate the computational
efficiency and accuracy of the CWENO7 approach on examples involving both smooth and discontinuous output functions $\bm U$.

The remainder of the paper is organized as follows. In \S\ref{sec2}, we present a detailed mathematical formulation of the
proposed CWENO7-based surrogate model. Numerical experiments and comparative analyses are provided in \S\ref{sec3}, highlighting the
advantages and limitations of the CWENO7-based method and its gPC-based counterpart; the latter method is described in Appendix \ref{appA}.
Finally, \S\ref{sec4} summarizes the key findings and outlines potential future directions.

%%%%%%%%%%%%%
% SECTION 2 %
%%%%%%%%%%%%%
\section{Methodology}\label{sec2}
\subsection{One Random Variable ($s=1$)}
We begin by selecting a set of collocation points $\{\xi_\ell\}_{\ell=1}^L$ in the random space, where $L$ is usually limited by the
computational cost of simulating the model \eref{1.1} to obtain the corresponding model outputs $\{\bm U(\bm x,t;\xi_\ell)\}_{\ell=1}^L$. A
new surrogate model is then constructed using the CWENO7 interpolation of these output functions in the random space. The new model will
allow one to efficiently and accurately estimate statistical moments such as the mean, variance, and standard deviation for each component
$U$ of $\bm U$,
\begin{equation}
\mu[U]:=\int\limits_\Xi U(\bm x,t;\xi)p(\xi)\,{\rm d}\xi,\quad{\rm Var}[U]:=\mu[U^2]-(\mu[U])^2,\quad\sigma[U]:=\sqrt{{\rm Var}[U]},
\label{2.1}
\end{equation}
as well as the corresponding PDFs. While evaluating the quantities in \eref{2.1} requires an accurate computation of the integral in the
formula for $\mu[U]$, the PDF reconstruction of $U$ can be performed using the histogram method. In the numerical examples presented in
\S\ref{sec3}, we have used \texttt{numpy.histogram} function in \texttt{Python} with the \texttt{auto} binning strategy.

\subsubsection{CWENO7 Interpolation}
\begin{figure}[ht!]
\centerline{\includegraphics[trim=0.0cm 0.0cm 0.0cm 0.0cm, clip, width=0.35\textwidth]{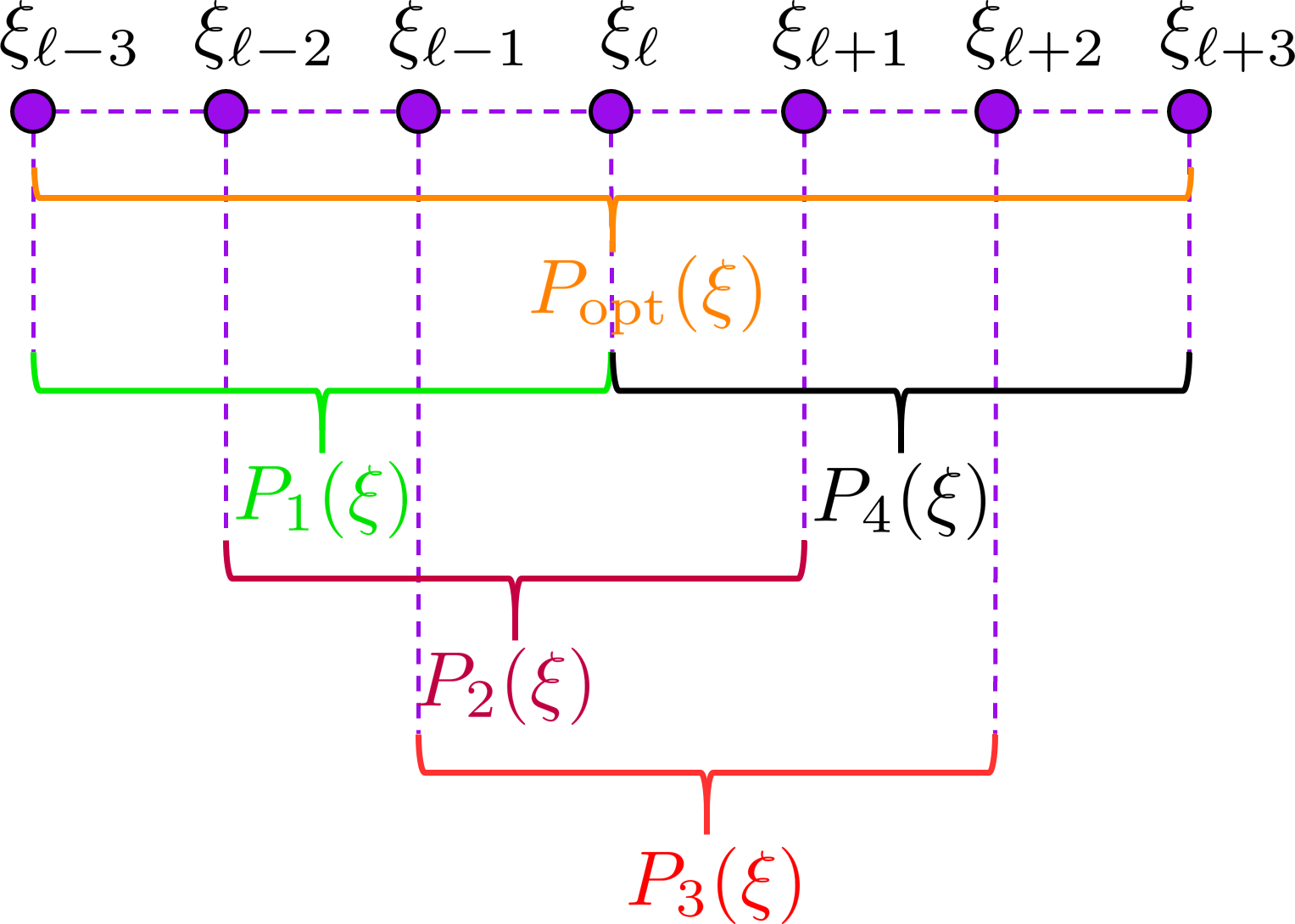}}
\caption{\sf CWENO7 stencil structure for degree 6 polynomial $P_{\rm opt}$ and four sub-stencils for cubic polynomials $P_1$, $P_2$, $P_3$,
$P_4$.\label{fig21}}
\end{figure}

Let us consider a stencil consisting of $7$ equidistant points
$\{\xi_{\ell-3},\xi_{\ell-2},\xi_{\ell-1},\xi_\ell,\xi_{\ell+1},\xi_{\ell+2},\xi_{\ell+3}\}$ with
$\xi_{i+1}-\xi_i\equiv\Delta\xi~\forall i$ (we stress that CWENO7 interpolations can be constructed on any set of nodes, and we consider 
them to be uniform for the sake of simplicity only). There exists a unique interpolating polynomial $P_{\rm opt}(\xi)$ of degree up to $6$
satisfying the interpolation conditions $P_{\rm opt}(\xi_i)=U(\bm x,t;\xi_i),~i=\ell-3,\dots,\ell+3$. However, it is well-known that this
polynomial may be oscillatory especially when the underlying function $U(\cdot,\cdot;\xi)$ is discontinuous. We therefore construct a CWENO7
polynomial
\begin{equation}
{\cal R}_\ell(\xi)=\sum_{k=0}^4\omega_kP_k(\xi),
\label{2.2}
\end{equation}
which is expected to be essentially non-oscillatory on the interval $[\xi_{\ell-\hf},\xi_{\ell+\hf}]$. Here, $P_k$ for $k=1,2,3,4$ are cubic
interpolating polynomials constructed over the sub-stencils outlined in \fref{fig21} and satisfying the corresponding interpolation
conditions $P_k(\xi_i)=U(\bm x,t;\xi_i),~i=\ell+k-4,\dots,\ell+k-1$.

In \eref{2.2}, the polynomial $P_0$ is obtained using a CWENO approach and given by
\begin{equation*}
P_0(\xi)=\frac{1}{d_0}\Big(P_{\rm opt}(\xi)-\sum_{k=1}^4d_kP_k(\xi)\Big),
\end{equation*}
where the coefficients $d_i\in(0,1),~i=0,\dots,4$ are selected to satisfy $d_0+d_1+d_2+d_3+d_4=1$. Our particular choice is
$d_0=\frac{3}{4}$ and $d_1=d_2=d_3=d_4=\frac{1}{16}$. The nonlinear weights $\omega_k$ are calculated based on smoothness indicators
$\beta_k$, effectively minimizing the influence of stencils, in which the discontinuity may be located. Specifically,
\begin{equation*}
\omega_k=\frac{\alpha_k}{\alpha_0+\alpha_1+\alpha_2+\alpha_3+\alpha_4},\quad k=0,\dots,4,
\end{equation*}
where the unnormalized weights $\alpha_k$ are obtained in the CWENOZ manner (see \cite{STP}) and given by
\begin{equation}
\alpha_k=d_k\left[1+\left(\frac{\tau}{\varepsilon+\beta_k}\right)^p\,\right],\quad k=0,\dots,4.
\label{2.3}
\end{equation}
with a typical parameter value $p=2$. In \eref{2.3}, $\varepsilon=(\Delta\xi)^q$ is a small parameter introduced to avoid division by zero
with a typical parameter value $q=3$, and the smoothness indicators $\beta_k$ are obtained as in \cite{JS}:
\begin{equation}
\beta_k=\sum_{i=1}^6(\Delta\xi)^{2i-1}\int\limits_{\xi_{\ell-\hf}}^{\xi_{\ell+\hf}}
\left(\frac{{\rm d}^i}{{\rm d}\xi^i}P_k(\xi)\right)^2{\rm d}\xi,\quad k=0,\dots,4.
\label{2.4}
\end{equation}
The integrals in \eref{2.4} can be evaluated exactly, and this will result in explicit expressions for $\beta_k$, which we omit for the sake
of brevity.
%{\color{blue}
%$$
%\begin{aligned}
%&\beta_0=,\\
%&\beta_1=,\\
%&\beta_2=,\\
%&\beta_3=,\\
%&\beta_4=.
%\end{aligned}
%$$
%these expressions have too many terms (several dozens/hundred each), I do not think we need to include them
%}
Finally, the global smoothness indicator $\tau$ in \eref{2.3} is defined as follows (see \cite{cravero2019}):
\begin{equation*}
\tau=\left|-\beta_1-3\beta_2+3\beta_3+\beta_4\right|.
\end{equation*}

Equipped with the CWENO7 interpolation, we estimate the mean and standard deviation by replacing $U(\bm x,t;\xi)$ in \eref{2.1} with the
piecewise polynomial $\sum_{\ell=1}^L{\cal R}_\ell(\xi)\mbox{{\Large$\chi$}}_{[\xi_{\ell-\hf},\xi_{\ell+\hf}]}(\xi)$, where
$\mbox{{\Large$\chi$}}_{[\xi_{\ell-\hf},\xi_{\ell+\hf}]}$ is a characteristic function of the interval $[\xi_{\ell-\hf},\xi_{\ell+\hf}]$.
This leads to
\begin{equation*}
\begin{aligned}
&\widetilde\mu=\int\limits_{\xi_1}^{\xi_{\frac{3}{2}}}{\cal R}_1(\xi)p(\xi)\,{\rm d}\xi+
\sum_{\ell=2}^{L-1}\int\limits_{\xi_{\ell-\hf}}^{\xi_{\ell+\hf}}{\cal R}_\ell(\xi)p(\xi)\,{\rm d}\xi+
\int\limits_{\xi_{L-\hf}}^{\xi_L}{\cal R}_L(\xi)p(\xi)\,{\rm d}\xi,\\
&\widetilde\sigma=\Bigg(\int\limits_{\xi_1}^{\xi_{\frac{3}{2}}}({\cal R}_1(\xi)-\widetilde\mu)^2p(\xi)\,{\rm d}\xi+
\sum_{\ell=2}^{L-1}\int\limits_{\xi_{\ell-\hf}}^{\xi_{\ell+\hf}}({\cal R}_\ell(\xi)-\widetilde\mu)^2p(\xi)\,{\rm d}\xi+
\int\limits_{\xi_{L-\hf}}^{\xi_L}({\cal R}_L(\xi)-\widetilde\mu)^2p(\xi)\,{\rm d}\xi\Bigg)^\hf,
\end{aligned}
\end{equation*}
which can be evaluated either exactly or highly accurately using a proper Gaussian quadrature. For $\widetilde\mu$, such quadrature reads as
(the quadrature for $\widetilde\sigma$ can be obtained similarly)
$$
\widetilde\mu\approx\sum_{\ell=1}^L\sum_{j=1}^J\gamma_{\ell_j}{\cal R}_{\ell}(\xi_{\ell_j})p(\xi_{\ell_j}),
$$
where $\gamma_{\ell_j}$ and $\xi_{\ell_j}$ are the coefficients and nodes of the Gaussian quadrature. Note that $J$ should be taken
sufficiently large to ensure that the quadrature errors are smaller than the interpolation errors. In the numerical experiments reported in
\S\ref{sec3}, we have taken $J=4$, which corresponds to the eighth-order Gaussian quadrature.
\begin{rmk}
We emphasize that the proposed algorithm is only applicable for $\ell=4,\dots,L-3$. For $\ell<4$ or $\ell>L-3$, one has the following two 
options. First, a one-sided CWENO7 interpolation can be used; see, e.g., \cite{STP}, where a one-sided CWENO approach is discussed. Second, 
if the output function is smooth near the boundary, then one can introduce ghost points across the boundary ($\ell=0,-1,-2$ and
$\ell=L,L+1,L+2$), in which the values of $\bm U$ are obtained using a seventh-order accurate extrapolation.
\end{rmk}

\subsection{Two Random Variables ($s=2$)}
In the case of two random variables $\bm\xi=(\xi,\eta)$, we select a set of collocation points $\{(\xi_\ell,\eta_m)\}$, $\ell=1,\dots,L$,
$m=1,\dots,M$, which form a Cartesian mesh in the random space, and obtain the corresponding output function values
$\{\bm U(\bm x,t;\xi_\ell,\eta_m)\}$. We then estimate the mean of each component $U$ of $\bm U$,
\begin{equation}
\mu[U]:=\iint\limits_\Xi U(\bm x,t;\xi,\eta)p(\xi,\eta)\,{\rm d}\xi{\rm d}\eta,
\label{2.5}
\end{equation}
by applying a Gaussian quadrature in the ``dimension-by-dimension'' manner. To this end, we apply the CWENO7 interpolations in the random
space in the ``dimension-by-dimension'' manner as well.

For given $\bm x$ and $t$, the point values $U(\xi_\ell,\eta_m):=U(\bm x,t;\xi_\ell,\eta_m)$ are available. We first fix $\eta_m$ and
perform the CWENO7 interpolations in the $\xi$-direction for each $m$ to compute the values $U(\xi_{\ell_j},\eta_m)$ at the Gaussian nodes
$\xi_{\ell_j}$ in the $\xi$-direction. We then fix $\xi_{\ell_j}$ and perform the CWENO7 interpolations in the $\eta$-direction for each
$\ell_j$ to obtain $U(\xi_{\ell_j},\eta_{m_r})$, where $\eta_{m_r}$ are the corresponding Gaussian nodes in the $\eta$-direction. Finally,
equipped with the values $U(\xi_{\ell_j},\eta_{m_r})$ at the two-dimensional (2-D) Gaussian nodes $(\xi_{\ell_j},\eta_{m_r})$, we apply a
Gaussian quadrature of order $2J$ to the integral in \eref{2.5} to obtain
$$
\widetilde\mu\approx\sum_{\ell=1}^L\sum_{m=1}^M\sum_{j=1}^J\sum_{r=1}^J\gamma_{\ell_j}\gamma_{m_r}
U(\xi_{\ell_j},\eta_{m_r})p(\xi_{\ell_j},\eta_{m_r}).
$$

%%%%%%%%%%%%%
% SECTION 3 %
%%%%%%%%%%%%%
\section{Numerical Examples}\label{sec3}
We now test the proposed surrogate modeling approach on several functions $\bm U$, which are assumed to represent discrete approximations of
the solutions of \eref{1.1}. In the first five examples, we will assume that the exact solution of \eref{1.1} is given by scalar functions
$U$ and they are either two smooth functions of one variable,
\begin{equation}
U(\xi)=3\cos(\pi\xi)
\label{3.1}
\end{equation}
and
\begin{equation}
U(\xi)=\tanh(9\xi)+0.5\xi,
\label{3.2}
\end{equation}
or a smooth function of two variables,
\begin{equation}
U(\xi,\eta)=3\cos(\pi\xi)\cos(\pi\eta),
\label{3.3}
\end{equation}
or a discontinuous function of one variable,
\begin{equation}
U(\xi)=\begin{cases}\phantom{-}3\cos(\pi\xi),&\xi<0.1,\\-3\cos(\pi\xi),&\xi>0.1,\end{cases}
\label{3.4}
\end{equation}
or a discontinuous function of two variables,
\begin{equation}
U(\xi,\eta)=\begin{cases}\phantom{-}3\cos(\pi\xi)\cos(\pi\eta),&\xi<0.1,\,\eta<0.1,\\
-3\cos(\pi\xi)\cos(\pi\eta),&\mbox{otherwise}.
\end{cases}
\label{3.5}
\end{equation}
We will construct surrogate models $\widetilde U$ for $U$ using both the gPC expansion and CWENO7 interpolation, and then will measure the
differences $\widetilde U-U$ as well as the differences between the approximated and exact means, standard deviations, and PDFs. Notice that
to compare the PDFs, we will use the histogram method on a very large number ($3\times10^7$ in the case of a single random variable and
$10^4\times10^4$ in the case of two random variables) of bins applied to both $\widetilde U$ and $U$.

In the final two examples, $\bm U:=(h,hu)^\top$ will be obtained as a numerical solution of the one-dimensional (1-D) Saint-Venant system of
shallow water equations with uncertainties,
\begin{equation}
\bm U_t+\bm F(\bm U)_x=\bm S(\bm U,x;\bm\xi),\quad\bm F(\bm U)=\Big(hu,hu^2+\frac{g}{2}h^2\Big)^\top,\quad\bm S=(0,-ghZ_x)^\top,
\label{3.6}
\end{equation}
where $h(x,t;\bm\xi)$ is the water depth, $u(x,t;\bm\xi)$ is the velocity, $Z(x;\bm\xi)$ is the bottom topography, and $g$ is the constant
acceleration due to gravity (in the examples below, we take $g=1$).

We will select uniform collocation points $\bm\xi_\ell$ and at each of them, we will numerically solve \eref{3.6} using the second-order
semi-discrete central-upwind scheme from \cite{KPshw} on a uniform spatial mesh with the nodes denoted by $x_j$. This way, at every grid
point $x_j$ in space and at the final computational time $T$, we will generate the discrete function
$w(x_j,T;\bm\xi_\ell):=h(x_j,T;\bm\xi_\ell)+Z(x_j;\bm\xi_\ell)$ representing the water surface. We will compute the mean and standard
deviation for $w$ and also approximate its PDF.

\subsubsection*{Example 1---Smooth Function \eref{3.1}}
We examine two distinct cases involving the uniformly and normally distributed random variable $\xi$.

\paragraph{Test 1 (Uniform Distribution).} We consider a random variable $\xi$ uniformly distributed over the interval $[-1,1]$, that is,
$\xi\sim{\cal U}(-1,1)$. In this case, the analytical PDF of $U$ can be derived using the transformation method \cite{papoulis2002} as
follows:
\begin{equation*}
p(\xi)=\frac{1}{\pi\sqrt{9-U^2(\xi)}},
\end{equation*}
which is singular at $u=\pm3$. The exact values of mean and standard deviation of $U$ are then $\mu=0$ and $\sigma=\sqrt{4.5}$,
respectively.

\begin{figure}[ht!]
\centerline{\includegraphics[trim=0.2cm 0.3cm 0.2cm 0.1cm, clip, width=0.31\textwidth]{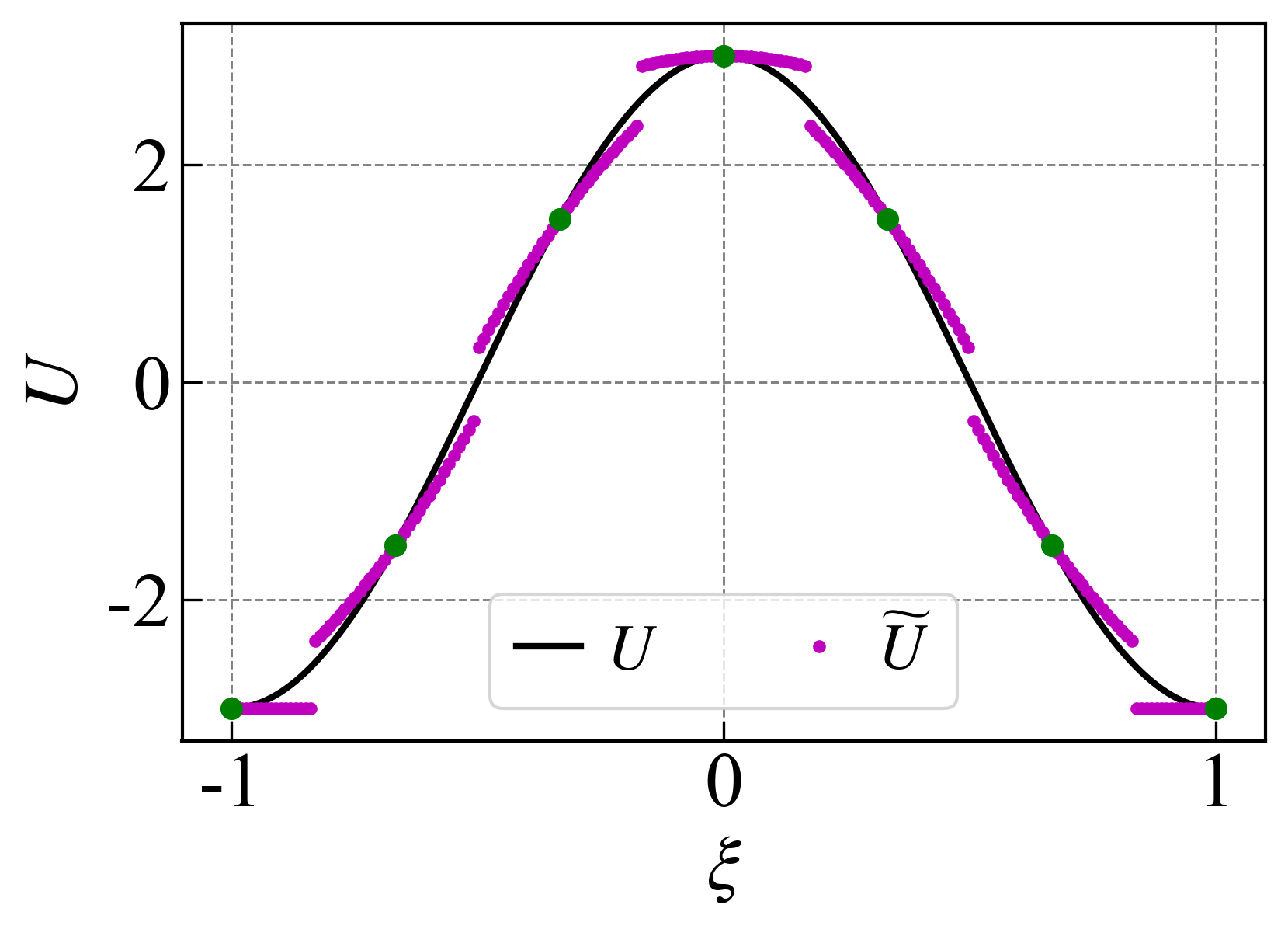}\hspace*{0.3cm}
            \includegraphics[trim=0.2cm 0.3cm 0.2cm 0.1cm, clip, width=0.31\textwidth]{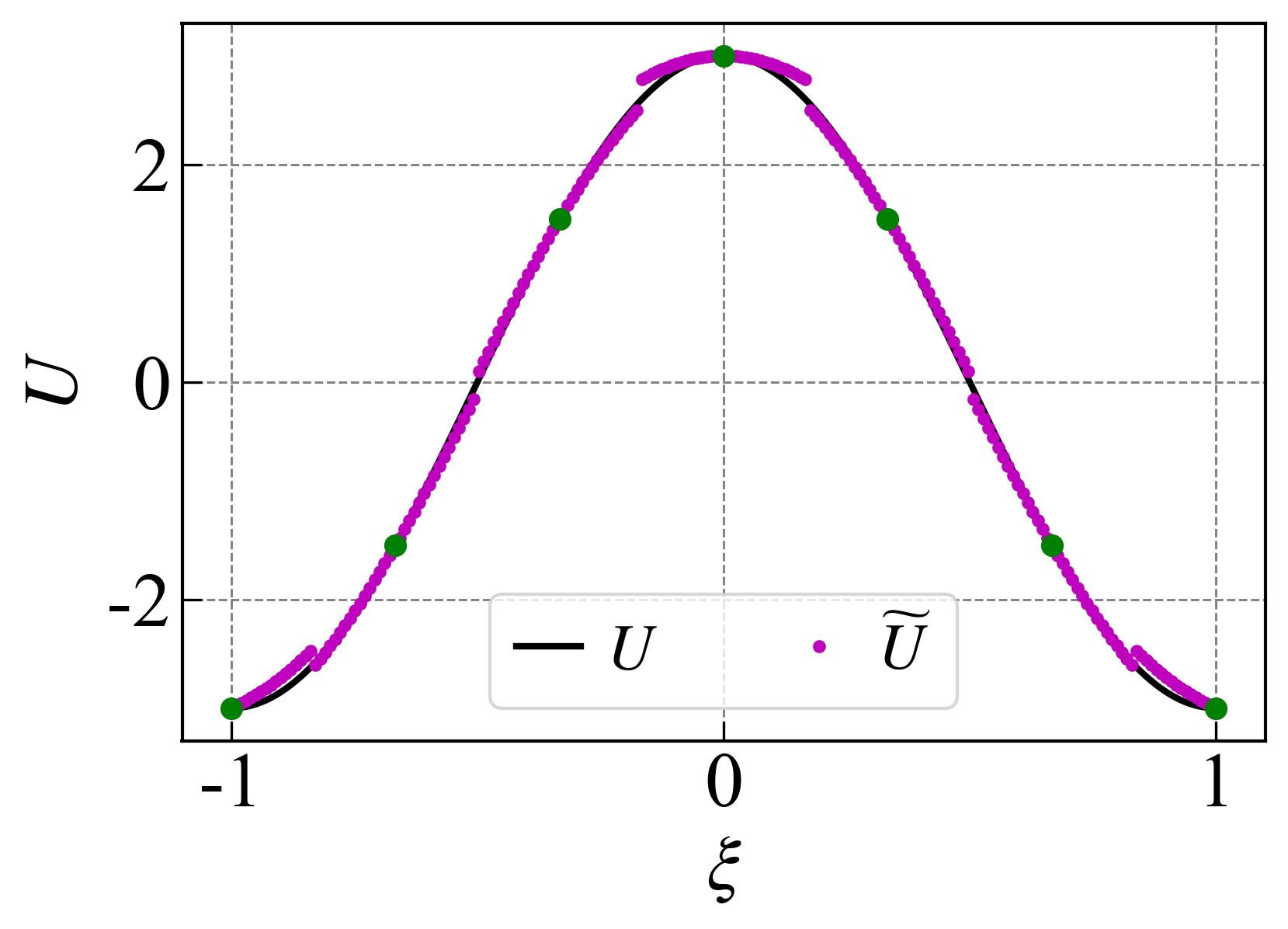}\hspace*{0.3cm}
            \includegraphics[trim=0.2cm 0.3cm 0.2cm 0.1cm, clip, width=0.31\textwidth]{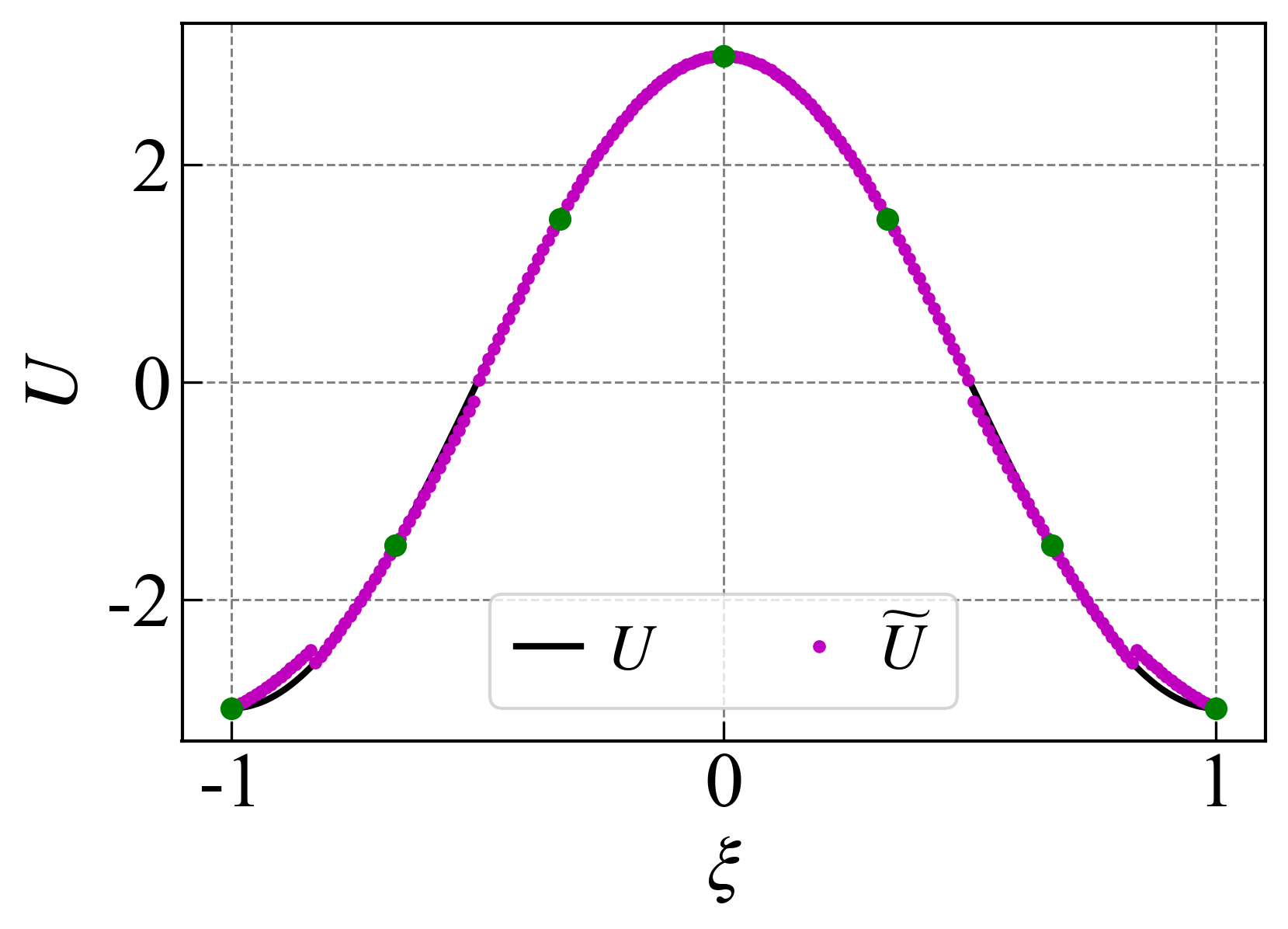}}
\caption{\sf Example 1, Test 1: Interpolations obtained for $L=7$ using CWENO3 (left), CWENO5 (middle), and CWENO7 (right). Green dots
represent the output data being interpolated.\label{fig33}}
\end{figure}

We first provide a motivation on why a high-order CWENO interpolation is needed to construct a good CWENO-based surrogate model. To this
end, we use the third-, fifth-, and seventh-order CWENO interpolations (CWENO3, CWENO5, and CWENO7) and plot the results obtained on a
uniform mesh with $L=7$ in \fref{fig33}. As one can see, when the mesh in $\xi$ is coarse, the higher-order CWENO7 interpolation is visibly
``smoother'' as its jumps at the cell interfaces $\xi=\xi_\lph$ are substantially smaller compared to those in the lower-order CWENO3 and
CWENO5 interpolations. Even though the size of the jumps decreases when $L$ increases, the CWENO7 interpolation seems to be a reasonable
choice as it is sufficiently accurate and at the same time not increasingly computationally expensive.

We now construct surrogate models $\widetilde U(\xi)$, which are based on the corresponding sets of collocation points with $L=7,9,11,13,15,17$,
and $19$. Notably, the gPC approach employs Gauss-Legendre collocation points tailored for $\xi\sim{\cal U}(-1,1)$, while the CWENO7 method
utilizes uniformly distributed points for $\xi_\ell$.

\begin{figure}[ht!]
\centerline{\includegraphics[trim=0.2cm 0.3cm 0.2cm 0.1cm, clip, width=0.35\textwidth]{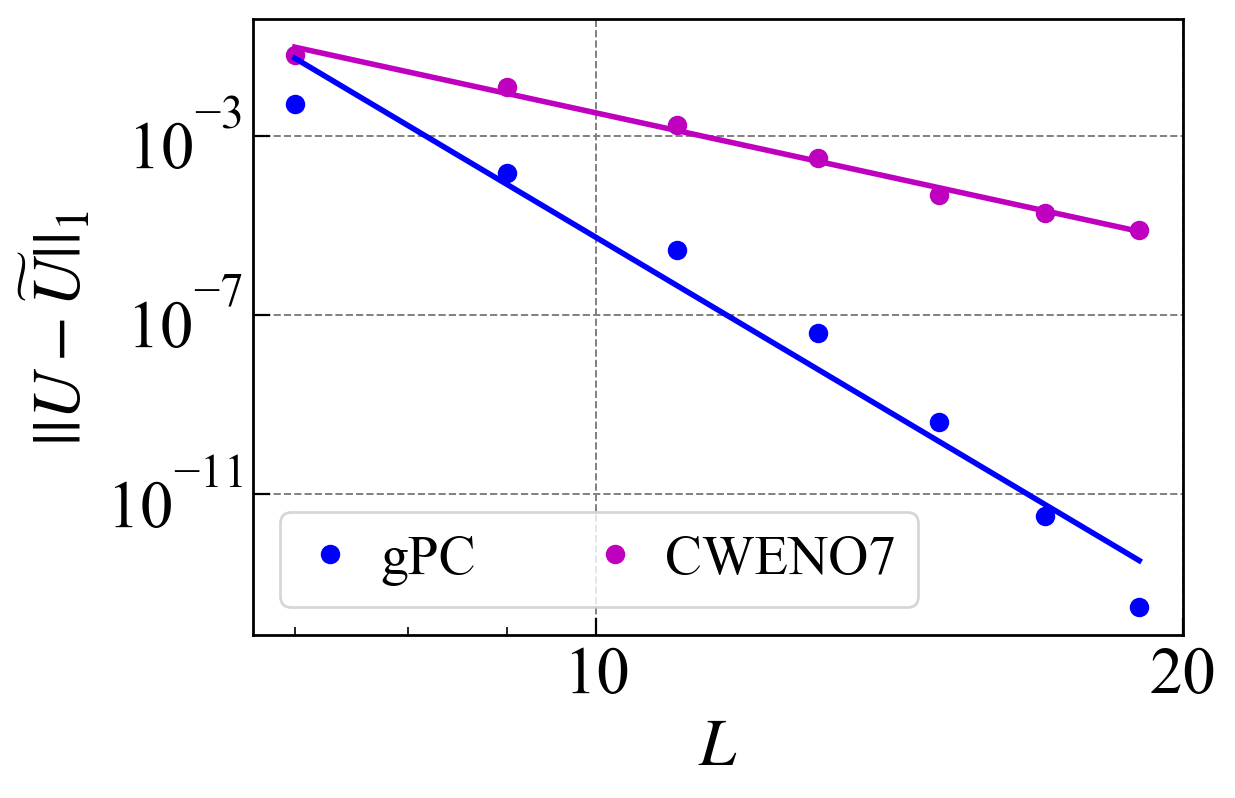}}
\caption{\sf Example 1, Test 1: $L^1$-errors for the gPC expansion and CWENO7 interpolation as functions of $L$ and the corresponding
power-law fits (solid lines).\label{fig31}}
\end{figure}

\smallskip
\underline{\it Convergence of the surrogate models}. We measure the discrete $L^1$-norm $\|U-\widetilde U\|_1$ by the Simpson rule using
$20000$ uniform subintervals on $[-1,1]$. We then use the line fitting \cite{DFS19} to find $k$ such that
\begin{equation*}
\|U-\widetilde U\|_1\approx CL^{-k}.
\end{equation*}
The obtained results are illustrated in \fref{fig31}, where one can see that, as expected for smooth output functions, the gPC expansion
demonstrates a substantially higher convergence rate compared to the CWENO7 interpolation: the corresponding exponents $k$ are about $25.8$
and $9.5$. 

\begin{figure}[ht!]
\centerline{\includegraphics[trim=0.2cm 0.3cm 0.2cm 0.1cm, clip, width=0.35\textwidth]{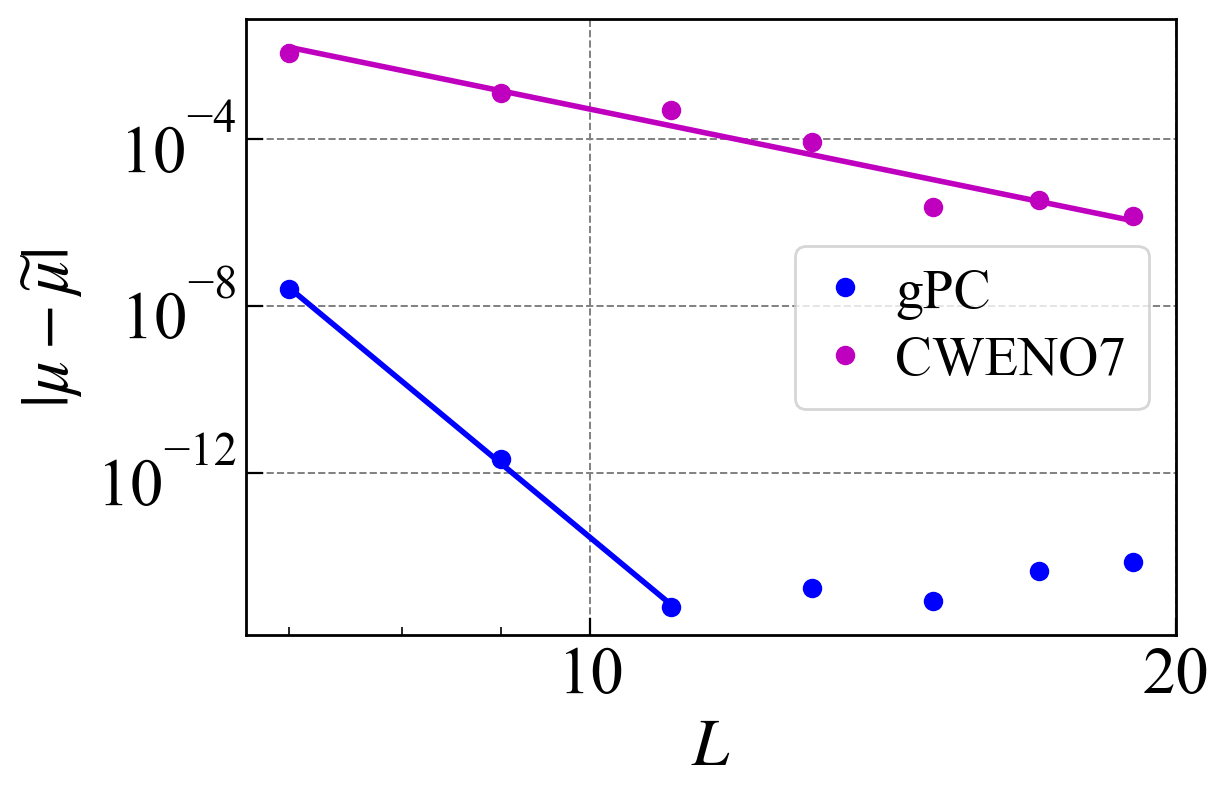}\hspace*{0.5cm}
	    \includegraphics[trim=0.2cm 0.3cm 0.2cm 0.1cm, clip, width=0.35\textwidth]{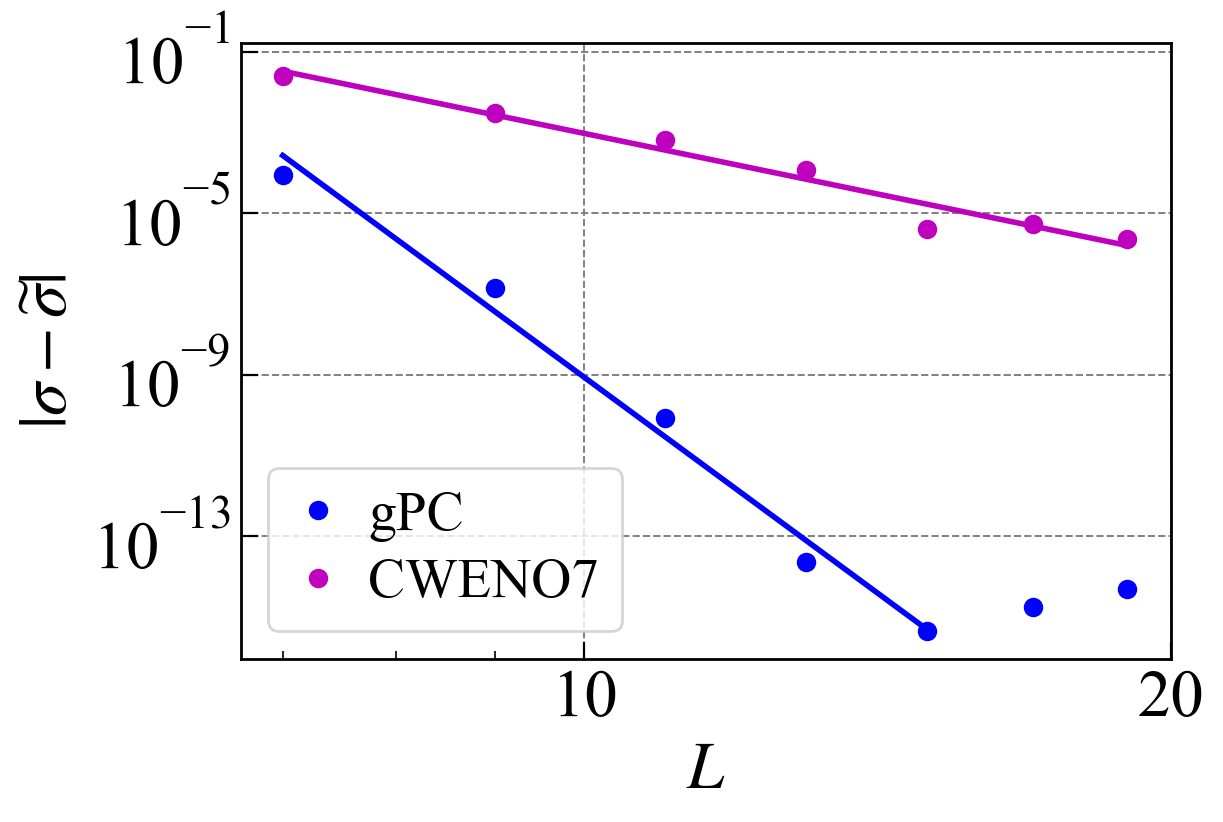}}
\caption{\sf Example 1, Test 1: Errors in $\mu$ (left) and $\sigma$ (right) for the gPC expansion and CWENO7 interpolation as functions of
$L$ and the corresponding power-law fits (solid lines).\label{fig32}}
\end{figure}

\smallskip
\underline{\it Convergence of the mean and standard deviation}. In addition, we check the convergence of the surrogate estimates for the
mean $\mu$ and standard deviation $\sigma$. The differences $|\mu-\widetilde\mu|$ and $|\sigma-\widetilde\sigma|$ as functions of $L$ are
plotted in \fref{fig32} along with the corresponding power-law fits (notice the error saturation for the gPC expansion, which occurs at 
$L=11$) with the exponents $38.8$ and $9.6$ (for $\mu$) and $35.5$ and $10$ (for $\sigma$).

\begin{figure}[ht!]
\centerline{\includegraphics[trim=0.2cm 0.3cm 0.2cm 0.2cm, clip, width=0.31\textwidth]{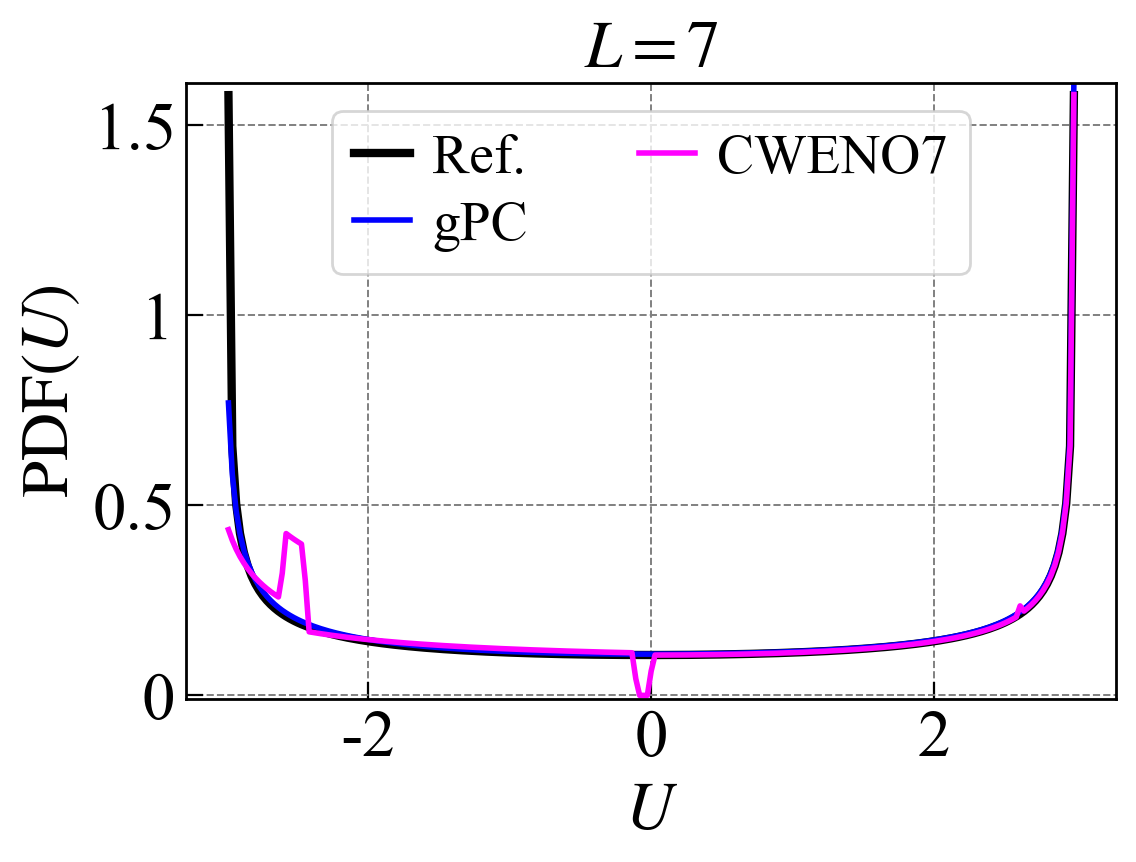}\hspace*{0.3cm}
            \includegraphics[trim=0.2cm 0.3cm 0.2cm 0.2cm, clip, width=0.31\textwidth]{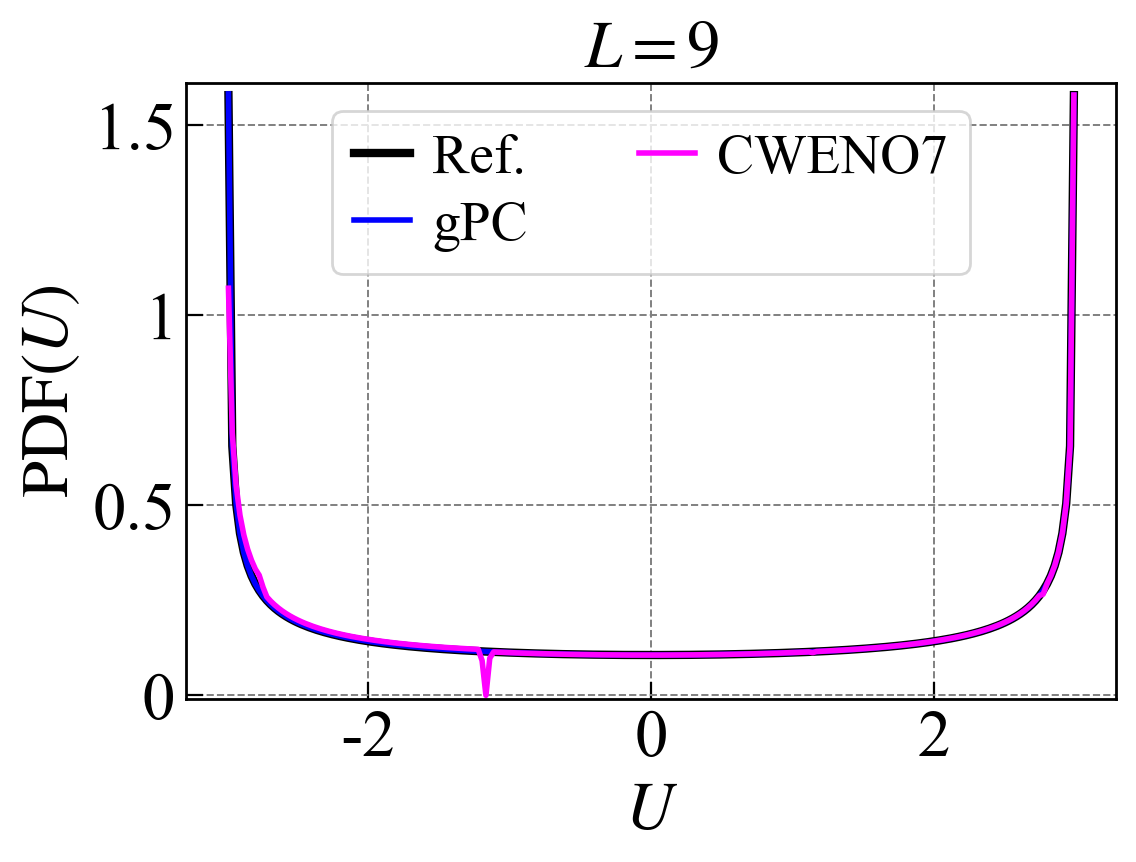}\hspace*{0.3cm}
	    \includegraphics[trim=0.2cm 0.3cm 0.2cm 0.2cm, clip, width=0.31\textwidth]{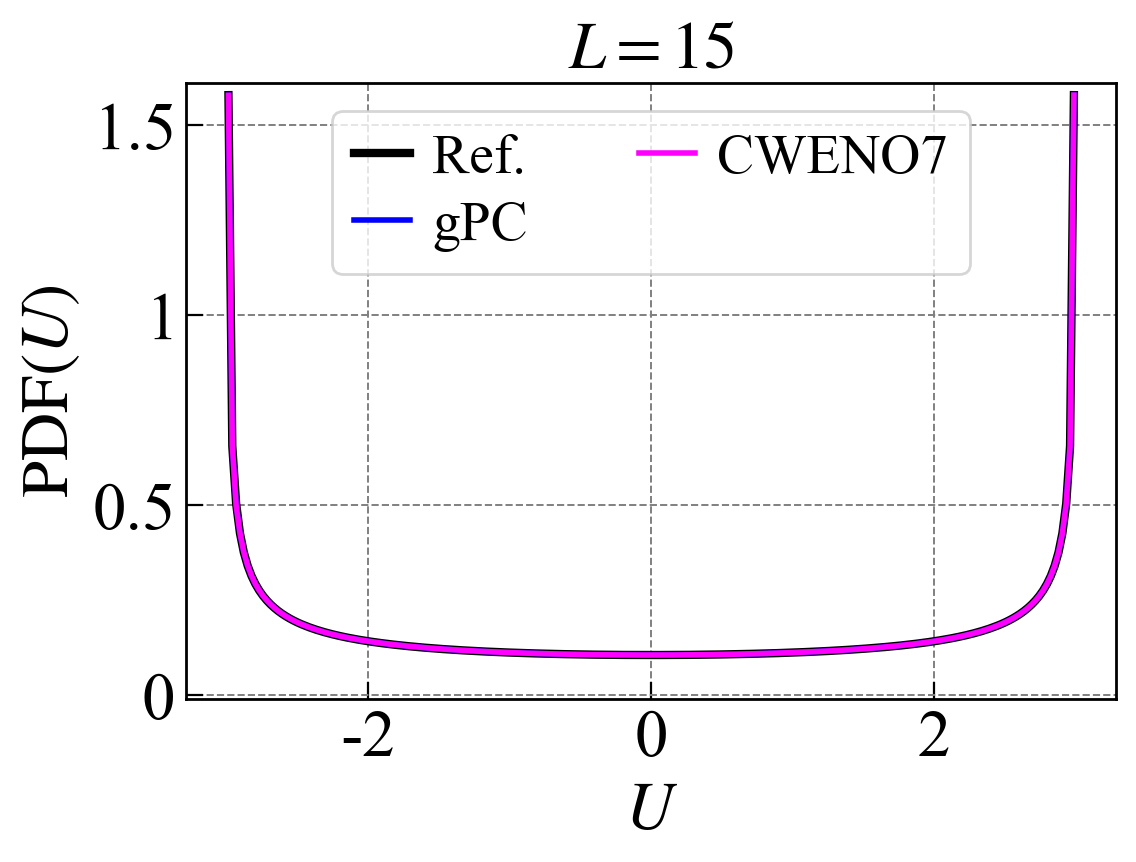}}
\caption{\sf Example 1, Test 1: Estimated PDFs for the gPC and CWENO7 interpolations $\widetilde U$ together with the reference PDF,
reconstructed from $U$ for $L=7$ (left), $9$ (middle), and $15$ (right).\label{fig34}}
\end{figure}

\smallskip
\underline{\it Convergence of PDF}. Next, we examine the convergence of the surrogate-based PDF approximations, which are shown in
\fref{fig34} for $L=7$, $9$, and $15$. As one can see, the PDFs computed by the CWENO7 interpolation using $L=7$ and $9$ contain
discrepancies, which are attributed to the fact that the CWENO7 interpolation is a piecewise polynomial that contains jumps at each cell
interface $\xi=\xi_\lph$. For instance, one can see in \fref{fig33} (right) the lack of monotonicity in $\widetilde U$, which causes the
discrepancy for the values of $U\in(-3,-2)$, and the gap in the values for $U\sim0$, which causes the drop of the PDF for $\widetilde U$
around $U=0$. The magnitude of the discrepancies is clearly smaller for $L=9$, and the PDF for $\widetilde U$ is visibly indistinguishable
from the PDF for $U$ when $L=15$.

\begin{figure}[ht!]
\centerline{\includegraphics[trim=0.2cm 0.3cm 0.2cm 0.2cm, clip, width=0.35\textwidth]{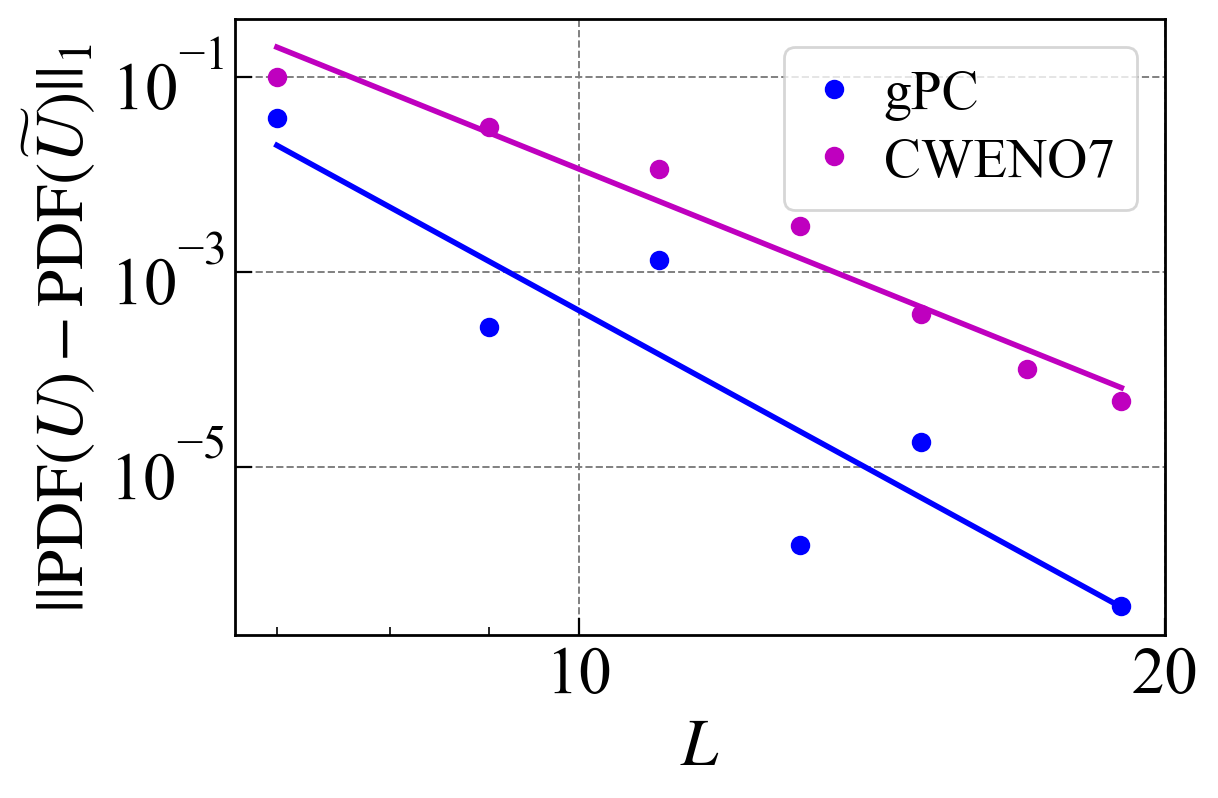}}
\caption{\sf Example 1, Test 1: $L^1$-error for the PDFs as functions of $L$ and the corresponding power-law fits (solid lines).
\label{fig35}}
\end{figure}

The performance of both gPC- and CWENO7-based surrogate models in terms of the accuracy of PDF approximations is further analyzed in
\fref{fig35}, where we plot the $L^1$-errors in PDFs as functions of $L$ along with the corresponding power-law fits with the exponents
$12.4$ (gPC) and $7.8$ (CWENO7). It is evident that the convergence rate for the PDFs is lower compared to the convergence rates observed
for mean and standard deviation. This reduction in convergence rate can be attributed to the inherent limitations of the histogram method,
particularly its sensitivity to finite bin widths and sampling density.

\paragraph{Test 2 (Normal Distribution).} Next, we consider normally distributed random variable $\xi\sim{\cal N}(0,1)$. 
Implementation-wise, this only affects the gPC expansion, as now it utilizes the Gauss-Hermite quadrature points, whereas CWENO7 
interpolation still employs equally spaced points $\xi_\ell$, which are now uniformly distributed over a larger interval $[-6,6]$. Unlike 
Test 1, the analytical expression for $p(\xi)$ is bulky and is therefore not provided. However, one can very accurately compute the mean 
and standard deviation, which are $\mu\approx0.021575650067$ and $\sigma\approx\sqrt{4.499534503363}$.

Overall, the results obtained in this test are similar to those from Test 1. However, due to the wider domain required by the normal
distribution, a larger number of collocation points is necessary to achieve a comparable accuracy, particularly with the CWENO7
interpolation.

\begin{figure}[ht!]
\centerline{\includegraphics[trim=0.2cm 0.3cm 0.2cm 0.1cm, clip, width=0.35\textwidth]{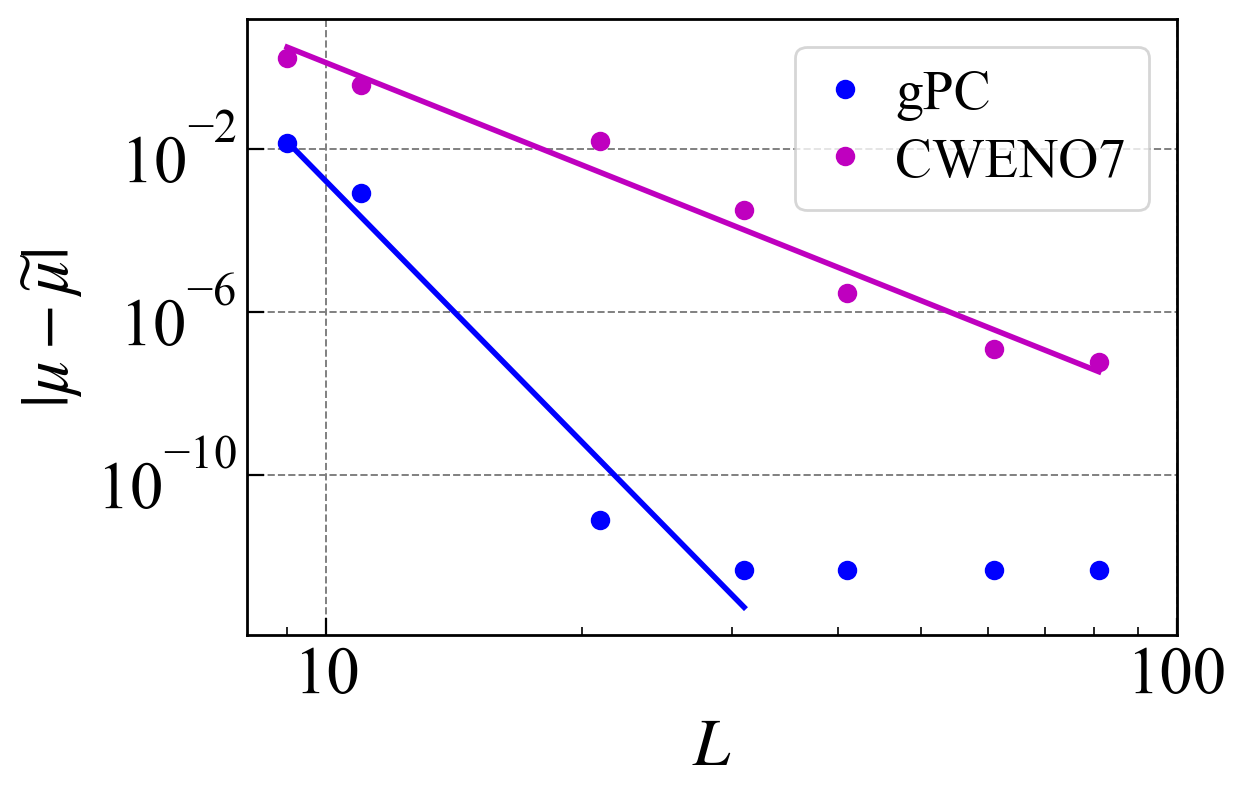}\hspace*{0.5cm}
            \includegraphics[trim=0.2cm 0.3cm 0.2cm 0.1cm, clip, width=0.35\textwidth]{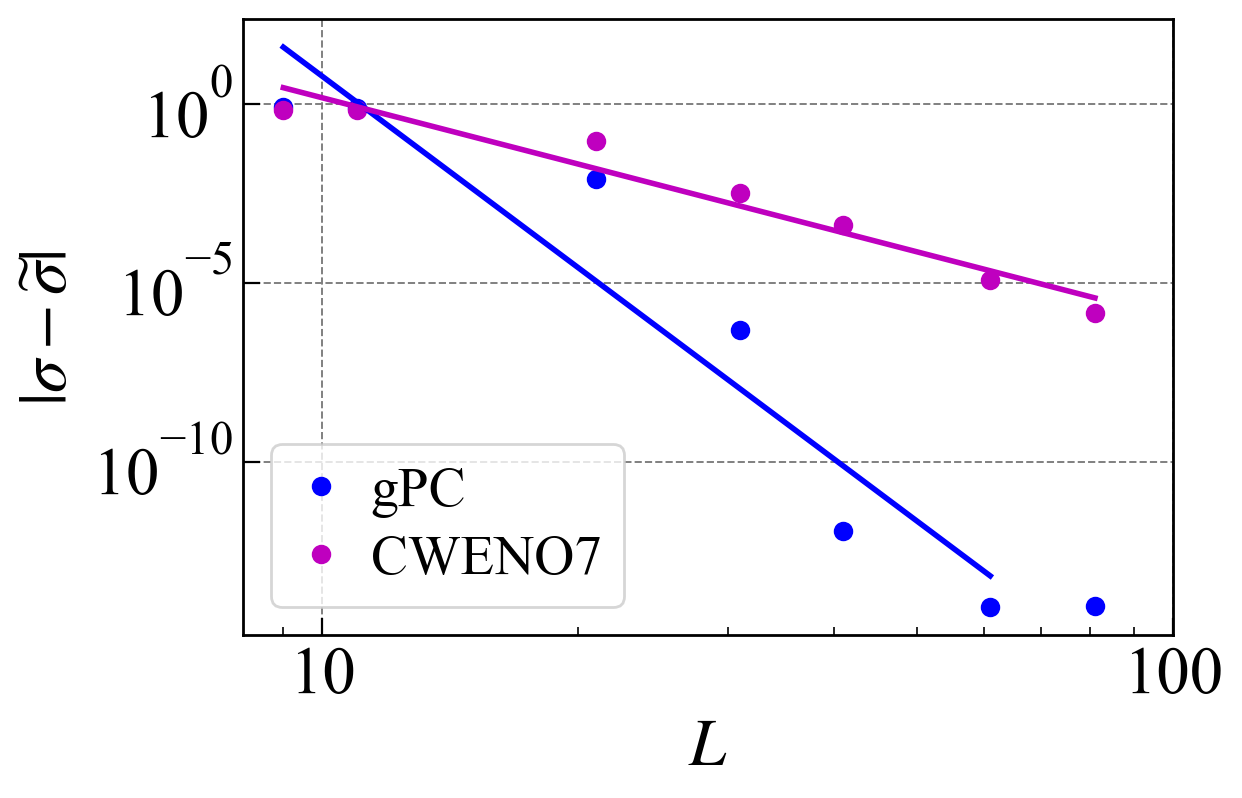}}
\caption{\sf Example 1, Test 2: Errors in $\mu$ (left) and $\sigma$ (right) for the gPC expansion and CWENO7 interpolation as functions of
$L$ and the corresponding power-law fits (solid lines).\label{fig36}}
\end{figure}

\smallskip
\underline{\it Convergence of the mean and standard deviation}. The differences $|\mu-\widetilde\mu|$ and $|\sigma-\widetilde\sigma|$ as
functions of $L$ (we take $L=9$, $11$, $21$, $31$, $41$, $61$, and $81$) are depicted in \fref{fig36} along with the corresponding
power-law fits. The corresponding exponents for the gPC expansion are $21.3$ (for $\mu$) and $17.8$ (for $\sigma$), while for the CWENO7
interpolation they are $8.3$ (for $\mu$) and $6.2$ (for $\sigma$).

\begin{figure}[ht!]
\centerline{\includegraphics[trim=0.2cm 0.3cm 0.2cm 0.2cm, clip, width=0.31\textwidth]{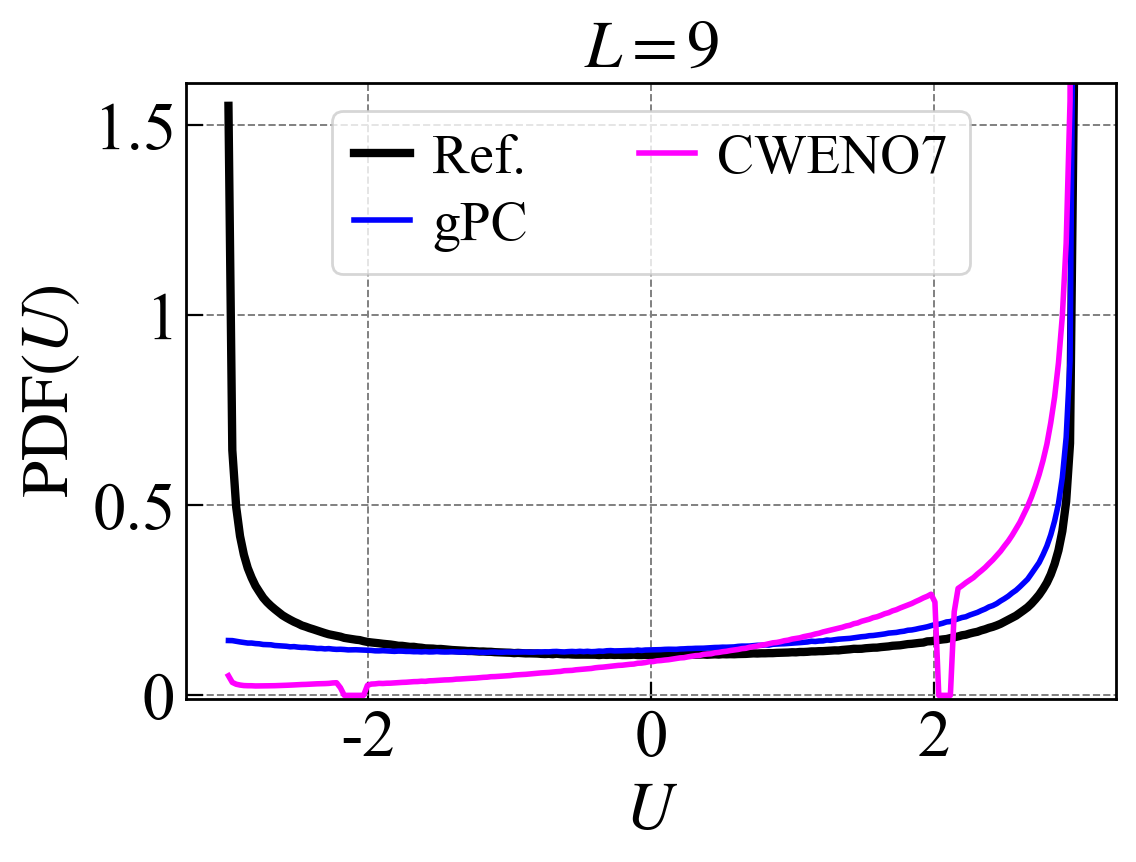}\hspace*{0.3cm}
            \includegraphics[trim=0.2cm 0.3cm 0.2cm 0.2cm, clip, width=0.31\textwidth]{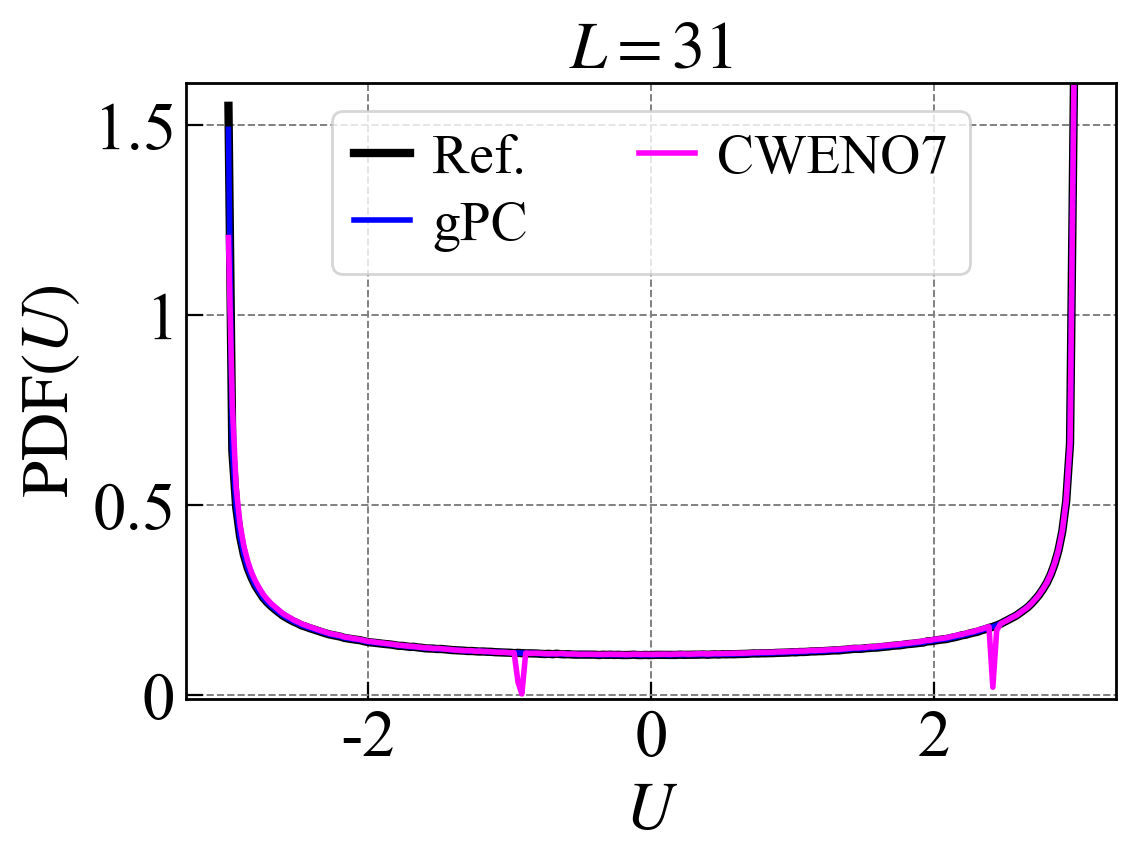}\hspace*{0.3cm}
            \includegraphics[trim=0.2cm 0.3cm 0.2cm 0.2cm, clip, width=0.31\textwidth]{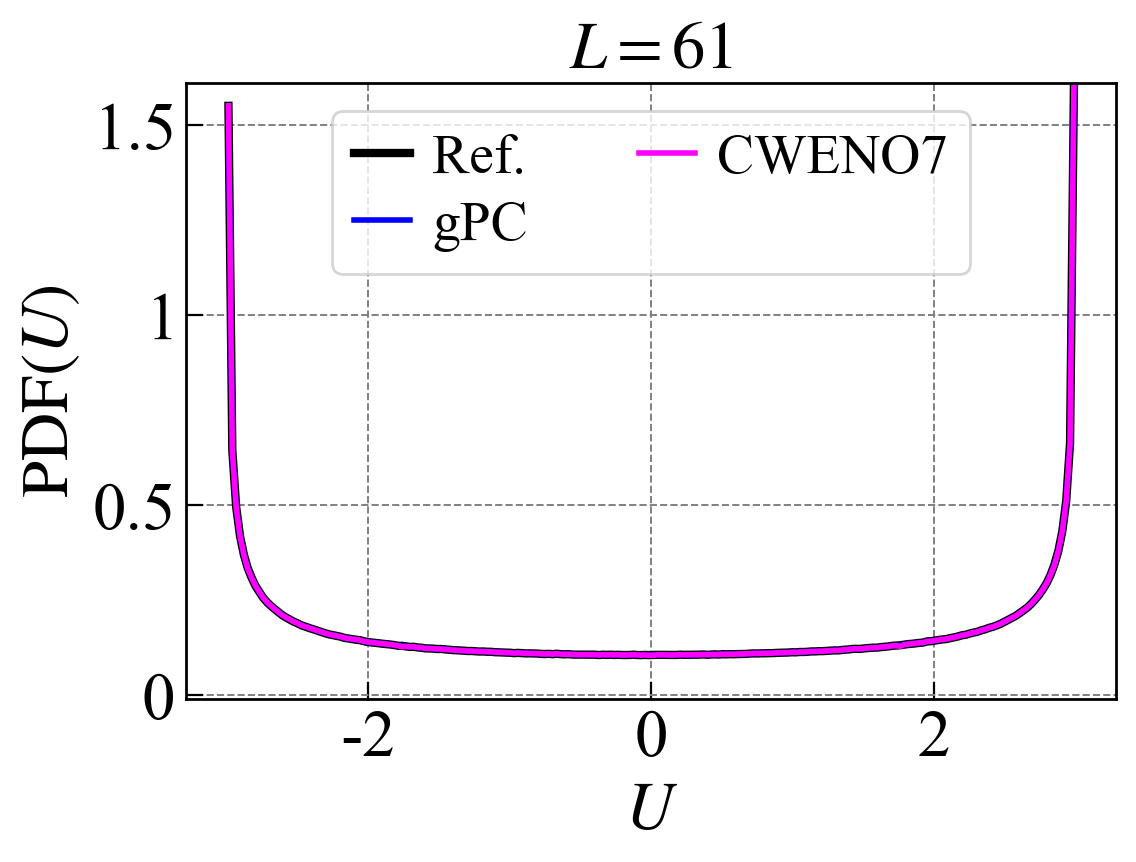}}
\caption{\sf Example 1, Test 2: Estimated PDFs for the gPC and CWENO7 interpolations $\widetilde U$ together with the reference PDF,
reconstructed from $U$ for $L=9$ (left), $31$ (middle), and $61$ (right).\label{fig37}}
\end{figure}

\begin{figure}[ht!]
\centerline{\includegraphics[trim=0.2cm 0.3cm 0.2cm 0.2cm, clip, width=0.35\textwidth]{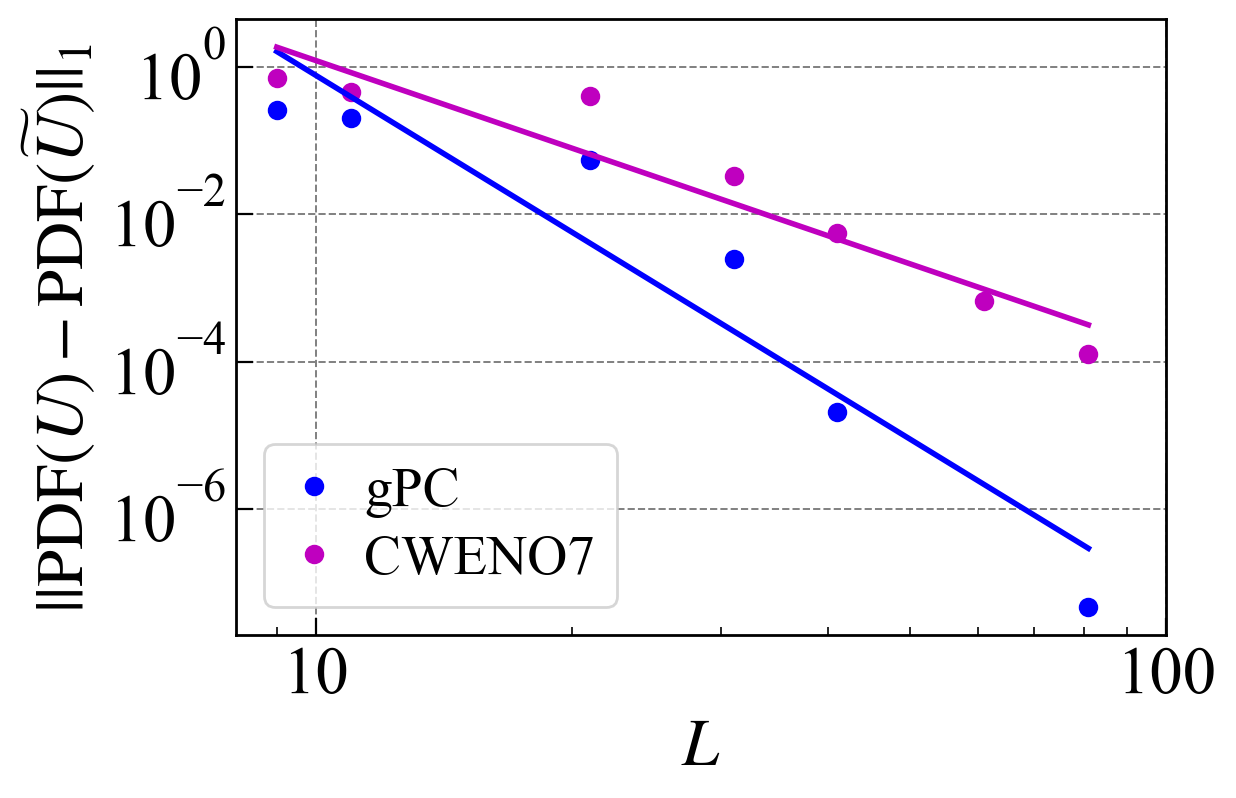}}
\caption{\sf Example 1, Test 2: $L^1$-error for the PDFs as functions of $L$ and the corresponding power-law fits (solid lines).
\label{fig38}}
\end{figure}

\smallskip
\underline{\it Convergence of PDF}. Next, we examine the convergence of the surrogate-based PDF approximations, which are shown in
\fref{fig37} for $L=9$, $31$, and $61$. As one can see, similarly to Test 1, the PDFs computed by the CWENO7 interpolation using small
($L=9$) or intermediate ($L=31$) number of collocation points, the resulting PDFs contain discrepancies, which disappear when the number of
collocation points is sufficiently large. \fref{fig38} displays the $L^1$-errors in PDFs as functions of $L$ along with the corresponding
power-law fits with the exponents $7.1$ (gPC) and $4$ (CWENO7). While both the gPC- and CWENO7-based surrogate models exhibit convergence
rates comparable to those observed in Test 1, the error magnitude for the CWENO7 interpolation is notably higher.

\subsubsection*{Example 2---Smooth Function \eref{3.2}}
We consider a random variable $\xi\sim{\cal U}(-1,1)$, for which the corresponding values of the mean and standard deviation are
$\mu[U]\approx0$ and $\sigma[U]\approx\sqrt{1.467145270396}$.

\begin{figure}[ht!]
\centerline{\includegraphics[trim=0.2cm 0.3cm 0.2cm 0.1cm, clip, width=0.35\textwidth]{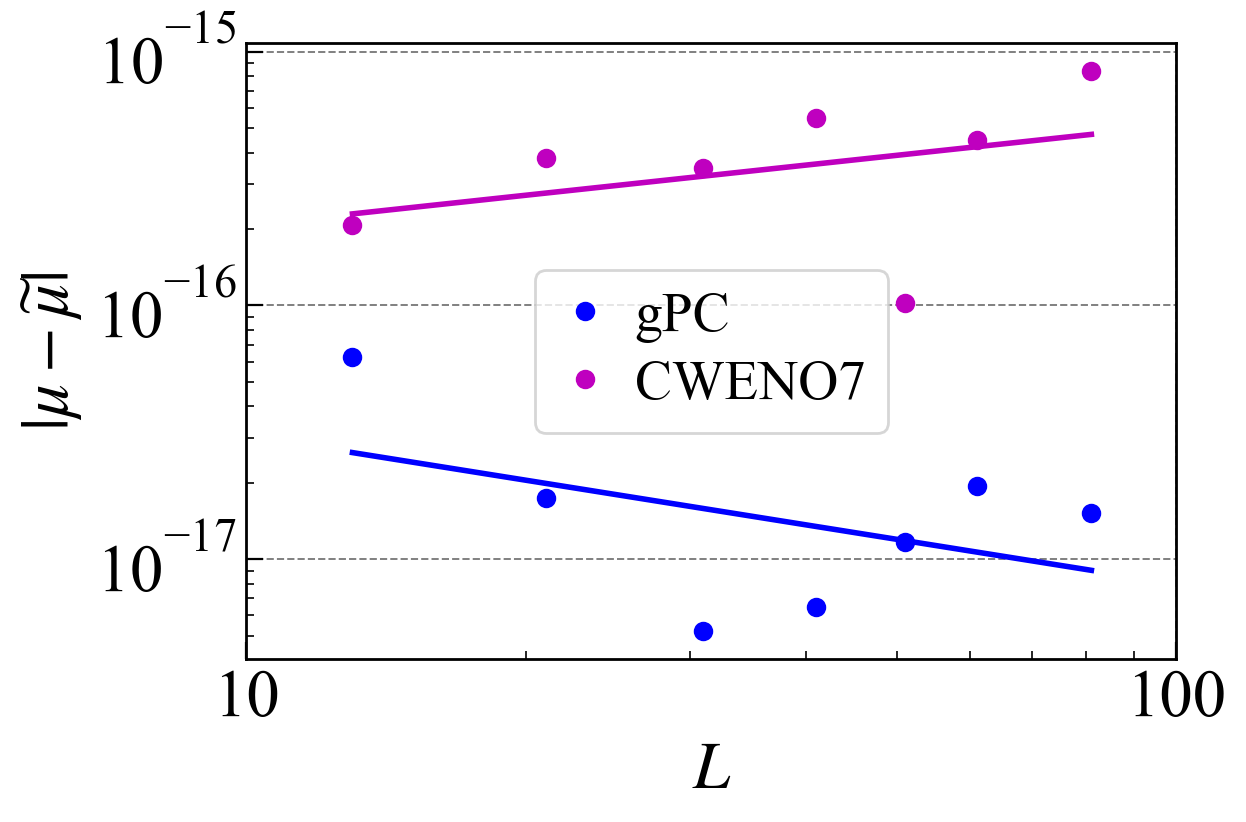}\hspace*{0.5cm}
            \includegraphics[trim=0.2cm 0.3cm 0.2cm 0.1cm, clip, width=0.35\textwidth]{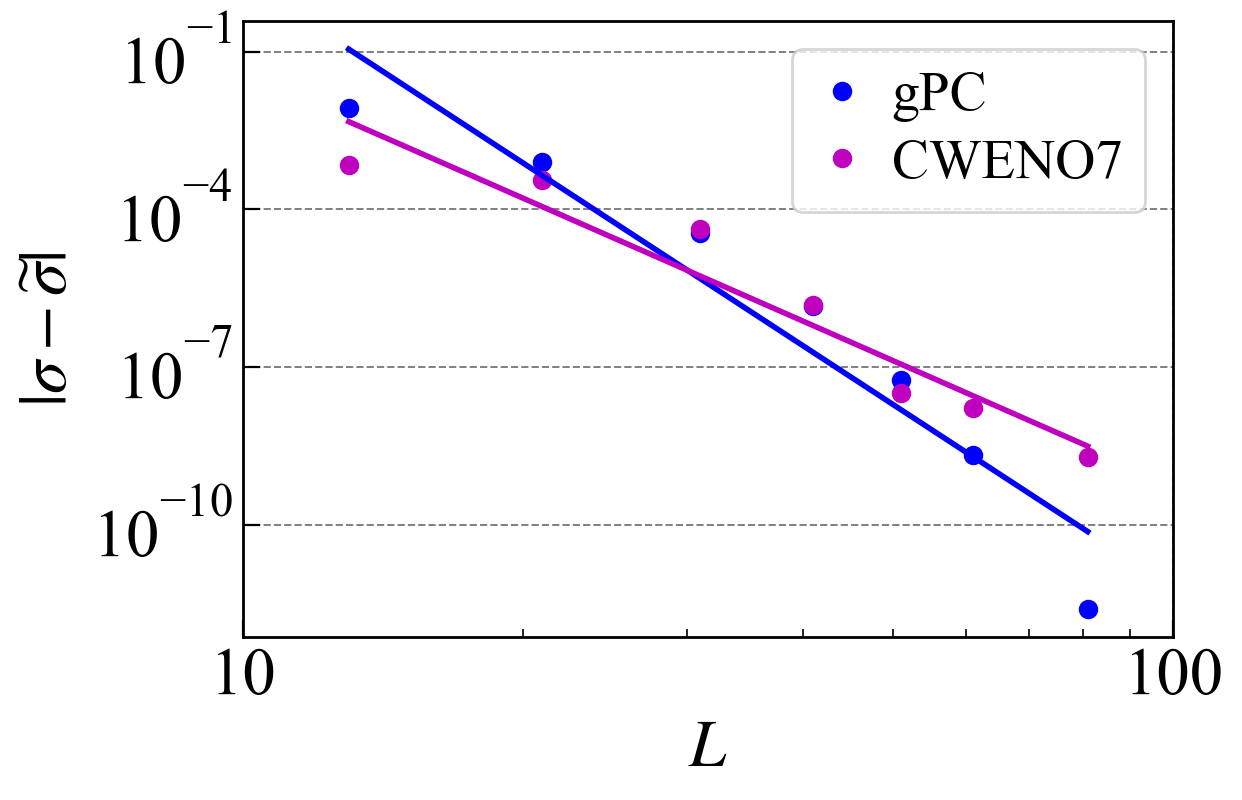}}
\caption{\sf Example 2: Errors in $\mu$ (left) and $\sigma$ (right) for the gPC expansion and CWENO7 interpolation as functions of $L$ and
the corresponding power-law fits (solid lines).\label{fig39}}
\end{figure}

\smallskip
\underline{\it Convergence of the mean and standard deviation}. The differences $|\mu-\widetilde\mu|$ and $|\sigma-\widetilde\sigma|$ as
functions of $L$ (we take $L=21$, $31$, $41$, $51$, $61$, and $81$) are depicted in \fref{fig39} along with the corresponding power-law
fits. As one can see, the errors in $\mu$ for both the gPC expansion and CWENO7 interpolation are at the level of machine zero even for 
$L=21$, while the calculated exponents for $\sigma$ are $11.6$ and $7.8$, respectively.

\begin{figure}[ht!]
\centerline{\includegraphics[trim=0.2cm 0.3cm 0.2cm 0.2cm, clip, width=0.31\textwidth]{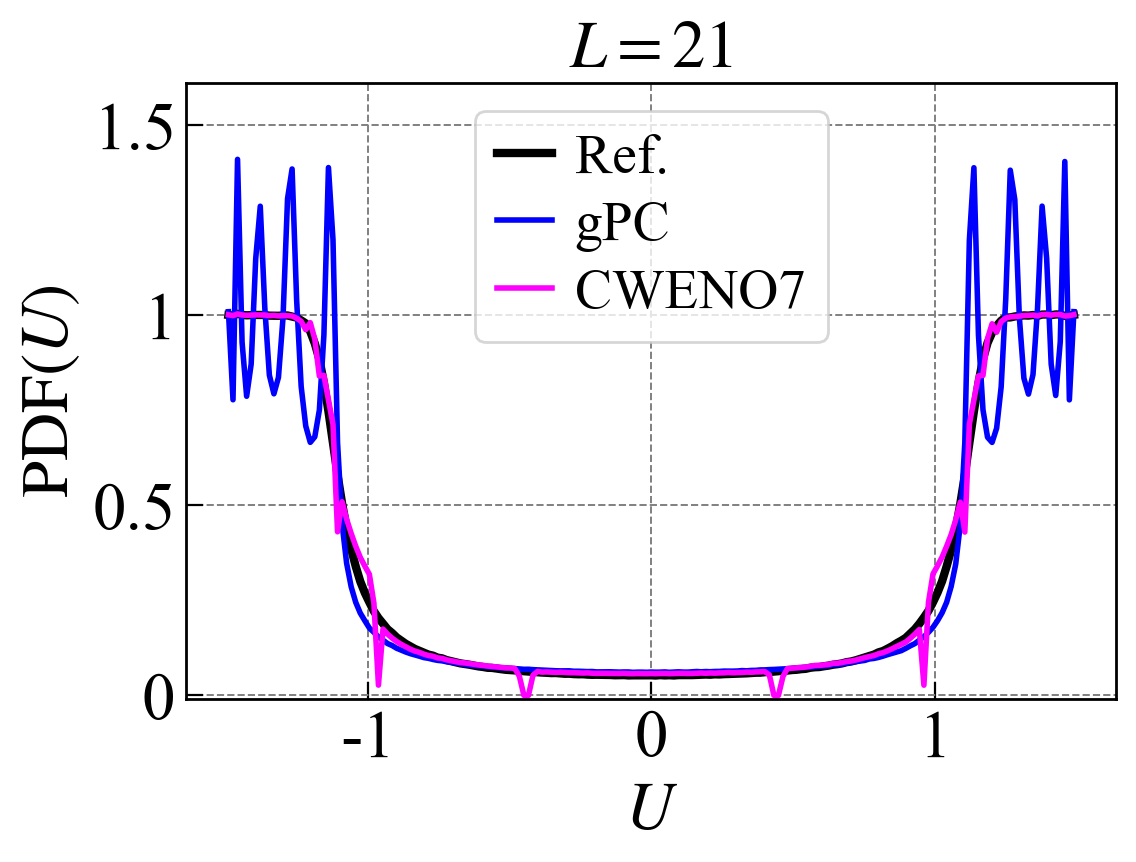}\hspace*{0.3cm}
            \includegraphics[trim=0.2cm 0.3cm 0.2cm 0.2cm, clip, width=0.31\textwidth]{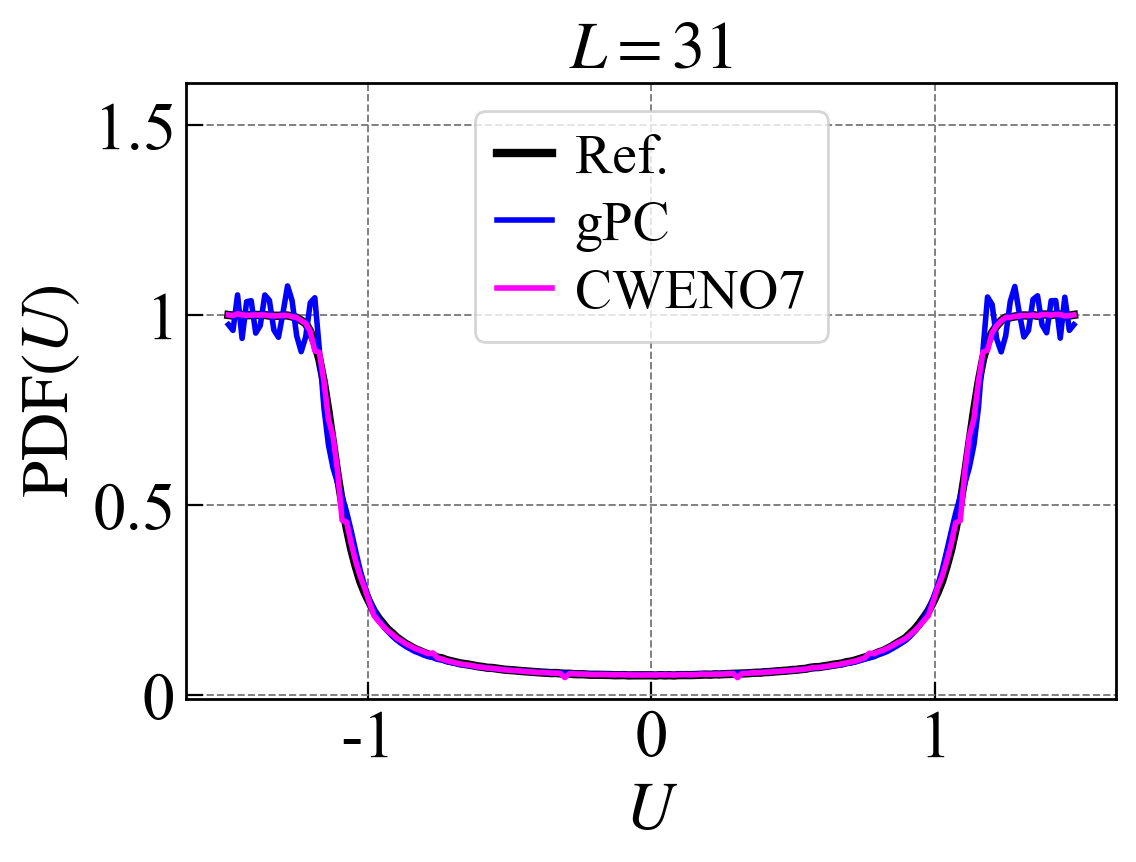}\hspace*{0.3cm}
            \includegraphics[trim=0.2cm 0.3cm 0.2cm 0.2cm, clip, width=0.31\textwidth]{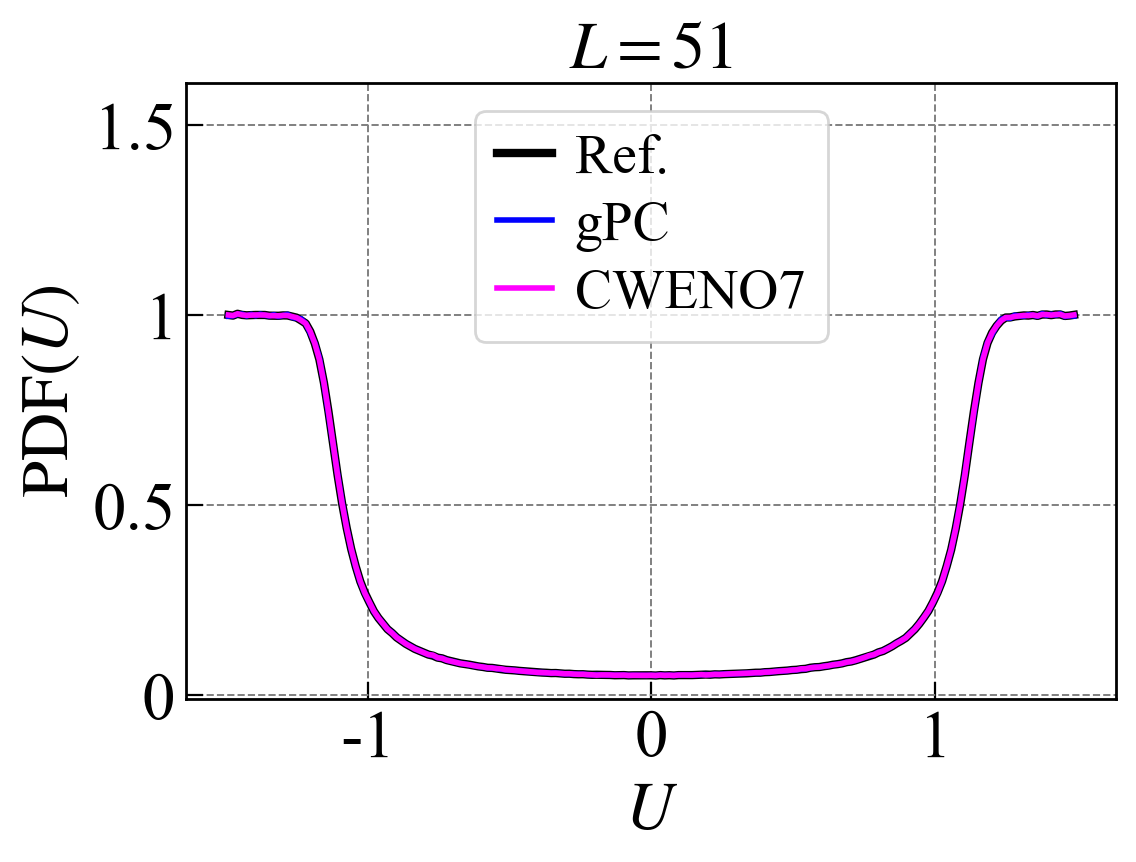}}
\caption{\sf Example 2: Estimated PDFs for the gPC and CWENO7 interpolations $\widetilde U$ together with the reference PDF, reconstructed
from $U$ for $L=21$ (left), $31$ (middle), and $51$ (right).\label{fig310}}
\end{figure}

\begin{figure}[ht!]
\centerline{\includegraphics[trim=0.2cm 0.3cm 0.2cm 0.2cm, clip, width=0.35\textwidth]{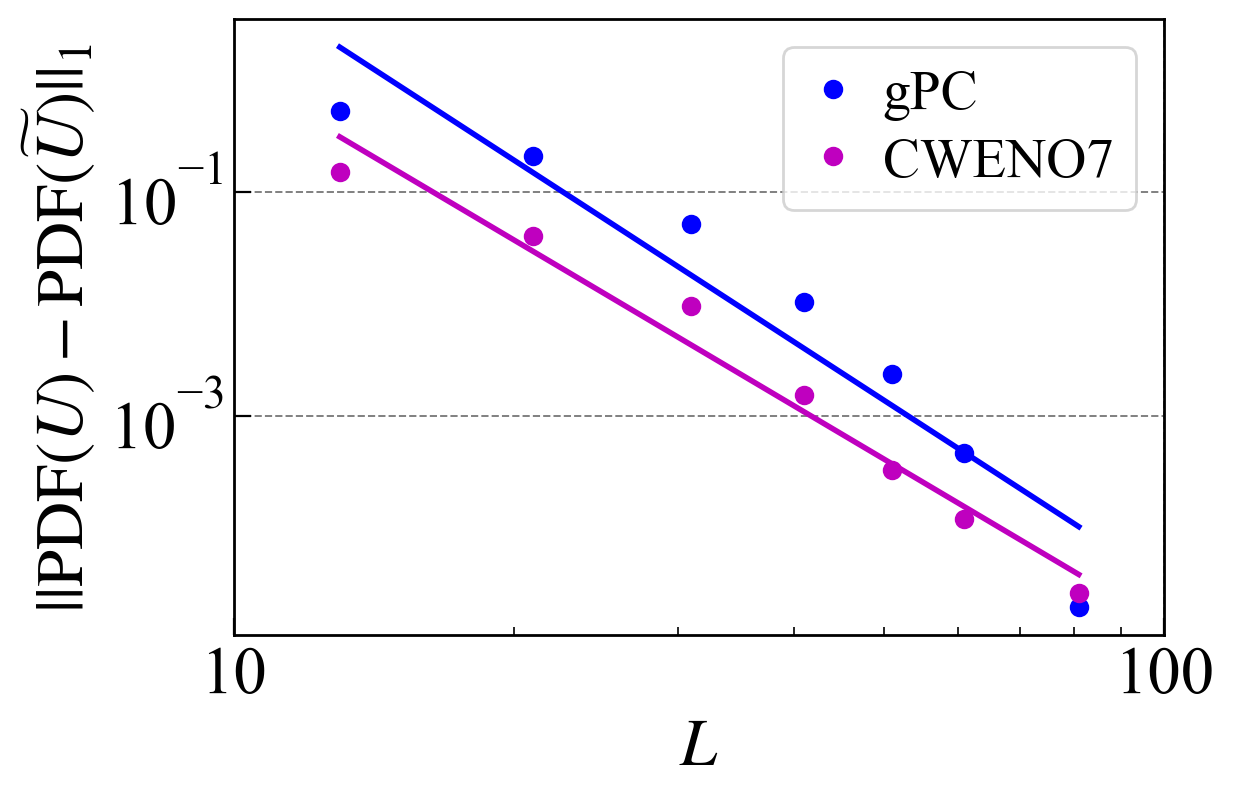}}
\caption{\sf Example 2: $L^1$-error for the PDFs as functions of $L$ and the corresponding power-law fits (solid lines).\label{fig311}}
\end{figure}

\smallskip
\underline{\it Convergence of PDF}. Next, we illustrate the convergence of the surrogate-based PDF approximations, which are shown in
\fref{fig310} for $L=21$, $31$, and $51$. Though we have observed convergence for both $\mu$ and $\sigma$, the PDF reconstructed using the
gPC-based surrogate model contained large oscillations when $L=21$. These oscillations decay when $L$ increases and disappear when $L=51$.
At the same time, the PDF reconstructed using the CWENO7-based surrogate model is oscillation-free, but as in Example 1, it contains 
discrepancies when $L=21$. \fref{fig311} displays the $L^1$-errors in PDFs as functions of $L$ along with the corresponding power-law fits
with the exponents $5.3$ (gPC) and $4.9$ (CWENO7). While both the gPC- and CWENO7-based surrogate models exhibit convergence rates
comparable to those observed in Example 1, the error magnitude for the gPC expansion is notably higher than for the CWENO7 interpolation 
due to the oscillations.

\subsubsection*{Example 3---Smooth Function \eref{3.3}}
In this example, we consider two differently distributed random variables $\xi\sim{\cal U}(-1,1)$ and $\eta\sim{\cal N}(0,1)$, and use
uniform collocation points for $\xi_\ell\in[-1,1]$ and $\eta_m\in[-6,6]$.

\begin{figure}[ht!]
\centerline{\includegraphics[trim=0.2cm 0.3cm 0.2cm 0.2cm, clip, width=0.31\textwidth]{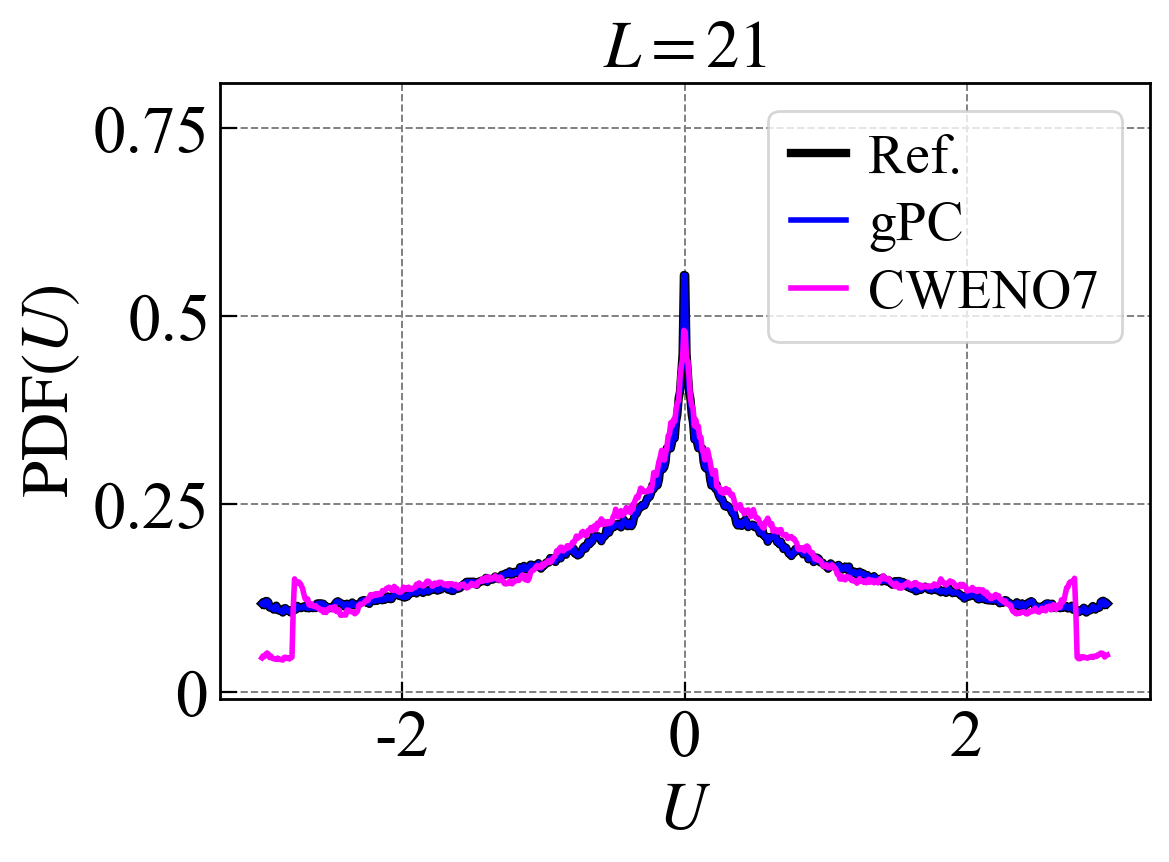}\hspace*{0.3cm}
            \includegraphics[trim=0.2cm 0.3cm 0.2cm 0.2cm, clip, width=0.31\textwidth]{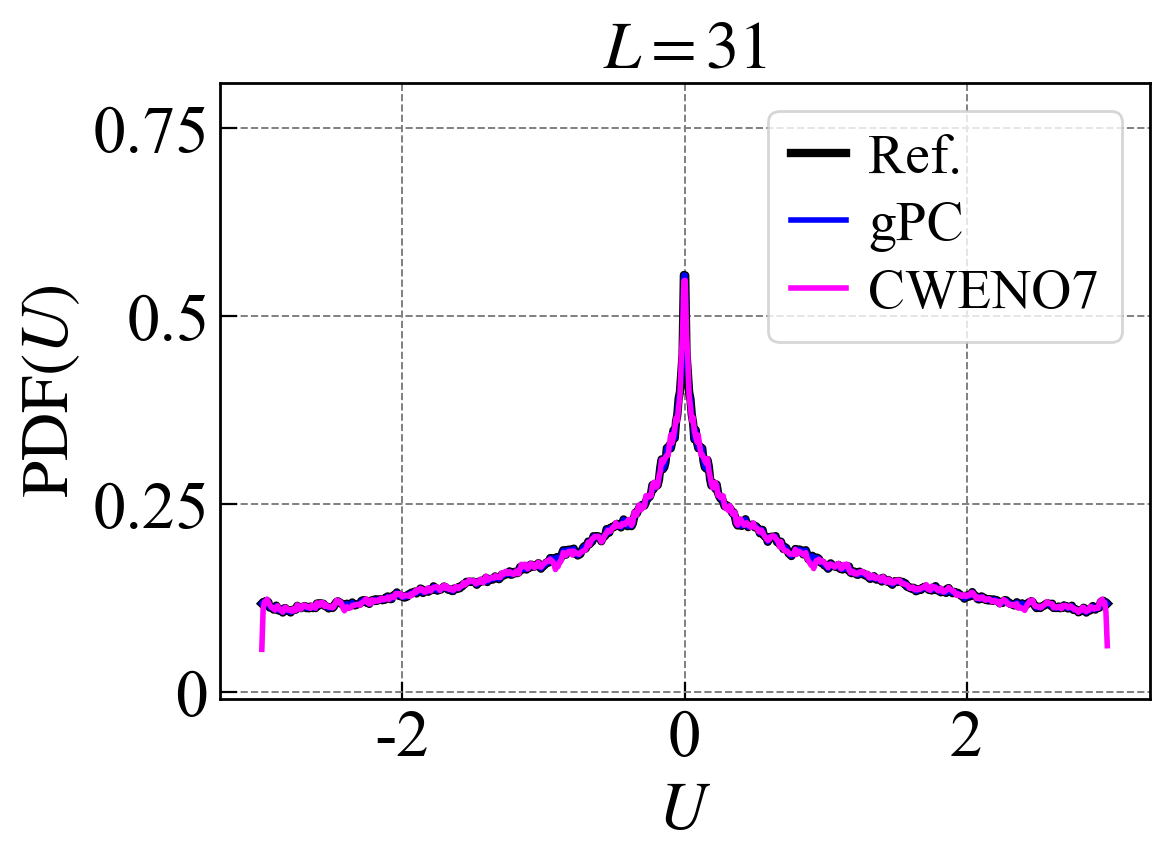}\hspace*{0.3cm}
            \includegraphics[trim=0.2cm 0.3cm 0.2cm 0.2cm, clip, width=0.31\textwidth]{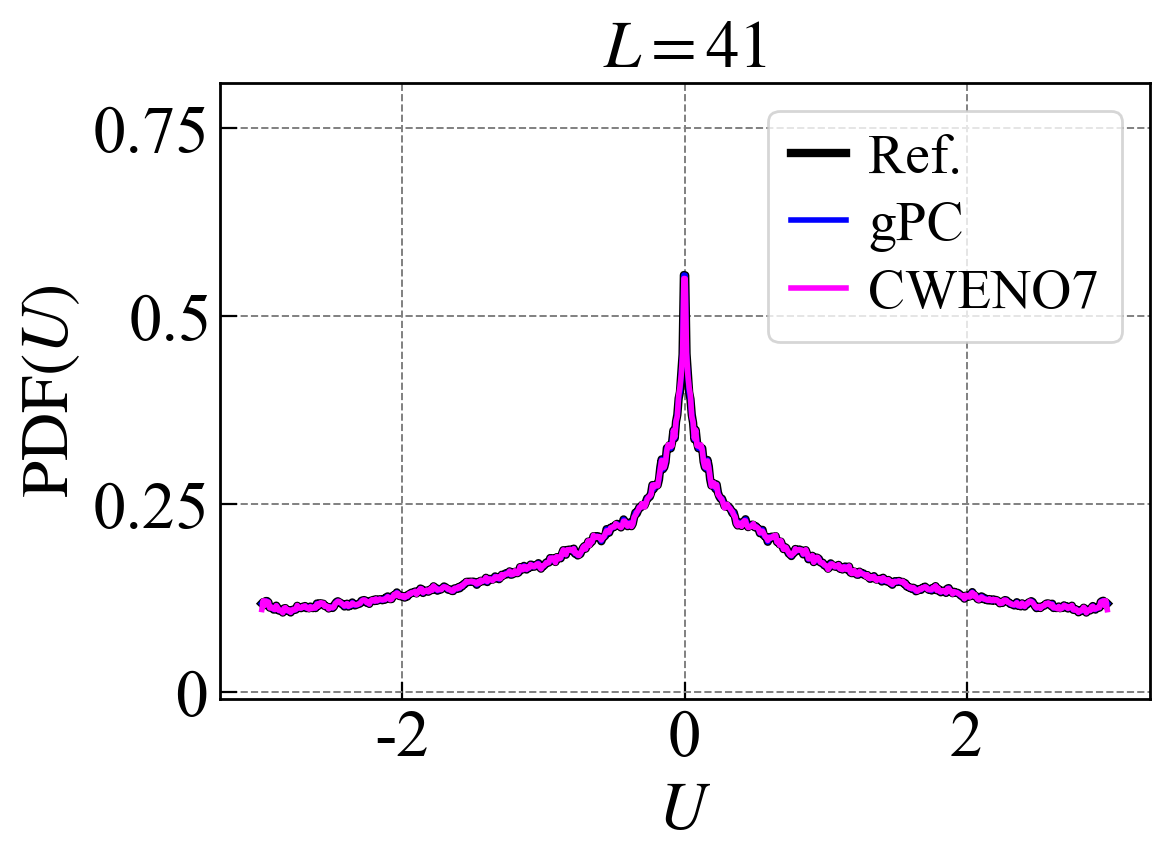}}
\caption{\sf Example 3: Estimated PDFs for the gPC and CWENO7 interpolations $\widetilde U$ together with the reference PDF, reconstructed
from $U$ for $L=M=21$ (left), $31$ (middle), and $41$ (right).\label{fig312}}
\end{figure}

We first show in \fref{fig312} the surrogate-based PDF approximations for $L=M=21$, $31$, and $41$ (in principle, one can take different
values for $L$ and $M$, but this is not essential for the conducted convergence study). As in the case of one random variable,
considered in Examples 1 and 2, one can observe the discrepancies in the CWENO7-based PDFs computed with $L=M=21$ and $31$. As before, 
these discrepancies disappear when the number of collocation points increases.

\begin{figure}[ht!]
\centerline{\includegraphics[trim=0.2cm 0.3cm 0.2cm 0.2cm, clip, width=0.35\textwidth]{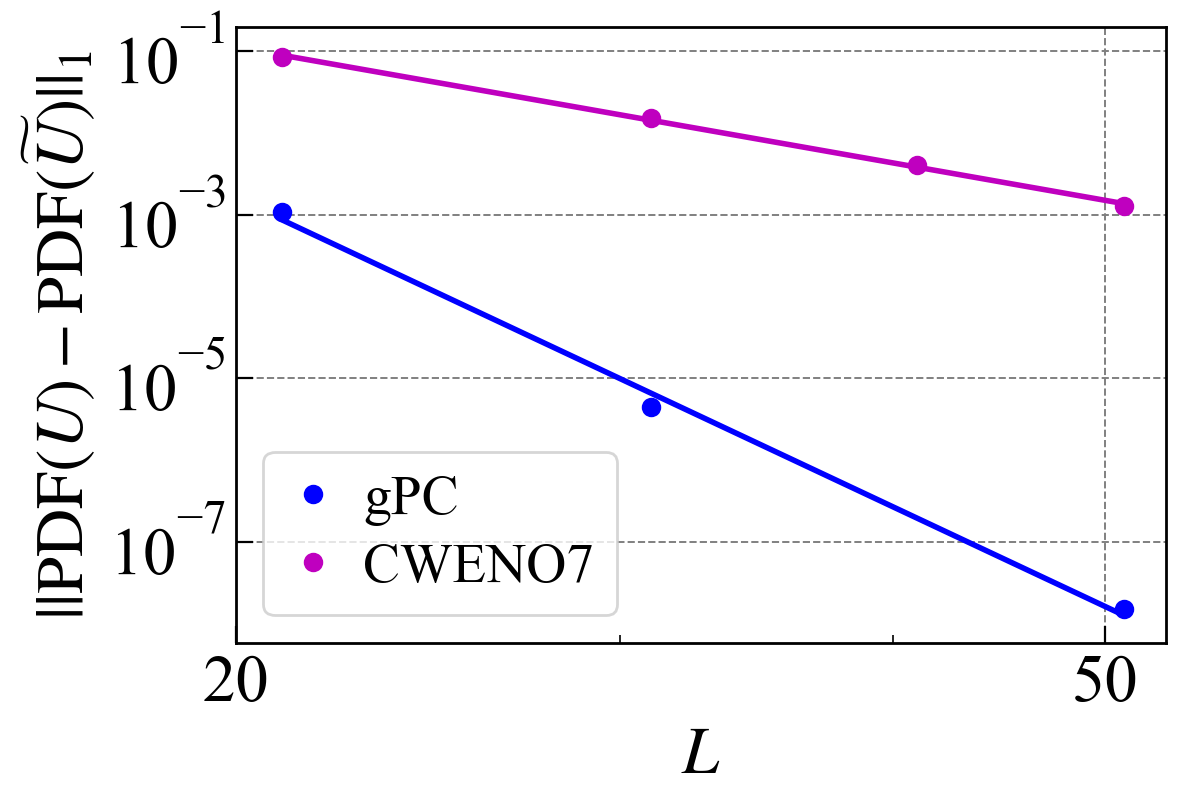}}
\caption{\sf Example 3: $L^1$-error for the PDFs as functions of $L$ and the corresponding power-law fits (solid lines).\label{fig313}}
\end{figure}

\fref{fig313} shows the $L^1$-errors in PDFs as functions of $L$ (for $L=M=21$, $31$, $41$, and $51$) along with the corresponding
power-law fits with the exponents $11.6$ (gPC) and $4.6$ (CWENO7). As expected, in this example, the gPC-based surrogate model outperforms
the CWENO7-based one due to the smoothness of $U$.

\subsubsection*{Example 4---Nonsmooth Function \eref{3.4}}
In this example, we consider a discontinuous function $U$ of a random variable $\xi\sim{\cal U}[-1,1]$, for which the CWENO7-based
surrogate model is expected to be superior to the gPC-based one.

\begin{figure}[ht!]
\centerline{\includegraphics[trim=0.2cm 0.3cm 0.2cm 0.2cm, clip, width=0.31\textwidth]{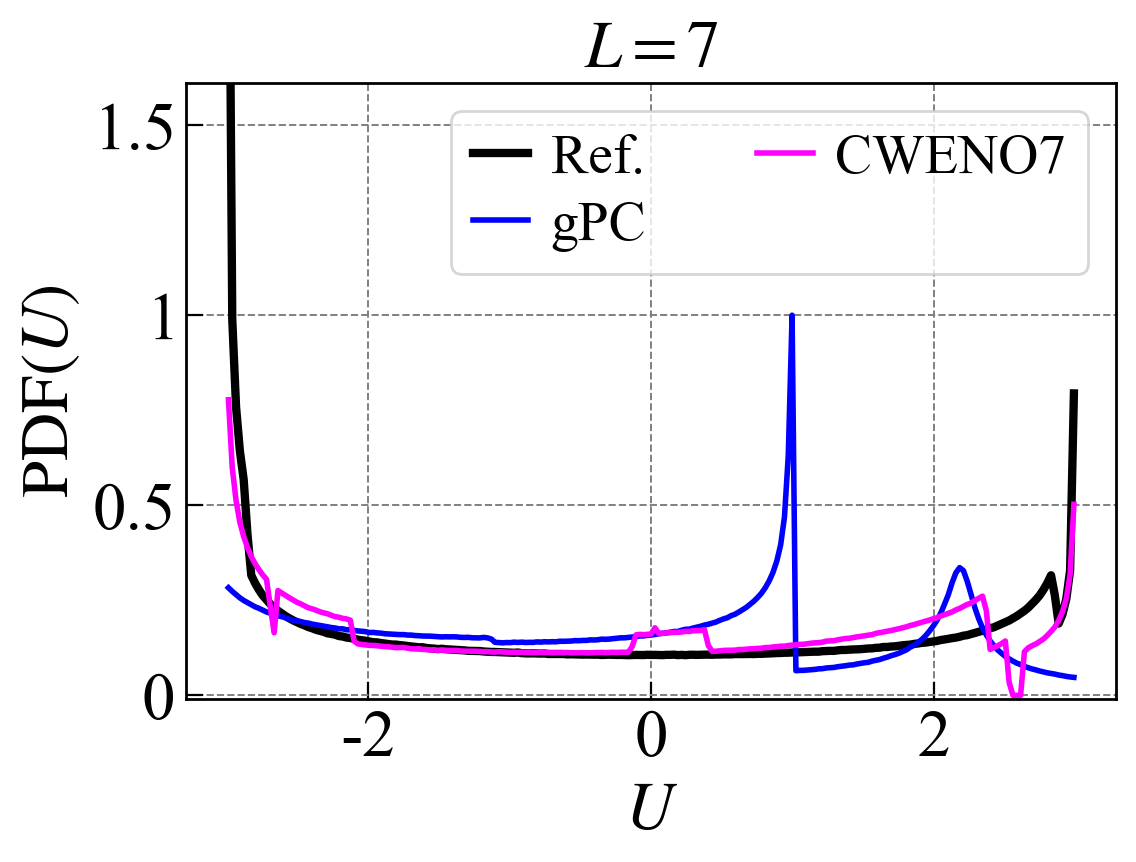}\hspace*{0.3cm}
            \includegraphics[trim=0.2cm 0.3cm 0.2cm 0.2cm, clip, width=0.31\textwidth]{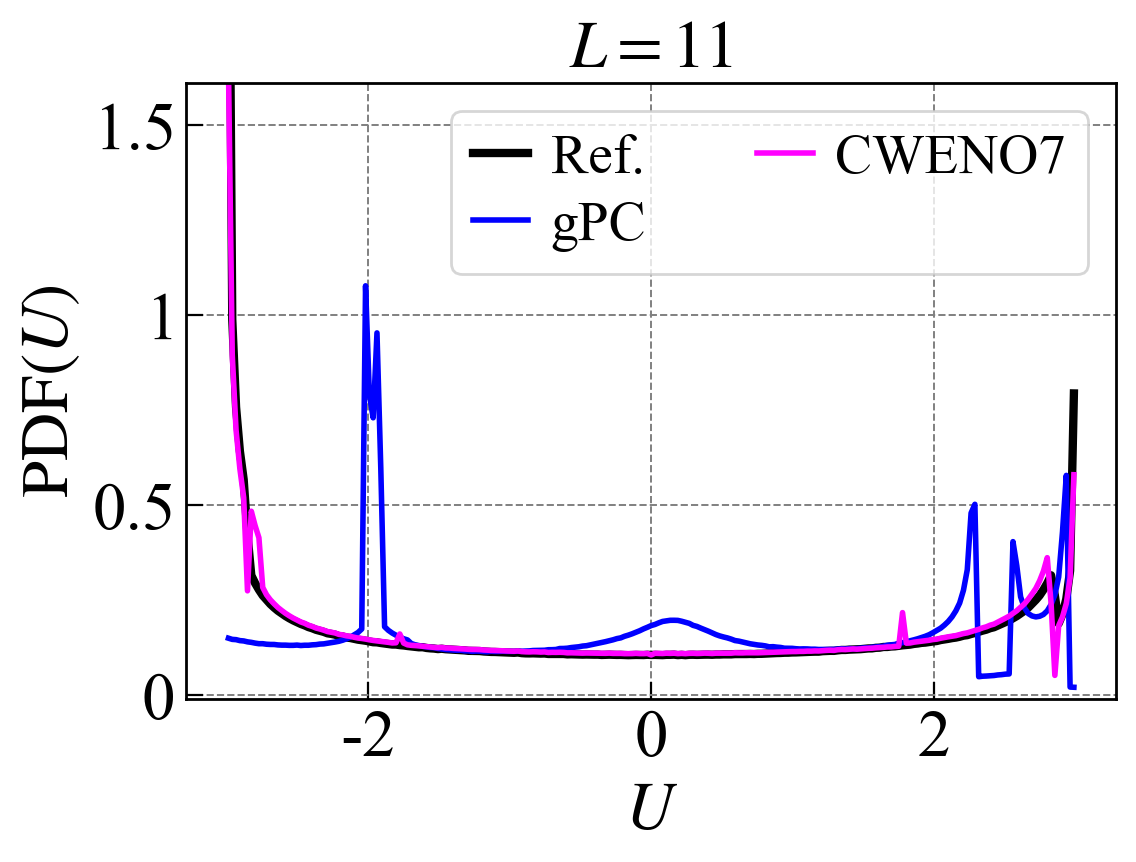}\hspace*{0.3cm}
            \includegraphics[trim=0.2cm 0.3cm 0.2cm 0.2cm, clip, width=0.31\textwidth]{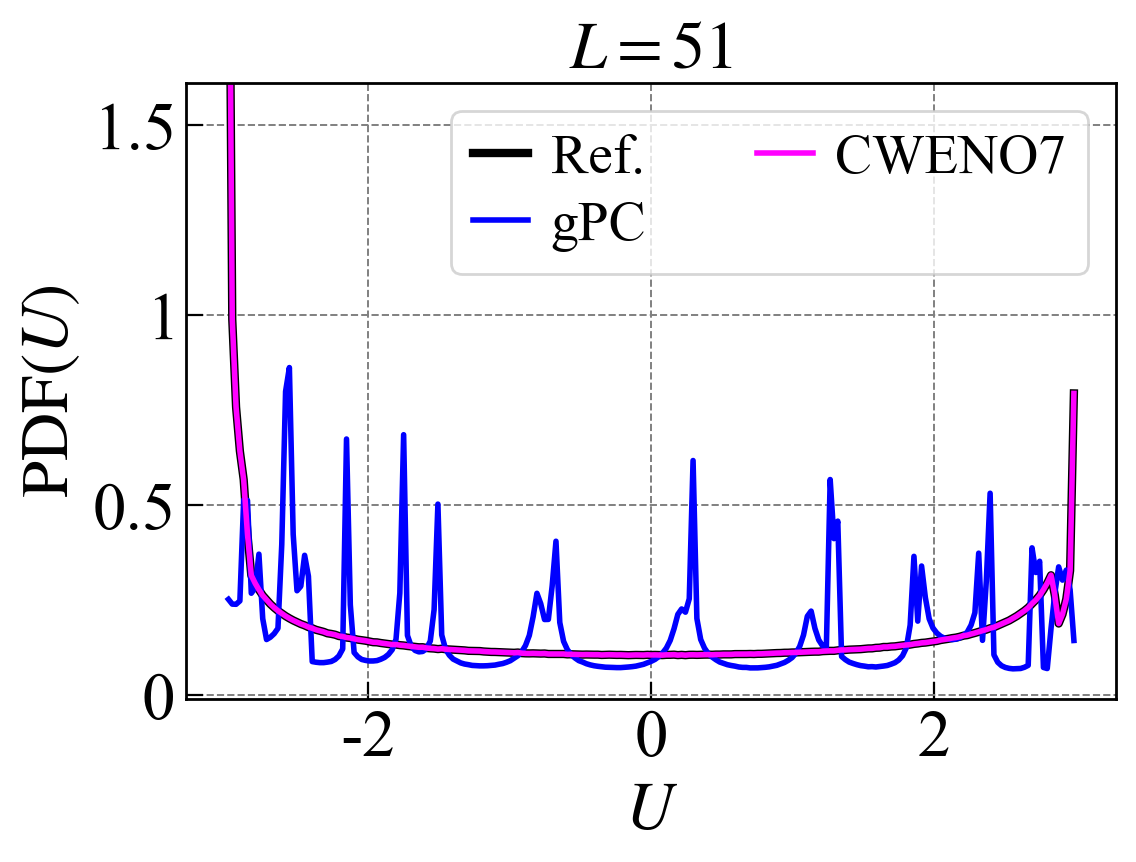}}
\caption{\sf Example 4: Estimated PDFs for the gPC and CWENO7 interpolations $\widetilde U$ together with the reference PDF, reconstructed
from $U$ for $L=7$ (left), $11$ (middle), and $51$ (right).\label{fig314}}
\end{figure}

We first show in \fref{fig314} the surrogate-based PDF approximations for $L=7$, $11$, and $51$. As one can see that the gPC expansion
yields oscillatory behavior and fails to restore the PDF even when the number of modes increases (in fact, the oscillations spread all over
the $U$-domain when $L$ is large). This result is expected since spectral-type approximations are known to suffer from the Gibbs phenomenon
and thus do not apply to discontinuous functions unless an appropriate filtering strategy is implemented. The CWENO7 interpolation, on the
other hand, provides a robust and accurate approximation of the PDF as $L$ increases.

\begin{figure}[ht!]
\centerline{\includegraphics[trim=0.2cm 0.3cm 0.2cm 0.2cm, clip, width=0.35\textwidth]{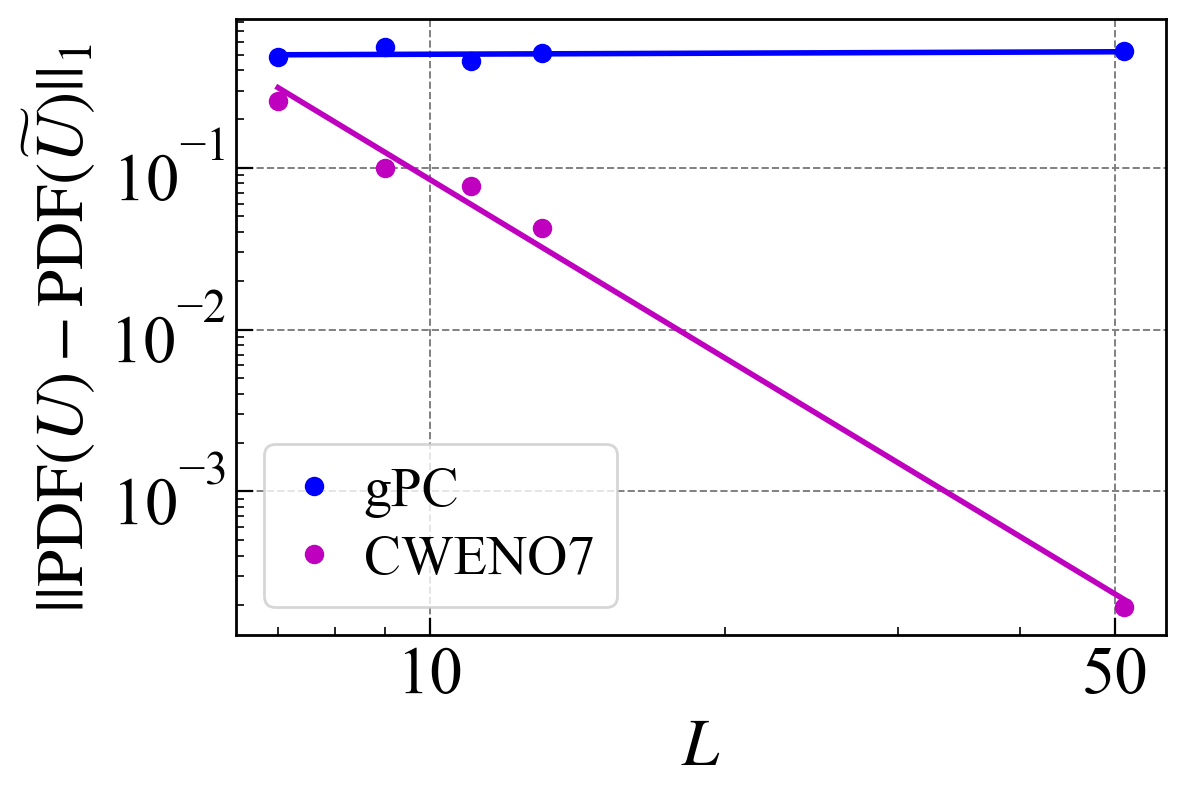}}
\caption{\sf Example 4: $L^1$-error for the PDFs as functions of $L$ and the corresponding power-law fits (solid lines).\label{fig315}}
\end{figure}

\fref{fig315} shows the $L^1$-errors in PDFs as functions of $L$ (for $L=7$, $9$, $11$, $13$, and $51$) along with the corresponding
power-law fits. One can observe that no convergence is achieved when the gPC expansion is used, while the use of the CWENO7 interpolation
yields a convergence rate of about $3.4$.

\subsubsection*{Example 5---Nonsmooth Function \eref{3.5}}
In this example, we consider a discontinuous function $U$ of two differently distributed random variables $\xi\sim{\cal U}(-1,1)$ and
$\eta\sim{\cal N}(0,1)$ and, as in Example 3, we take the equally spaced collocation points for $\xi\in[-1,1]$ and $\eta\in[-6,6]$.

\begin{figure}[ht!]
\centerline{\includegraphics[trim=0.2cm 0.3cm 0.2cm 0.2cm, clip, width=0.31\textwidth]{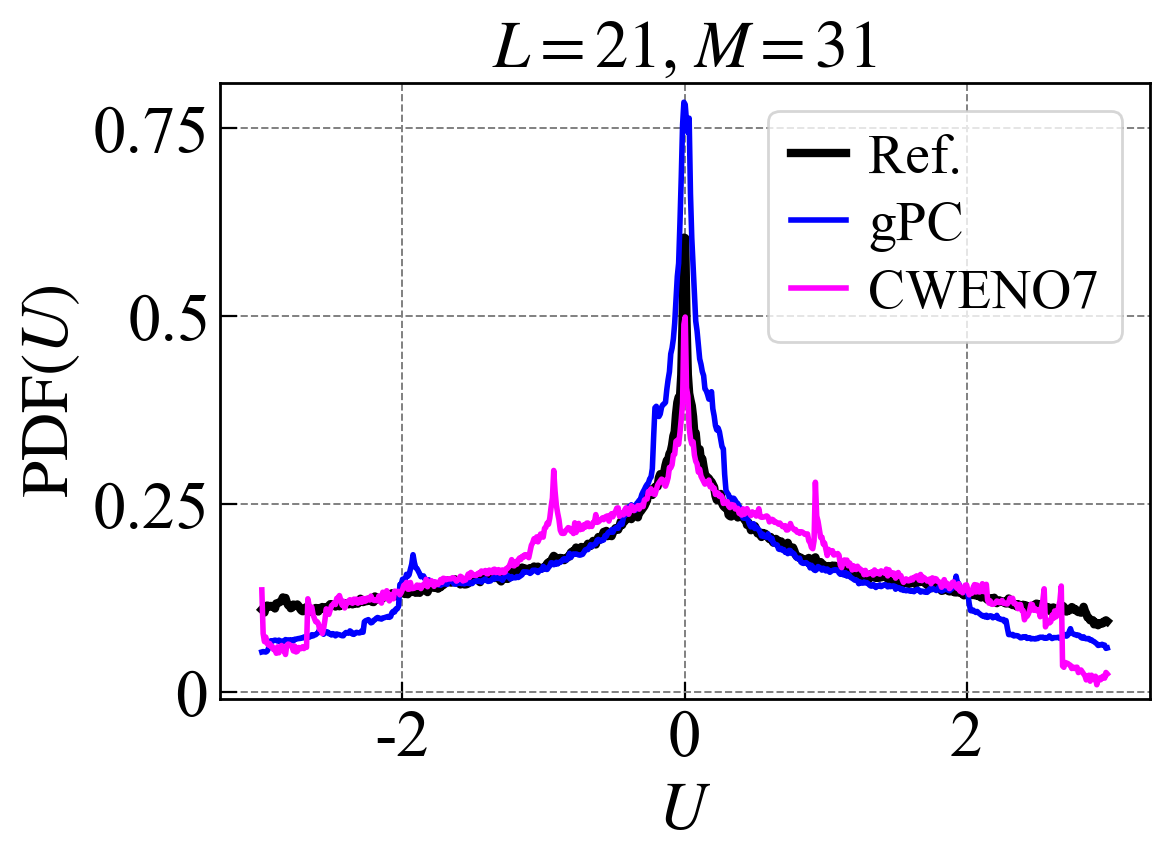}\hspace*{0.3cm}
            \includegraphics[trim=0.2cm 0.3cm 0.2cm 0.2cm, clip, width=0.31\textwidth]{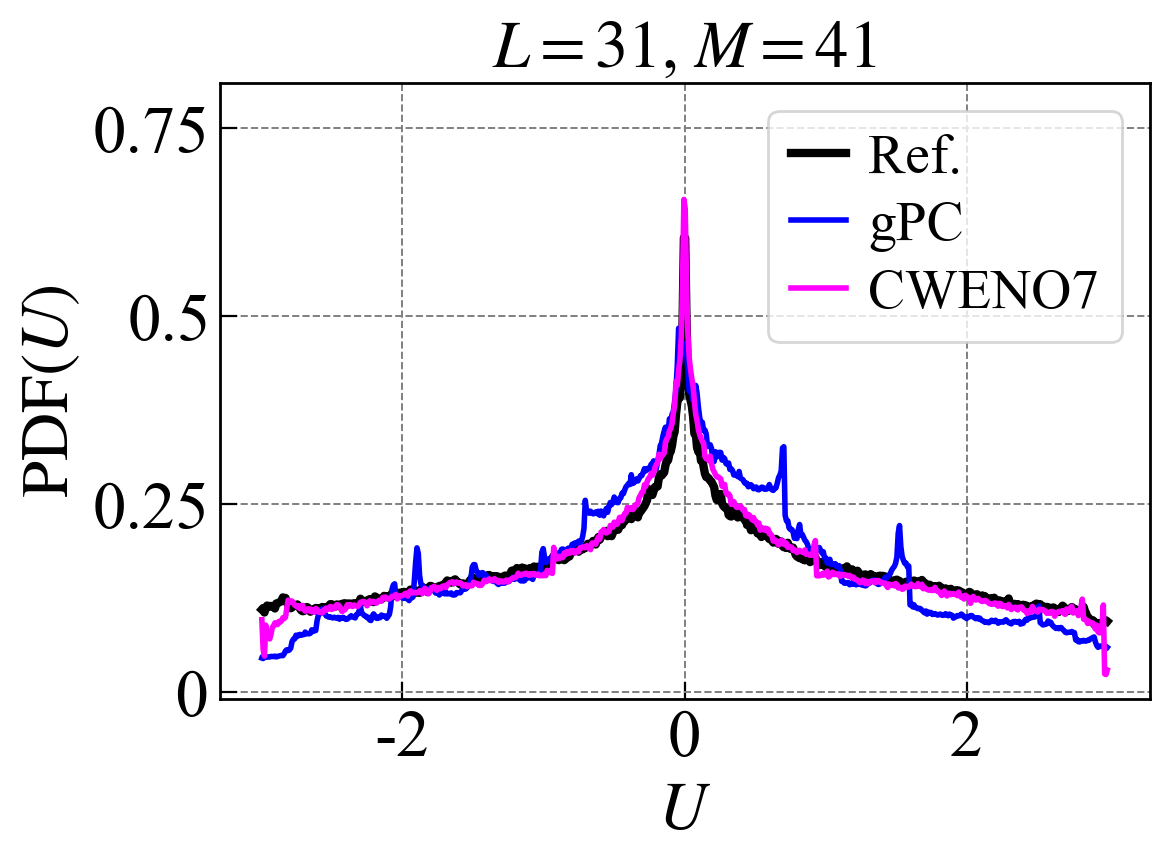}\hspace*{0.3cm}
            \includegraphics[trim=0.2cm 0.3cm 0.2cm 0.2cm, clip, width=0.31\textwidth]{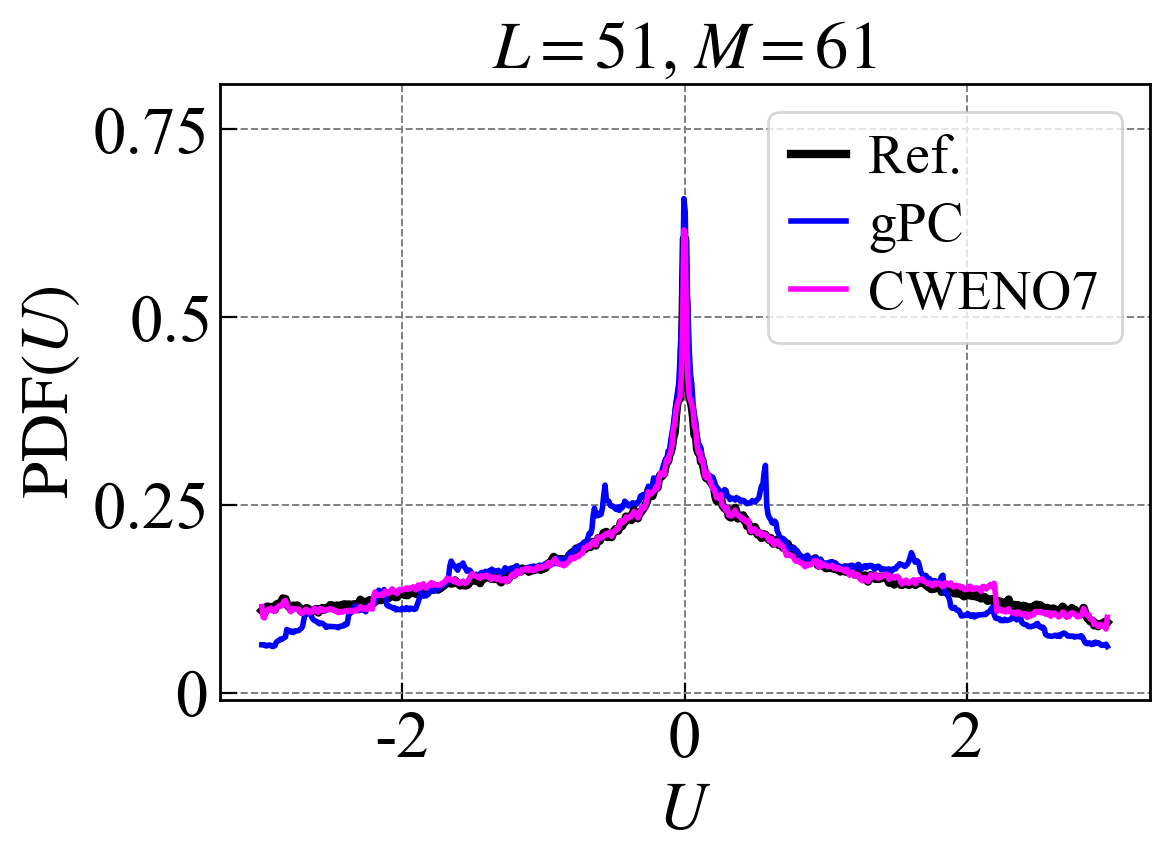}}
\caption{\sf Example 5: Estimated PDFs for the gPC and CWENO7 interpolations $\widetilde U$ together with the reference PDF, reconstructed
from $U$ for $(L,M)=(21,31)$ (left), $(31,41)$ (middle), and $(51,61)$ (right).\label{fig316}}
\end{figure}

\begin{figure}[ht!]
\centerline{\includegraphics[trim=0.2cm 0.3cm 0.2cm 0.2cm, clip, width=0.35\textwidth]{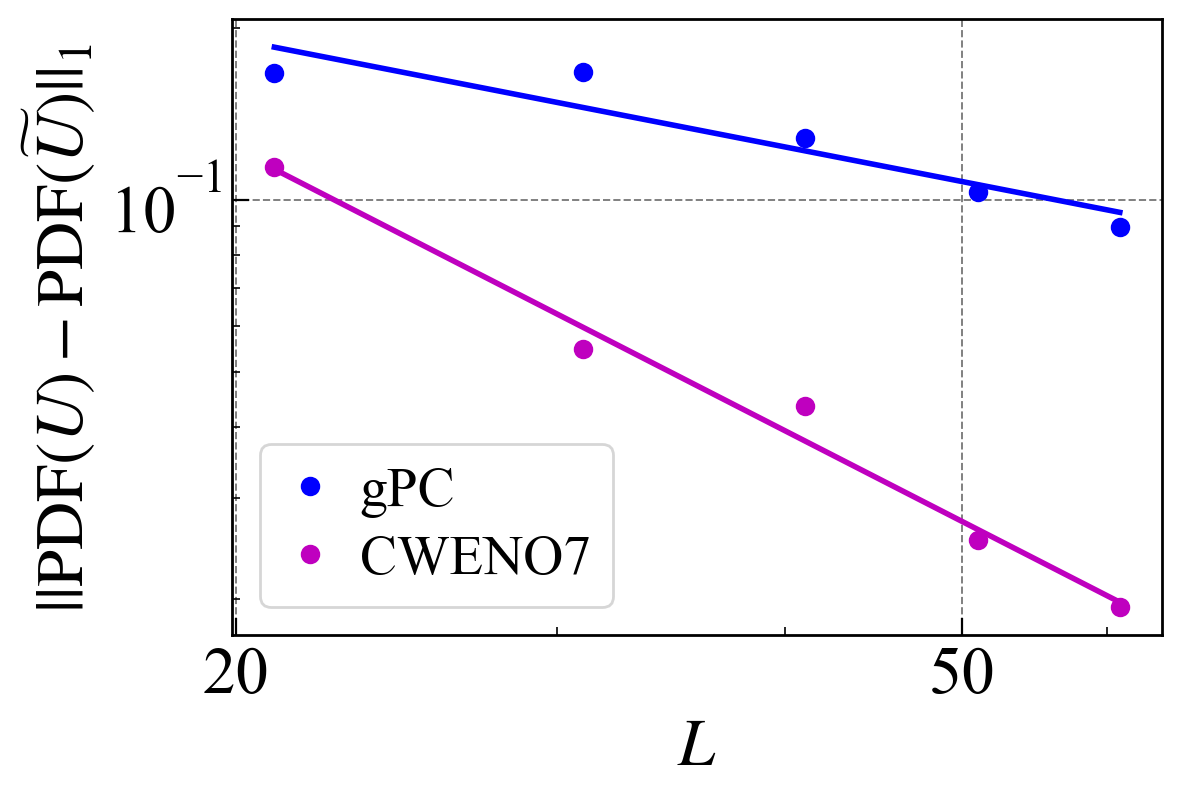}}
\caption{\sf Example 5: $L^1$-error for the PDFs as functions of $L$ and the corresponding power-law fits (solid lines).\label{fig317}}
\end{figure}

We now take $L\ne M$ and a particular choice is $(L,M)=(11,21)$, $(21,31)$, $(31,41)$, $(41,51)$, and $(51,61)$. The PDF approximations for
three of these pairs of values are shown in \fref{fig316}. While both of the studied surrogate models exhibit an oscillatory behavior for a
small number of collocation points, the CWENO7 interpolation successfully restores the PDF, provided a sufficiently large number of
collocation points is used. The gPC expansion, on the contrary, still suffers from the Gibbs phenomenon and its observed $L^1$ convergence
rate is very small ($0.6$ vs. $1.5$, observed for the CWENO7 interpolation); see \fref{fig317}, in which we plot the $L^1$-error as a
function of $L$.

\subsubsection*{Example 6---Dam Break, One Random Variable}
In this example, we take a scalar random variable $\xi\sim{\cal U}[-1,1]$ and obtain the discrete function $w(x_j,0.8;\xi)$ by numerically
solving \eref{3.6} for the deterministic initial data,
\begin{equation}
w(x,0;\xi)=\left\{\begin{aligned}&1,&&x<0,\\&0.5,&&x>0,\end{aligned}\right.\qquad u(x,0;\xi)\equiv0,
\label{3.7}
\end{equation}
and stochastic bottom topography,
\begin{equation*}
Z(x;\xi)=\begin{cases}0.125\xi+0.125(\cos(5\pi x)+2), &|x|<0.2,\\0.125\xi+0.125,&\mbox{otherwise},\end{cases}
\end{equation*}
prescribed in the spatial computational domain $x\in[-1,1]$ subject to the free boundary conditions. We take $800$ spatial finite-volume
cells and $L=32$ collocation points.

\begin{figure}[ht!]
\centerline{\includegraphics[trim=0.2cm 0.3cm 0.2cm 0.2cm, clip, width=0.31\textwidth]{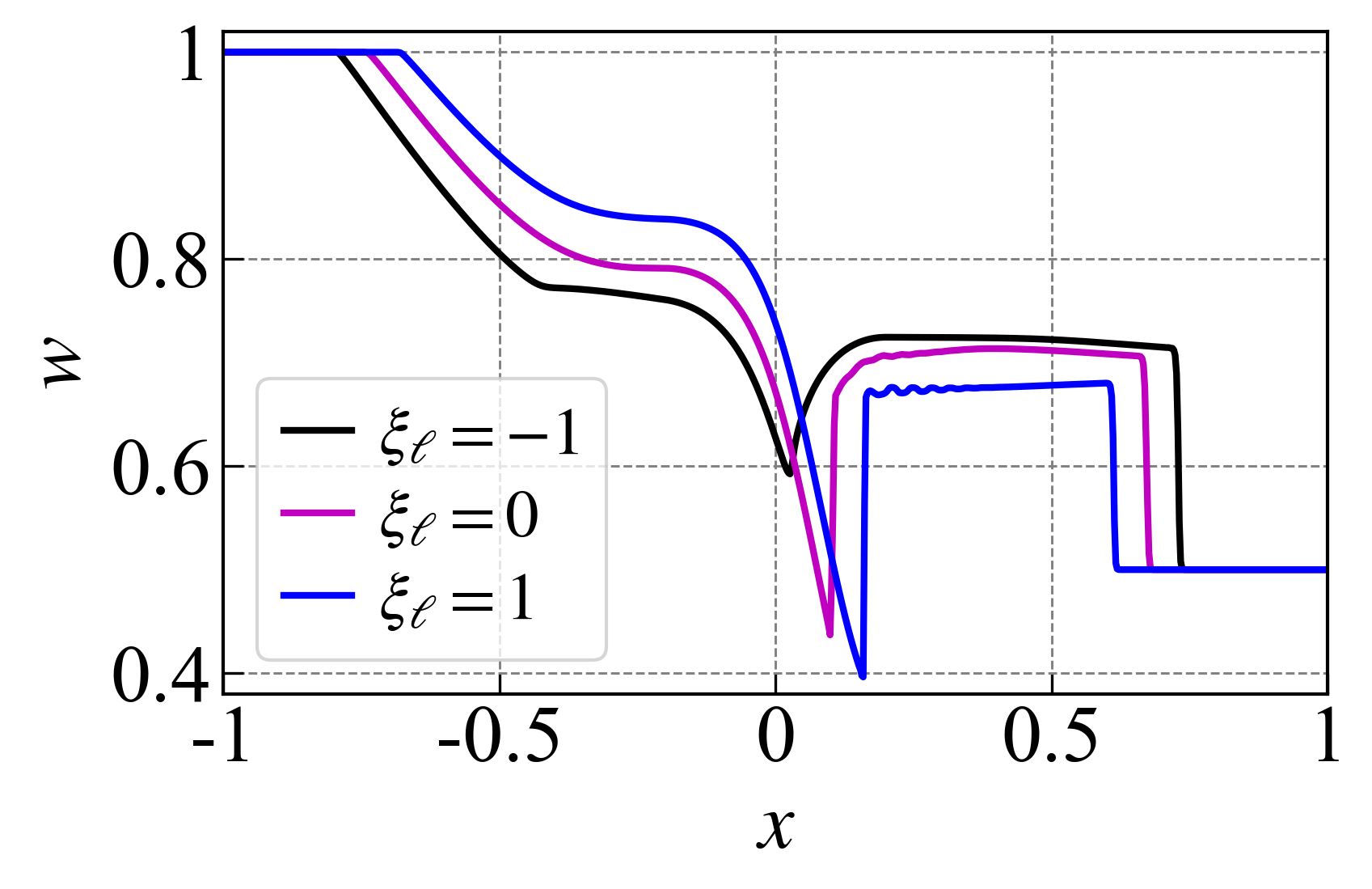}\hspace*{0.3cm}
            \includegraphics[trim=0.2cm 0.3cm 0.2cm 0.2cm, clip, width=0.318\textwidth]{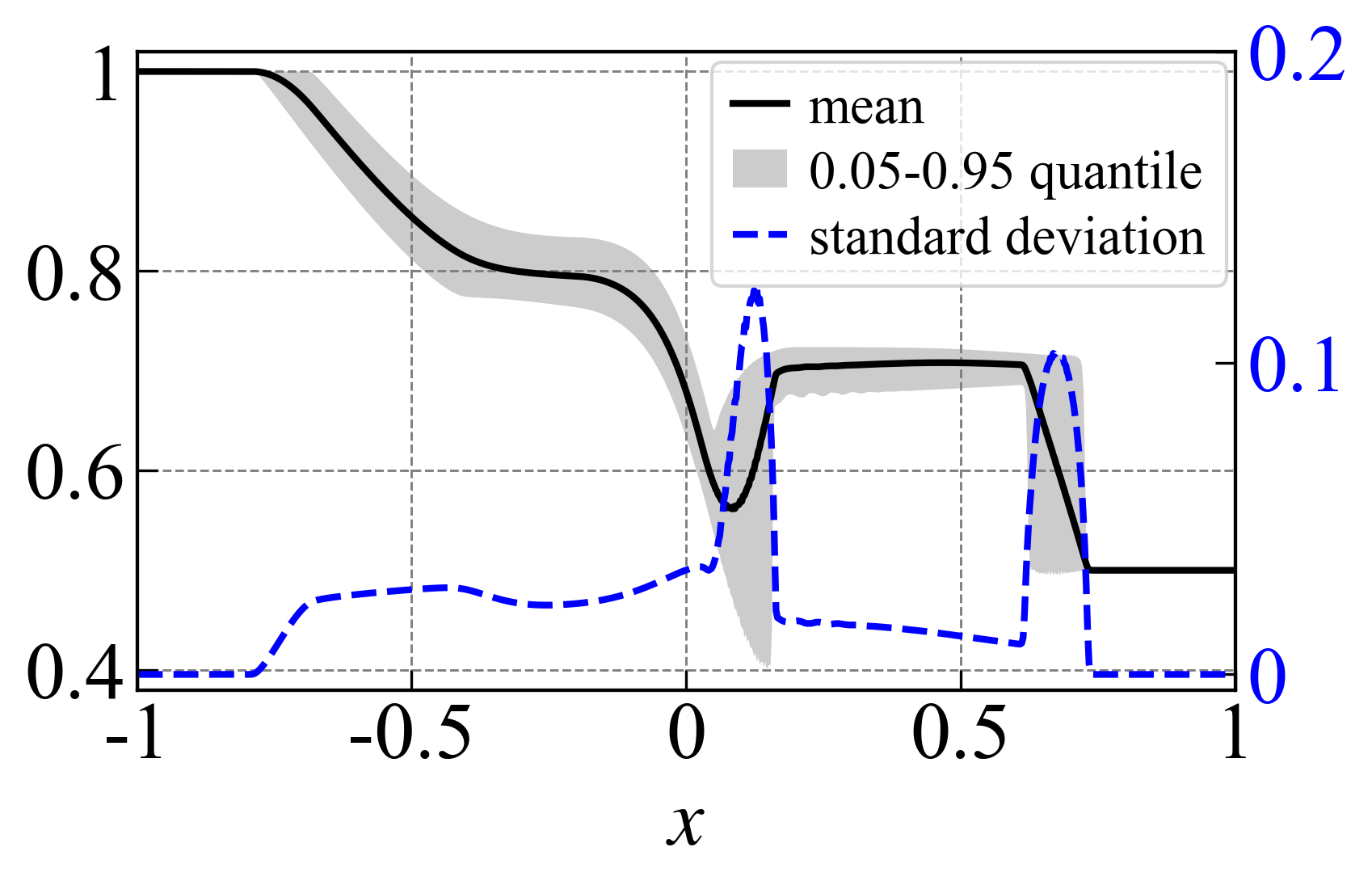}\hspace*{0.3cm}
            \includegraphics[trim=0.2cm 0.3cm 0.2cm 0.2cm, clip, width=0.318\textwidth]{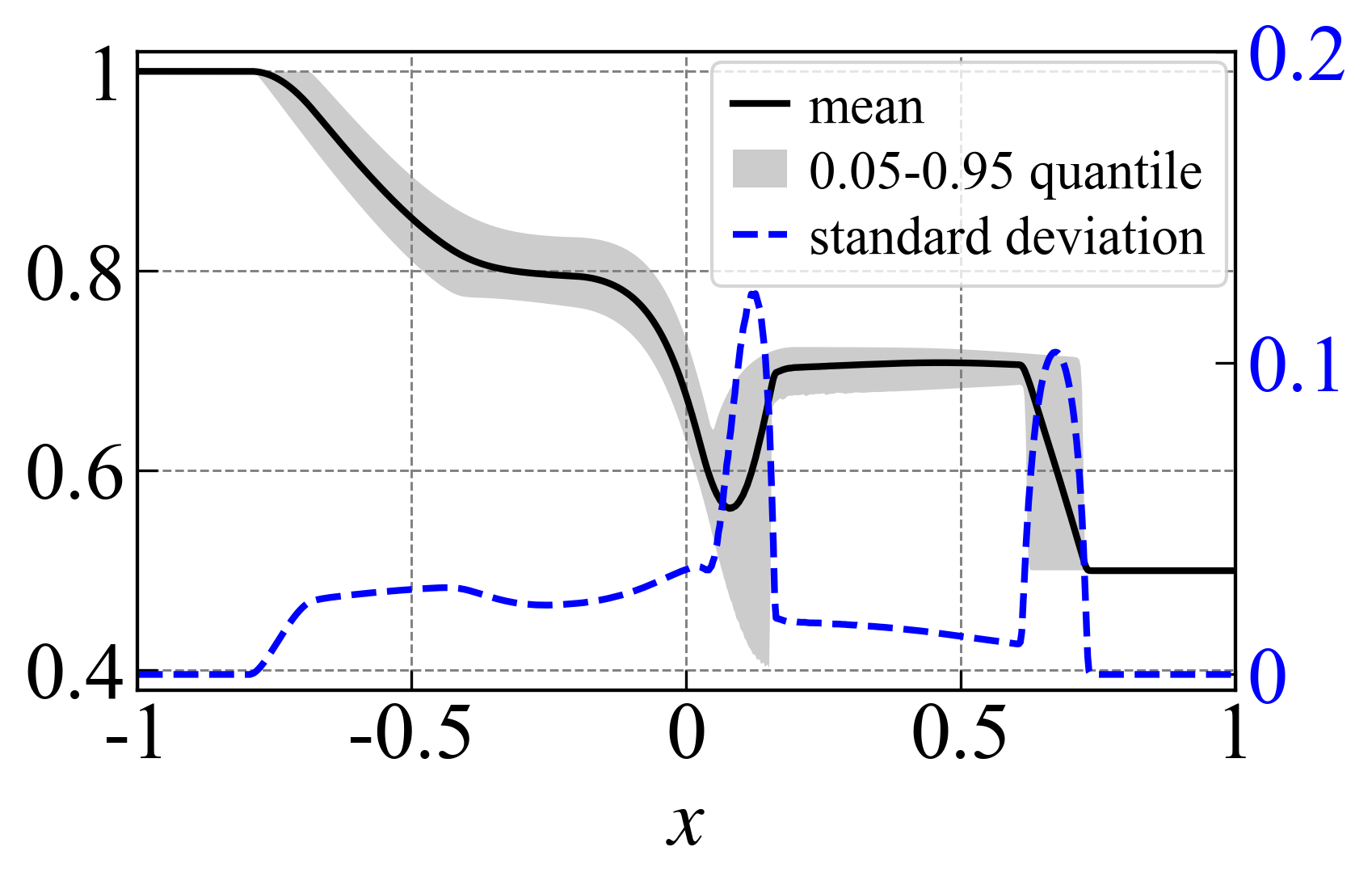}}
\caption{\sf Example 6: Computed water surface $w(x,0.8;\xi_\ell)$ for $\xi_\ell=-1$, $0$, and $1$ (left) along with the mean and standard
deviation obtained using the gPC expansion (middle) and CWENO7 interpolation (right).\label{fig318}}
\end{figure}

In \fref{fig318} (left), we plot the computed solution profiles $w(x,0.8;\xi_\ell)$ for three different values of the parameter $\xi_\ell$.
As one can see, for $\xi_\ell=-1$, the water surface over the bottom hump is smoothed out, while it contains a hydraulic jump for
$\xi_\ell=0$ and $1$ (the size of the jump increases when $\xi_\ell$ increases). The figure also displays the mean and standard deviation of
the water surface, obtained using both the gPC expansion (center) and the CWENO7 interpolation (right), demonstrating that the two methods
yield almost identical results.

\begin{figure}[ht!]
\centerline{\includegraphics[trim=0.2cm 0.3cm 0.2cm 0.2cm, clip, width=0.31\textwidth]{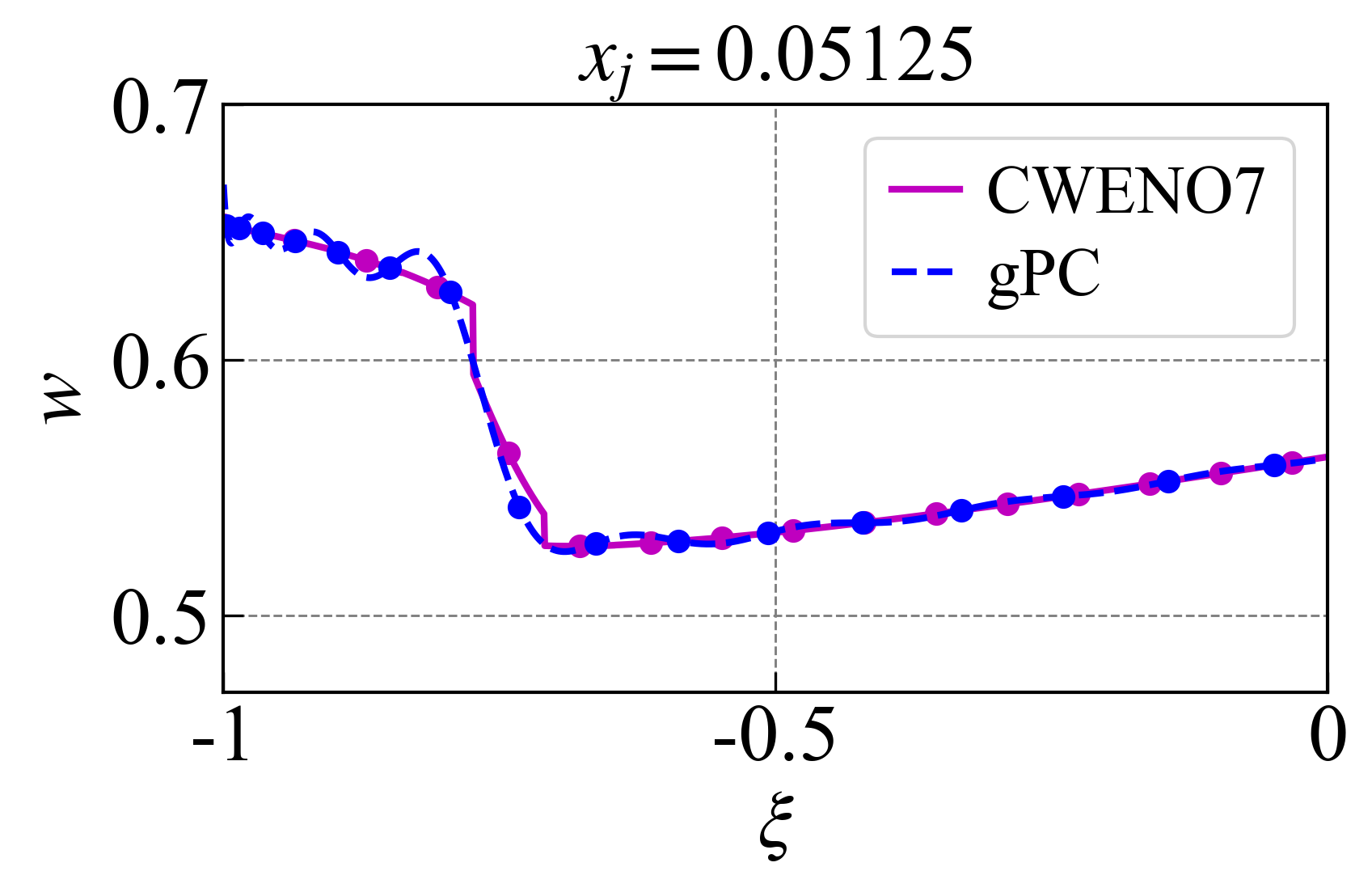}\hspace*{0.3cm}
            \includegraphics[trim=0.2cm 0.3cm 0.2cm 0.2cm, clip, width=0.31\textwidth]{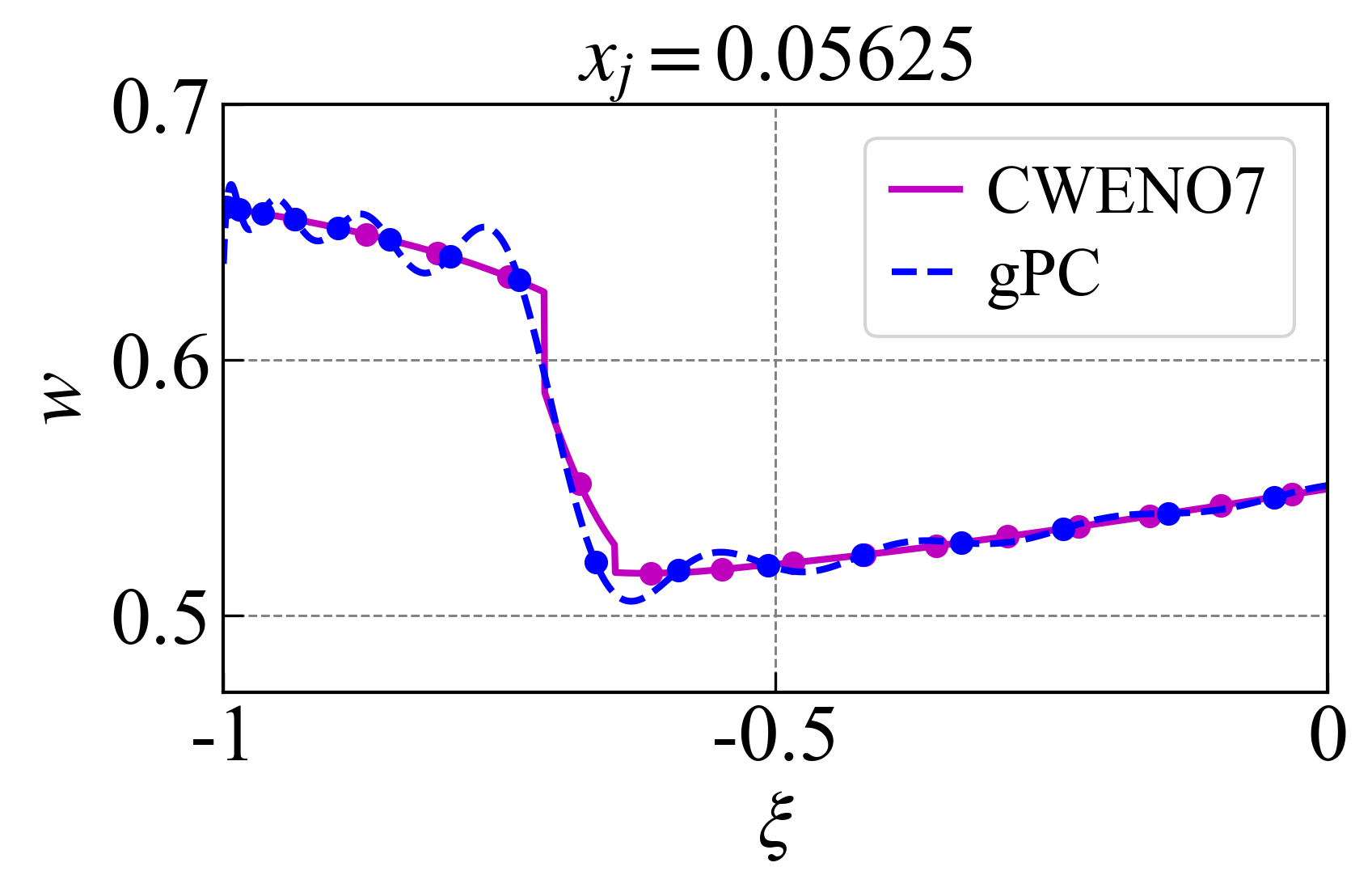}\hspace*{0.3cm}
            \includegraphics[trim=0.2cm 0.3cm 0.2cm 0.2cm, clip, width=0.31\textwidth]{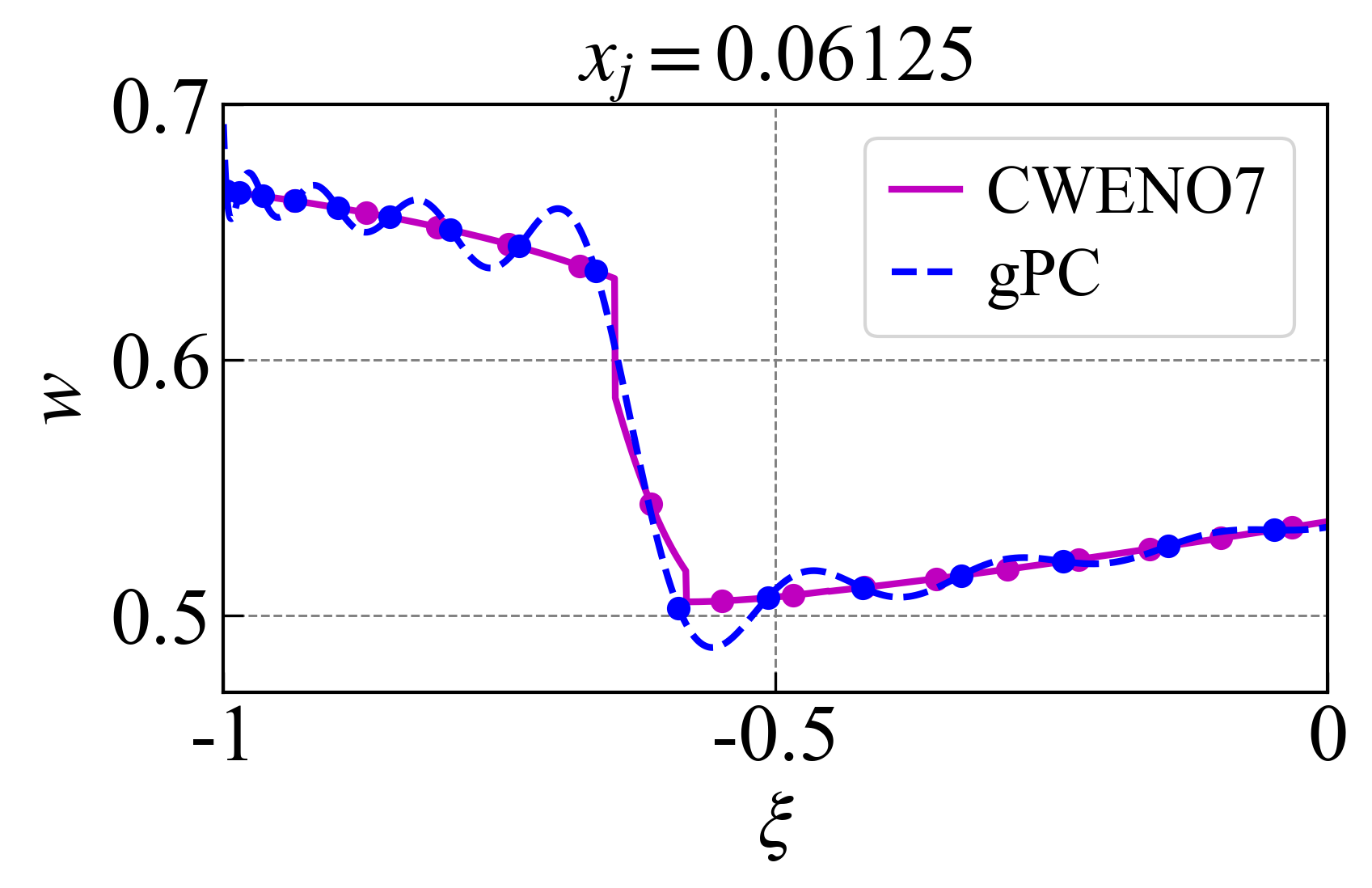}}
\vskip8pt
\centerline{\includegraphics[trim=0.2cm 0.3cm 0.2cm 0.2cm, clip, width=0.31\textwidth]{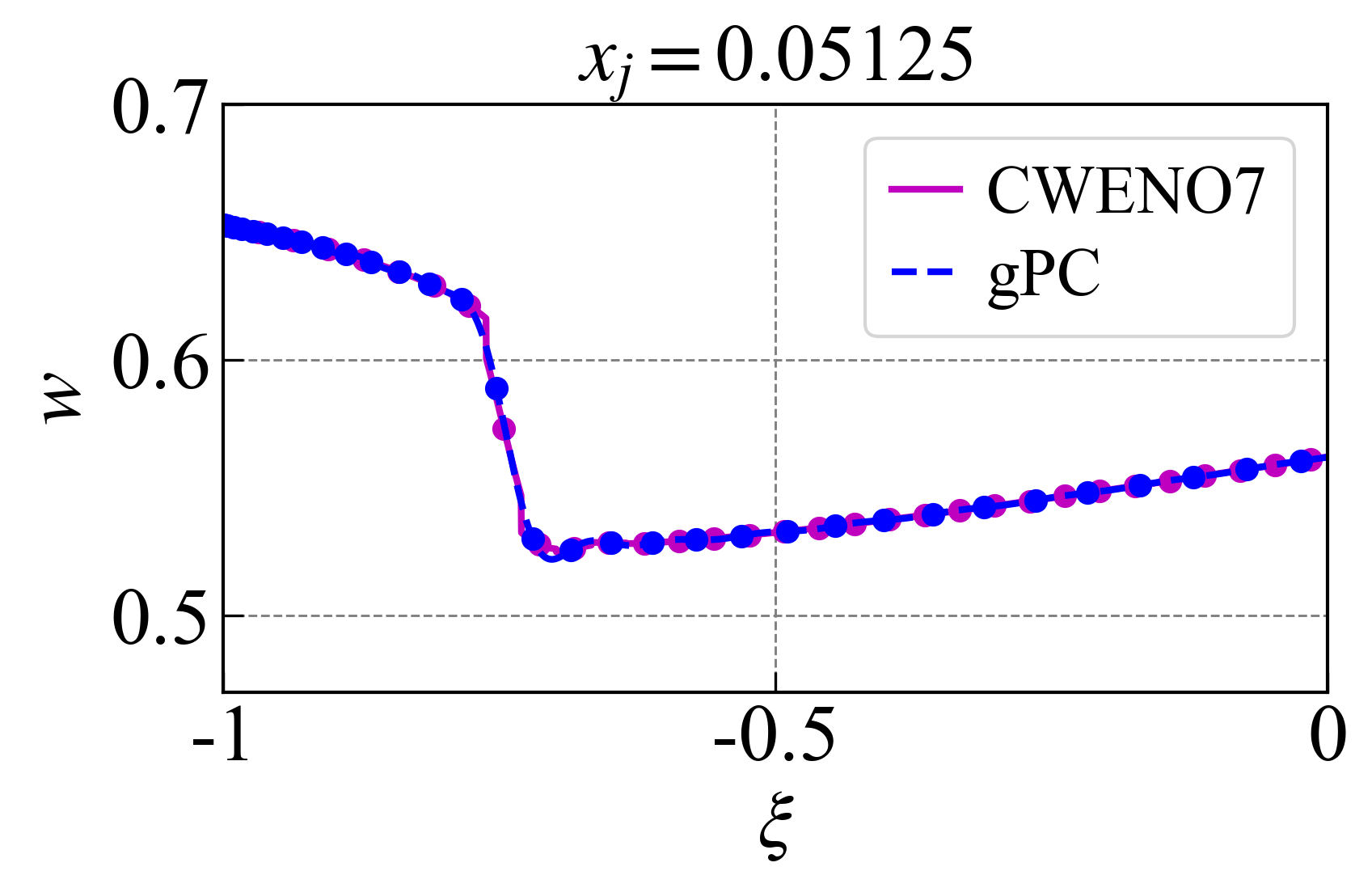}\hspace*{0.3cm}
            \includegraphics[trim=0.2cm 0.3cm 0.2cm 0.2cm, clip, width=0.31\textwidth]{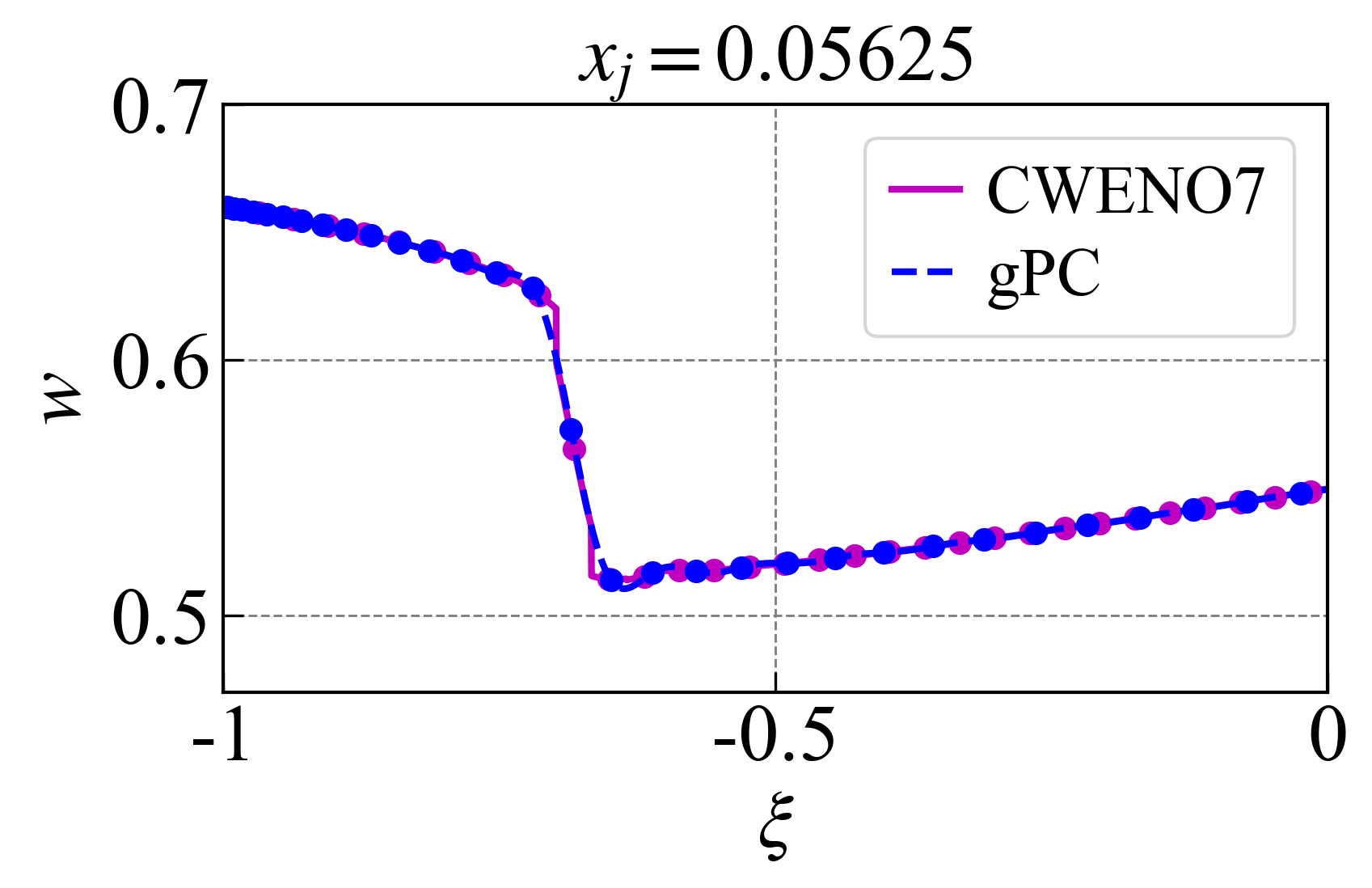}\hspace*{0.3cm}
            \includegraphics[trim=0.2cm 0.3cm 0.2cm 0.2cm, clip, width=0.31\textwidth]{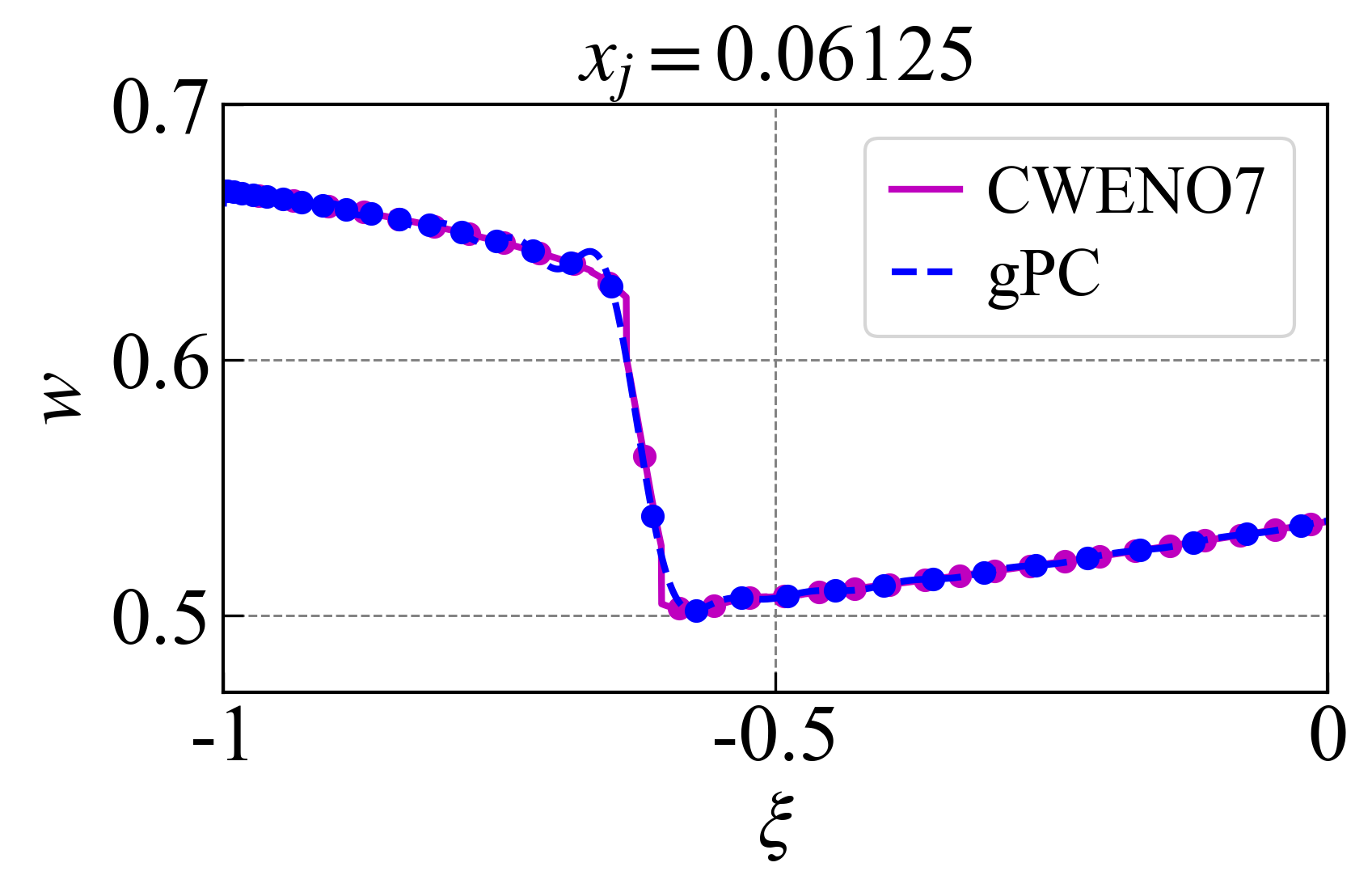}}
\caption{\sf Example 6: 1-D slices of $w(x_j,0.8;\xi)$ (zoomed at $\xi\in[-1,0]$) reconstructed using the gPC and CWENO7
surrogate models at $x_j=0.05125$, $0.05625$, and $0.06125$ for $L=32$ (top row) and $L=64$ (bottom row).\label{fig319}}
\end{figure}

However, the advantage of the proposed CWENO7-based approach can be seen when the water surface and the corresponding PDF are
reconstructed. In the top row of \fref{fig319}, we depict 1-D slices of $w(x_j,0.8;\xi)$, reconstructed using the gPC and CWENO7 surrogate
models at three different cell centers $x_j=0.05125$, $0.05625$, and $0.06125$. As one can see, the CWENO7 interpolation is non-oscillatory,
while the gPC expansion contains oscillations for $L=32$, whose magnitude decay when we take twice bigger number of collocation points
($L=64$); see \fref{fig319} (bottom row).

\begin{figure}[ht!]
\centerline{\includegraphics[trim=0.2cm 0.3cm 0.2cm 0.2cm, clip, width=0.35\textwidth]{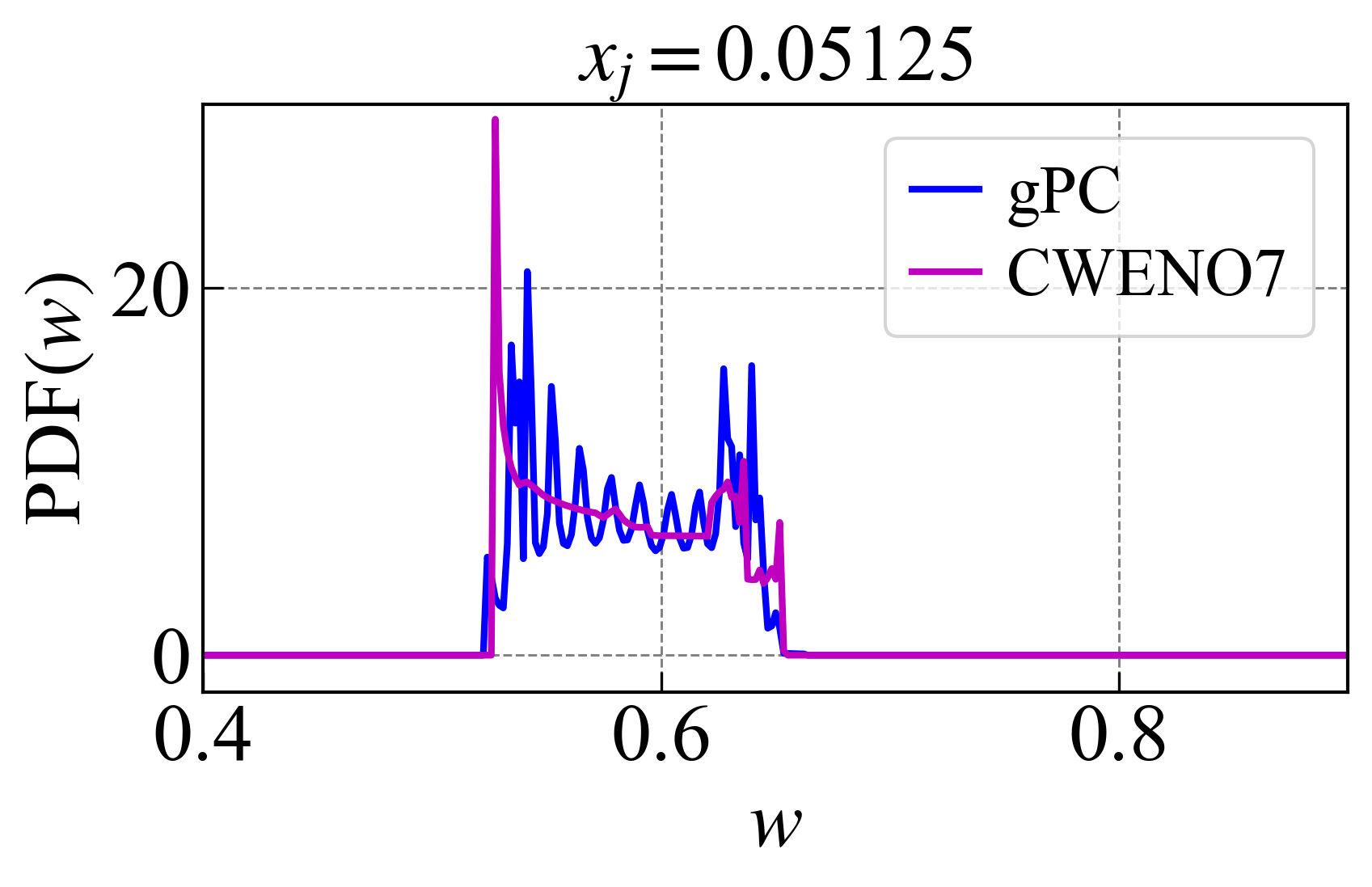}\hspace*{0.5cm}
            \includegraphics[trim=0.2cm 0.3cm 0.2cm 0.2cm, clip, width=0.35\textwidth]{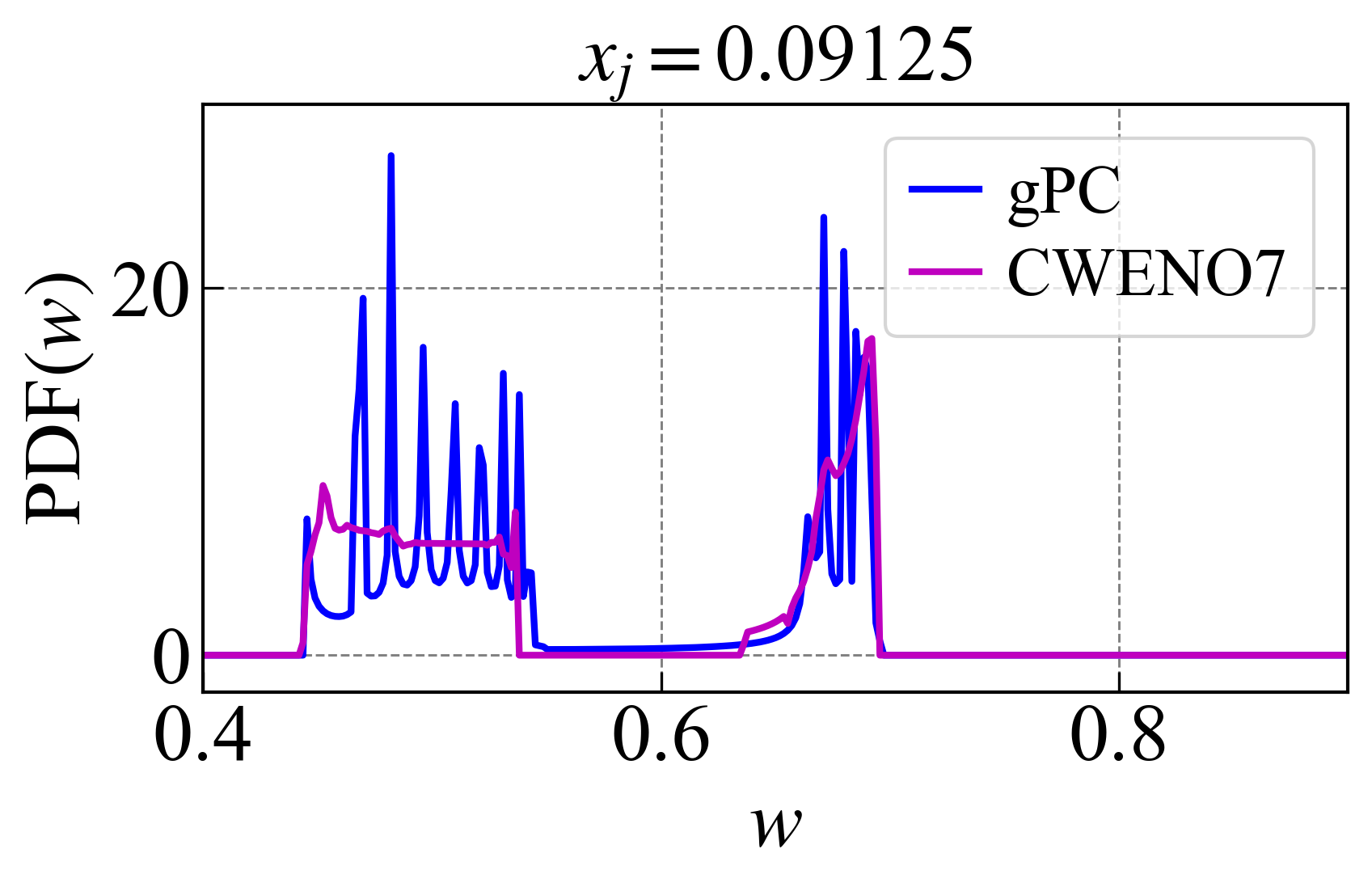}}
\vskip8pt
\centerline{\includegraphics[trim=0.2cm 0.3cm 0.2cm 0.2cm, clip, width=0.35\textwidth]{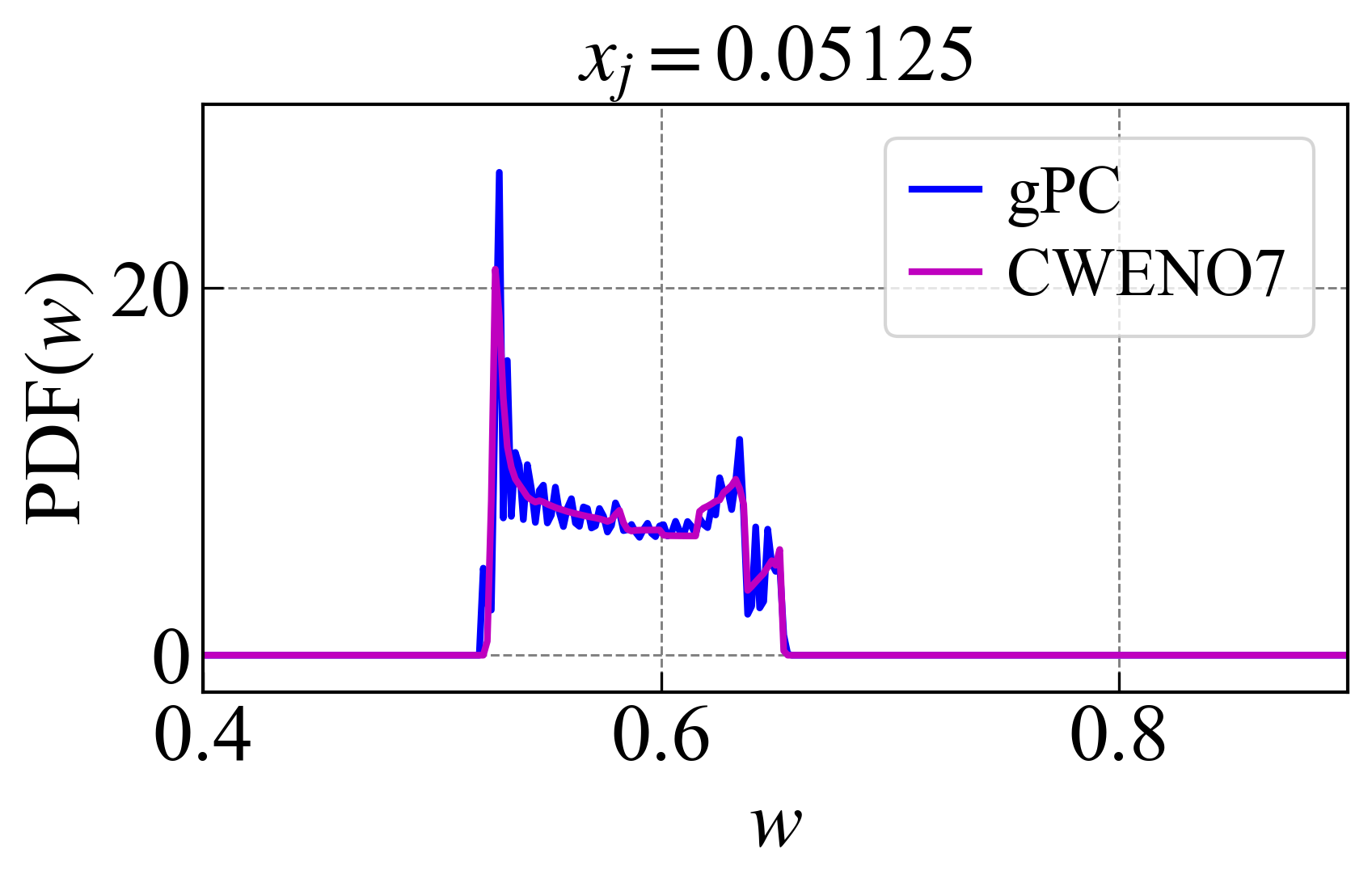}\hspace*{0.5cm}
            \includegraphics[trim=0.2cm 0.3cm 0.2cm 0.2cm, clip, width=0.35\textwidth]{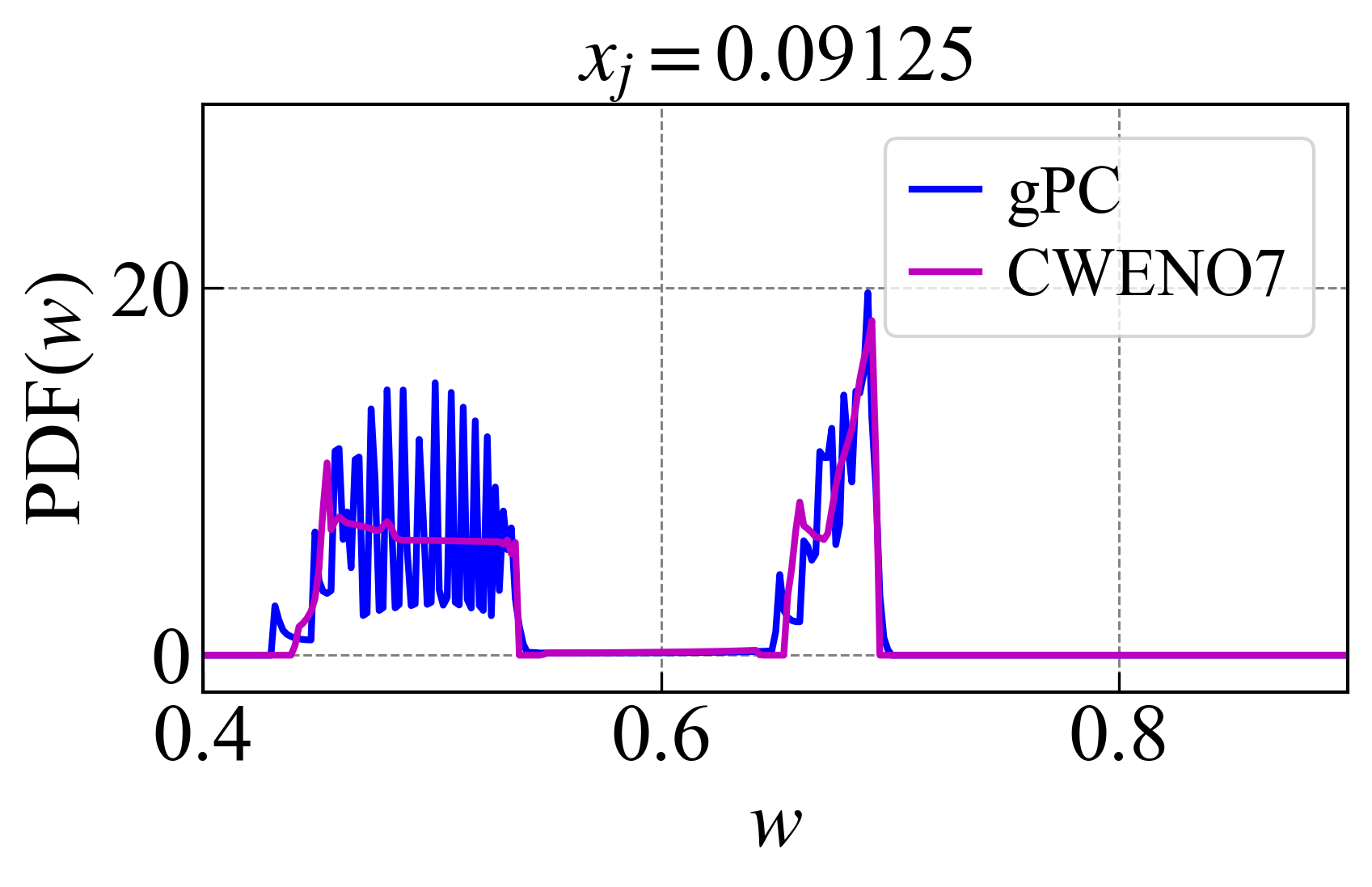}}
\caption{\sf Example 6: PDFs of $w(x,0.8;\xi)$ reconstructed using the gPC- and CWENO7-based surrogate models at $x_j=0.05125$ and
$0.09125$ for $L=32$ (top row) and $L=64$ (bottom row).\label{fig320}}
\end{figure}

In the top row of \fref{fig320}, we show the PDFs of $w(x_j,0.8;\xi)$ (for $x_j=0.05125$ and $x_j=0.09125$) approximated by the surrogate 
models based on the gPC and CWENO7 approaches with $L=32$. As one can see, the PDFs generated by the gPC-based model exhibit oscillatory 
behavior that can misrepresent the underlying distribution, while the PDFs derived from the CWENO7 interpolation are essentially non-
oscillatory. When the number of collocation points is doubled, the resolution achieved using the CWENO7-based model increases, but the PDFs 
approximated using the gPC expansion are still very oscillatory, especially for $x_j=0.09125$. This demonstrates that the proposed
CWENO7-based surrogate model is more accurate and robust than its gPC-based counterpart.

\subsubsection*{Example 7---Dam Break, Two Random Variables}
In the last example, which is a modification of Example 6, we consider a 2-D random variable $\bm\xi=(\xi,\eta)^\top$ with both
$\xi\sim{\cal U}[-1,1]$ and $\eta\sim{\cal U}[-1,1]$ and obtain the discrete function $w(x_j,0.8;\bm\xi)$ by numerically solving \eref{3.6}
subject to the following random initial conditions for the water surface $w$ instead of the deterministic one used in \eref{3.7}:
\begin{equation*}
w(x,0;\xi,\eta)=\left\{\begin{aligned}&1,&&x<0.1\eta,\\&0.5,&&x>0.1\eta,\end{aligned}\right.\qquad u(x,0;\xi,\eta)\equiv0.
\end{equation*}
The rest of the data are precisely the same as in Example 6 and the collocation points are equidistantly placed with $L=M=32$.

\begin{figure}[ht!]
\centerline{\includegraphics[trim=0.2cm 0.3cm 0.2cm 0.2cm, clip, width=0.32\textwidth]{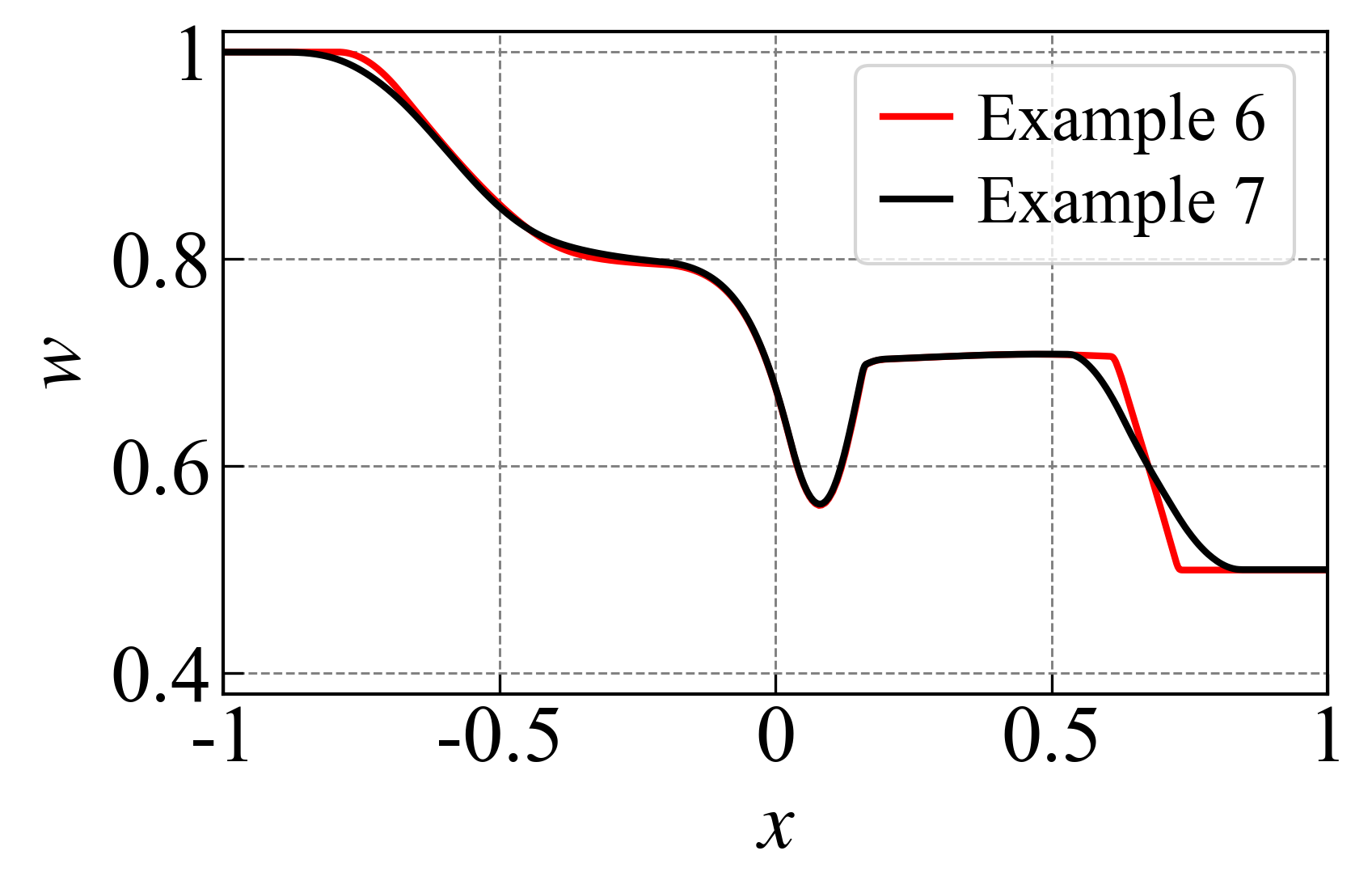}\hspace*{1.5cm}
            \includegraphics[trim=0.2cm 0.3cm 0.2cm 0.2cm, clip, width=0.35\textwidth]{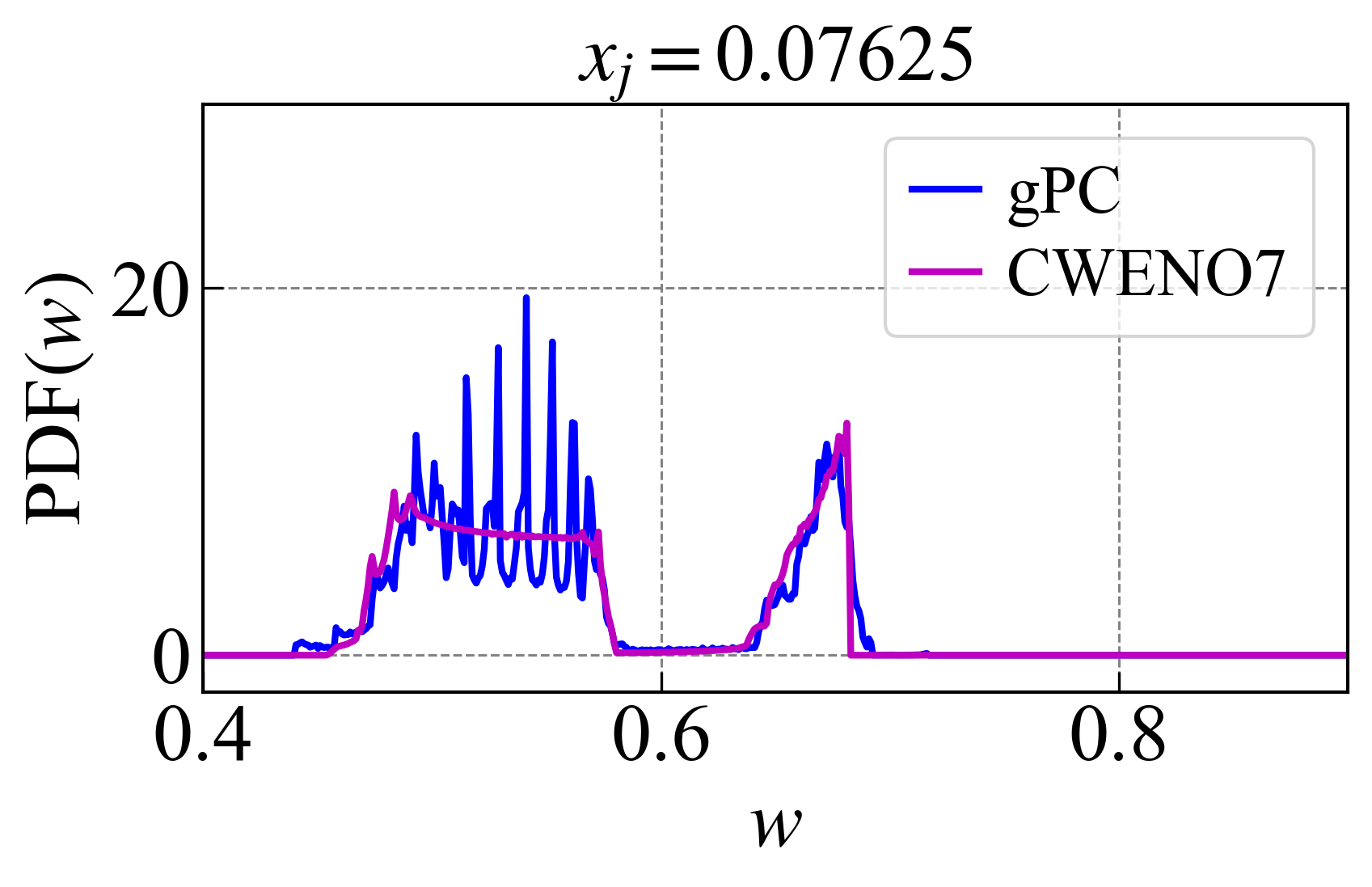}}
\caption{\sf Example 7: Mean of $w$ obtained using the CWENO7 interpolation along with the mean computed by the same surrogate model in
Example 6 (left) and PDFs of $w(x,0.8;\xi)$ reconstructed using the gPC- and CWENO7-based surrogate models at $x_j=0.07625$ (right).
\label{fig321}}
\end{figure}

In \fref{fig321}, we plot the water surface mean computed using the proposed CWENO7-based surrogate model along with the mean from
\fref{fig318} (right). As expected, the mean of $w$ is now substantially more smeared due to the uncertainty in the location of the initial
discontinuity. As in Example 6, the gPC-based surrogate model gives almost identical mean of $w$ and thus we do not show it. We stress that
even though the mean of $w$ is very smeared, the use of the gPC expansion leads to large spurios oscillations in the reconstructed PDFs. To
demonstrate this, we select a particular value $x_j=0.07625$ and plot the PDFs generated by the gPC- and CWENO7-based models.

%%%%%%%%%%%%%
% SECTION 4 %
%%%%%%%%%%%%%
\section{Conclusions}\label{sec4}
We have presented a new surrogate model for forward problems in uncertainty quantification. The proposed model is based on the
seventh-order central weighted essentially non-oscillatory (CWENO7) interpolation in the stochastic collocation framework. We compare the
CWENO7-based surrogate model with the surrogate model based on the generalized polynomial chaos (gPC) expansion on several one- and
two-dimensional test problems, including smooth and nonsmooth functions as well as data generated as solutions of the Saint-Venant system 
of shallow water equations. We have demonstrated that the CWENO7-based model consistently produced accurate results, and its approximation 
of the probability density functions remained stable for discontinuous data, for which the gPC-based model produces very large spurious
oscillations attributed to the Gibbs phenomenon. The proposed CWENO7 surrogate model also works well in higher dimensions with the help of
``dimension-by-dimension'' approach, that is, without requiring special adjustments. Future work will extend the method to its
shape-preserving versions, improve its handling of high-dimensional problems, and investigate higher-order CWENO interpolations.

%%%%%%%%%%%%%%%%%%%
% ACKNOWLEDGEMENT %
%%%%%%%%%%%%%%%%%%%
\begin{acknowledgment}
The work of A. Chertock was supported in part by NSF grant DMS-2208438. The work of A. Kurganov was supported in part by NSFC grants
12171226 and W2431004.

This article is an expanded version of paper \cite{chertock2025} presented at the International Conference on Mathematics and
Computational Methods Applied to Nuclear Science and Engineering (M\&C 2025) on April 29, 2025 in Denver, CO, USA. The conference paper is
available at {\tt https://www.ans.org/pubs/proceedings/article-58202/}.
\end{acknowledgment}

\appendix
\section{Generalized Polynomial Chaos (gPC) Expansion}\label{appA}
Generalized polynomial chaos (gPC) expansion in random variables is an expansion with respect to orthogonal polynomials, which correspond 
to the a priori known random variable probability distributions.

For a single random variable $\xi$ with probability density function (PDF) $p(\xi)$, the gPC expansion of a stochastic process $U(\xi)$ is 
\begin{equation}
U(\xi)\approx\sum_{\ell=0}^L\widehat U_\ell\Phi_\ell(\xi),\quad\widehat U_\ell=\int\limits_\Xi U(\xi)\Phi_\ell(\xi)p(\xi)\,{\rm d}\xi,
\label{A1}
\end{equation}
where $\Phi_{\ell}(\xi)$ are degree $\ell$ polynomials of orthonormal with respect to $p(\xi)$, that is,
$$
\int\limits_\Xi\Phi_\ell(\xi)\Phi_j(\xi)p(\xi)\,{\rm d}\xi=\delta_{\ell j}~~\mbox{for}~j,\ell=0,\dots,L,
$$
$L$ is the maximum polynomial order, and $\delta_{\ell j}$ denotes the Kronecker symbol. Using the coefficients $\widehat U_\ell$, the
approximations of the mean and standard deviation are given by
\begin{equation*}
\widetilde\mu=\widehat U_0,\quad\widetilde\sigma=\bigg(\sum_{\ell=1}^L\widehat U_\ell^2\bigg)^\hf.
\end{equation*}

In practice, the projection integrals in \eref{A1} cannot usually be evaluated analytically, and a numerical quadrature has to be employed 
instead. This leads to the stochastic collocation method, in which the coefficients $\widehat U_\ell$ are approximated from evaluations of
$U(\xi)$ at selected collocation points $\xi=\xi_\ell$ with corresponding quadrature weights. The choice of these points is dictated by the
distribution of the input random variable: for instance, Gauss-Hermite and Gauss-Legendre quadratures are used for Gaussian and uniform 
distributions.  

For problems involving 2-D random variables, the set of collocation points can be formed as tensor products of the univariate collocation
points.
%%%%%%%%%%%%%%%%
% BIBLIOGRAPHY %
%%%%%%%%%%%%%%%%
\bibliographystyle{siamnodash}
\bibliography{ref}

\end{document}